\documentclass[reqno,a4paper,11pt]{amsart}
\usepackage{amsmath,amsthm,amssymb,mathrsfs}
\usepackage[
left = 2.25cm,
right = 2.25cm,
top = 2.2 cm,
bottom = 2.2 cm
]{geometry}
\usepackage{graphicx,xcolor,tikz,mathtools}
\usepackage{bm,stmaryrd}
\usepackage{amsaddr}
\usepackage{enumitem}
\usepackage[toc,page]{appendix}
\usepackage[
pdfencoding=auto,
colorlinks = true,
citecolor = blue,
linkcolor = blue,
anchorcolor = blue,
urlcolor = blue,
filecolor=blue,
pdfkeywords={}
]{hyperref}
\usepackage{esint} % Average integral symbol

\usepackage{todonotes}
\makeatletter
\define@key{todonotes}{Helmut}[]{%
	\setkeys{todonotes}{author=\textbf{Helmut},inline,color=red!10}}%
\define@key{todonotes}{Andrea}[]{%
	\setkeys{todonotes}{author=\textbf{Andrea},inline,color=green!20}}%
\define@key{todonotes}{Yadong}[]{%
	\setkeys{todonotes}{author=\textbf{Yadong},inline,color=cyan!20}}%
\makeatother

\theoremstyle{plain}
\newtheorem{theorem}{Theorem}[section]
\newtheorem{corollary}[theorem]{Corollary}
\newtheorem{definition}[theorem]{Definition}
\newtheorem{lemma}[theorem]{Lemma}
\newtheorem{proposition}[theorem]{Proposition}
\newtheorem{remark}[theorem]{Remark}
\theoremstyle{definition}

\numberwithin{equation}{section}

\newcommand{\bu}{\mathbf{u}}
\newcommand{\bv}{\mathbf{v}}
\newcommand{\bw}{\mathbf{w}}

\newcommand{\bH}{\mathbf{H}}

\newcommand{\R}{\mathbb{R}}

\renewcommand{\d}{\mathrm{d}}
\newcommand{\dx}{\,\d x}
\newcommand{\dt}{\,\d t}

\newcommand{\dtau}{\,\d \tau}
\newcommand{\ddt}{\frac{\d}{\d t}}

\newcommand{\ptial}[1]{ \partial_{#1} }
\newcommand{\pt}{\ptial{t}}
\newcommand{\onehalf}{\frac{1}{2}}

\newcommand{\norm}[1]{\left\Vert #1 \right \Vert}
\newcommand{\normm}[1]{\Vert #1 \Vert}
\newcommand{\normmm}[1]{\left\vert #1 \right\vert}

\newcommand{\dist}{\mathrm{dist}}

\def\Ld{L^2(\Omega)}
\def\bLd{\mathbf L^2(\Omega)}
\def\bLds{\mathbf L^2_\sigma(\Omega)}
\def\N{\mathbb N}
\def\R{\mathbb R}
\def\vphi{\varphi}
\def\Hu{H^1(\Omega)}

\def\bHus{\bH^1_\sigma(\Omega)}
\def\lam{\lambda}
\allowdisplaybreaks[4]

\begin{document}
	
	\title[Conserved Navier--Stokes--Allen--Cahn equations]{Well-posedness and longtime behavior of the conserved Navier--Stokes--Allen--Cahn equations with unmatched viscosities and singular potential}
    % \title[Navier--Stokes Equations on an Elastic Surface]{Local Well-Posedness of the Surface Navier--Stokes Equations on Elastic Membranes}
	
	\author[M. Grasselli, C. Hurm, and A. Poiatti]{
		\small
		Maurizio Grasselli$^{\dagger}$, Christoph Hurm$^\ast$, and
		Andrea Poiatti$^\ddagger$
	}

	\address{
		$^\dagger$Dipartimento di Matematica, Politecnico di Milano, 20133 Milano, Italy
	}
	\email{maurizio.grasselli@polimi.it}

\address{
		$^\ast$Fakult\"at f\"ur Mathematik,
		Universit\"at Regensburg,
		93053 Regensburg, Germany
	}
	\email{christoph.hurm@ur.de}

	\address{
        $^\ddagger$Dipartimento di Scienze Matematiche, Fisiche e Informatiche, \\ Università degli Studi di Parma, 43124 Parma, Italy
	}
	\email{andrea.poiatti@unipr.it}

	\date{\today}
	
	

	\begin{abstract}
     % We consider an incompressible Navier--Stokes system nonlinearly coupled with a conserved Allen--Cahn equation with a singular potential (e.g., of Flory--Huggins type). This is a mass-conserving model for a two-phase flow. Assuming constant density and non-constant viscosity, we first prove, in dimension three, the existence of unique local (in time) strong solutions to an initial-boundary value problem characterized by no-slip boundary condition for the velocity and homogeneous boundary condition for the phase-field. Then we prove, using a relative energy approach, a conditional weak-strong uniqueness in dimension three and the unconditional uniqueness of weak solutions in dimension two. Finally, based on some recent seminal results by the first and third authors, for the first time we show that, in both dimensions and for general singular potentials, any global-in-time weak solution asymptotically separates from pure phases and converges to a unique equilibrium. This result holds under minimal assumptions on the viscosity coefficient. Moreover, if the viscosity coefficient is more regular, then we can resort to the asymptotic strict separation property and the conditional weak-strong uniqueness principle to show that each weak solution asymptotically regularizes. This entails convergence to equilibrium in higher-order norms.
     We consider an incompressible Navier--Stokes system nonlinearly coupled with a conserved Allen--Cahn equation with a singular potential (e.g., of Flory--Huggins type). This model describes a mass-conserving two-phase flow with constant density and non-constant viscosity.

First, in three spatial dimensions, we prove the existence and uniqueness of local-in-time strong solutions to the associated initial--boundary value problem, subject to no-slip boundary conditions for the velocity field and homogeneous boundary conditions for the phase field. Next, by means of a relative energy approach, we establish a conditional weak--strong uniqueness principle in three dimensions, as well as unconditional uniqueness of weak solutions in two dimensions.

Finally, building on recent seminal results by the first and third authors, we prove for the first time that, in both two and three dimensions and for general singular potentials, every global-in-time weak solution asymptotically separates from the pure phases and converges to a unique equilibrium. This result is obtained under minimal assumptions on the viscosity coefficient. Moreover, under additional regularity assumptions on the viscosity, we combine the asymptotic strict separation property with the conditional weak--strong uniqueness principle to show that weak solutions undergo asymptotic regularization. As a consequence, convergence to equilibrium also holds in higher-order norms.
	\end{abstract}
	
	\maketitle
	
	% 	\setcounter{tocdepth}{1}
	% 	\tableofcontents
	
	\section{Introduction}
    In this contribution we aim to study a mass-conserving model of two-phase incompressible flow with matched densities. Suppose that our binary fluid occupies a bounded open set $\Omega\subset \mathbb{R}^d$, $d\in\{2,3\}$, in a given time interval $(0,T)$, $T>0$. The model consists of the Navier--Stokes system coupled with a conserved Allen--Cahn equation endowed with suitable boundary and initial conditions.
    More precisely, denoting by $\bv$ the (volume) averaged velocity of the fluid mixture and by $\varphi$ the relative concentration of one component, the model reads as follows
    \begin{align}
        \partial_t\bv + (\bv\cdot\nabla)\bv - \text{div}(\nu(\varphi)D\bv) + \nabla \pi &= \mu\nabla\varphi &\quad&\text{ in }\Omega\times (0,T), \label{Intro:Eq:1}\\
        \nabla\cdot\bv &= 0 &\quad&\text{ in }\Omega\times (0,T), \label{Intro:Eq:2}\\
        \partial_t\varphi + \bv\cdot\nabla\varphi + \mu - \overline{\mu} &= 0 &\quad&\text{ in }\Omega \times (0,T), \label{Intro:Eq:3}\\
        \mu &= - \Delta\varphi + f^\prime(\varphi) &\quad&\text{ in }\Omega \times (0,T), \label{Intro:Eq:4}\\
        \bv = \partial_\mathbf{n}\varphi &= 0 &\quad&\text{ on }\partial\Omega \times (0,T), \label{Intro:Eq:5}\\
        \bv(0) = \bv_0,\quad \varphi(0) &= \varphi_0 &\quad&\text{ in }\Omega. \label{Intro:Eq:6}
    \end{align}
    Here the density and the mobility are equal to the unity, $\nu$ is the viscosity of the mixture, $\overline{w}$ stands for the integral average of an integrable
    function $w$ over $\Omega$, and $\mathbf{n}$ stands for the outward unit normal to $\partial\Omega$. Moreover, $f$ is a general singular potential, typically the well-known Flory--Huggins potential density
    \begin{equation}
    \label{FHP}
    f(s)=\frac{\theta}{2}((1+s)\text{ln}(1+s)+(1-s)\text{ln}(1-s)) - \frac{\theta_0}{2}\vphi^2
    \quad\forall\,s\in(-1,1).
    \end{equation}
    Here $\theta>0$ stands for the absolute temperature, while
    $\theta_0>0$ is related to the critical temperature below which $f$ has a double-well shape.

    Problem \eqref{Intro:Eq:1}-\eqref{Intro:Eq:6} was analyzed in \cite{GGW2022} in the case of nonconstant density (for details about the model, see its Introduction and the references therein, in particular,
    \cite{GKL2018}).
    There, the existence of a global weak solution in dimension three as well the existence of a strong solution in dimension two are established. Concerning uniqueness, a two-dimensional conditional weak-strong uniqueness holds.
    
    Here we prove the local well-posedness of strong solutions in dimension three as well as a conditional weak-strong uniqueness. These solutions are instantaneously uniformly strictly separated from pure phases. Let us point out that the condition to ensure weak-strong uniqueness, which is shown by means of a relative energy approach, is weaker than the more standard one of having two solutions uniformly strictly separated from pure phases (see, e.g., \cite{GGGP}). Indeed, exploiting a technique first developed in \cite{HKP}, we allow the weak solution to be not necessarily separated, but such that the quantity $F'(\varphi)$ is integrable in space-time with a sufficiently large  exponent. On the other hand, in dimension two, we show that weak solutions are
    unconditionally unique, which was left open in \cite{GGW2022}. Moreover, we establish that any global weak solution asymptotically separates from pure phases and converges to a unique equilibrium in dimension three (and, of course, two), without the necessity of further regularization. This allows to treat also the case of minimal regularity on the viscosity $\nu$, namely assuming $\nu\in C^0([-1,1])$. The fundamental argument is based on a new technique introduced in \cite{GP2} (see also the seminal
    \cite{GP} where this approach was first developed, and \cite{GLPW} for a further application to a coupled system).  Additionally, in case $\nu$ is Lipschitz continuous, thanks to the validity of the asymptotic strict separation property, together with the weak-strong uniqueness principle, we are even able to show that any weak solution asymptotically regularizes and thus converges to the unique equilibrium in higher-order norms.
    
    We recall that the existence of a global weak solution to system \eqref{Intro:Eq:1}-\eqref{Intro:Eq:4} endowed with initial and dynamic boundary conditions has been established in \cite{GGPC}, in the case of constant viscosity (see also \cite{LYY} for a smooth approximation of \eqref{FHP}). More recently, the sharp interface limit has been investigated in \cite{AM} by means of rigorous matched asymptotic expansions. The existence of a global attractor for the original model studied in \cite{GGW2022} has been proven in \cite{ST25}. A stochastic version of \eqref{Intro:Eq:1}-\eqref{Intro:Eq:6} with constant viscosity has been analyzed in \cite{DGS}. In conclusion, multi-component extensions of similar models have also been analyzed for instance in \cite{AGGmulti, GPCAC}, whereas a non-isothermal variant to problem \eqref{Intro:Eq:1}-\eqref{Intro:Eq:6}, together with its sharp interface limit, was studied in \cite{AMP}.

The plan of the paper goes as follows. In the next section we introduce notation and specify the basic assumptions. The main theorems are stated in Section \ref{MR}. The existence of a global weak solution and the local existence (and uniqueness) of global solution are proven in Sections \ref{weak} and \ref{strong}, respectively. Section \ref{weak-strong} is devoted to establish the conditional weak-strong uniqueness. Uniqueness of weak solutions in dimension two is demonstrated in Section \ref{weakuniq}. The last sections are essentially devoted to the longtime behavior of weak solutions. More precisely, preliminary results are proven in Sections \ref{proofequil} and \ref{goodeq}, while the proof of convergence to a single equilibrium is given in Section \ref{uniqeq}. Finally, the asymptotic regularization of a weak solution and its convergence to equilibrium in higher-order norms are shown in Section \ref{asymreg}.

    \section{Notation and basic assumptions}

  Here we first introduce some notation along with the functional spaces which will be used in the sequel.

 \begin{enumerate}[label=\textnormal{(N\arabic*)},leftmargin=*]
 	
 	\item \textbf{Notation for general Banach spaces.}
 	For any (real) normed space $X$%of scalar-valued functions
    , we denote its norm by $\|\cdot\|_X$,
 	its {dual space by $X'$}.
 	If $X$ is a Hilbert space, we write $(\cdot,\cdot)_X$ to denote the corresponding inner product. However, in the case $X=\mathbb{R}^N$,
    we use the dot.
 	Moreover, the corresponding spaces of vector-valued or matrix-valued functions with each component in $X$ are denoted by $\mathbf{X}$.
 	
 	\item \textbf{Lebesgue and Sobolev spaces.}
 	For $1 \leq p \leq \infty$ and $k \in \N$, the classical Lebesgue and Sobolev spaces defined on $\Omega$ are denoted by $L^p(\Omega)$ and $W^{k,p}(\Omega)$, and their standard norms are denoted by $\|\cdot\|_{L^p(\Omega)}$ and $\|\cdot\|_{W^{k,p}(\Omega)}$, respectively.
 	In the case $p = 2$, we set $H^k(\Omega) = W^{k,2}(\Omega)$ and, in particular, we define
 $$
 H_N^2(\Omega) := \{v \in H^2(\Omega): \; \,\partial_{\mathbf{n}}v=0 \; \text{ a.e. on } \partial\Omega\}.
 $$
  Note that the $L^2(\Omega)$ inner product is simply denoted by $(\cdot,\cdot)$ even in the case $\mathbf{L}^2(\Omega)$.
 	Also, for any interval $I\subset\R$, any Banach space $X$, $1 \leq p \leq \infty$ and $k \in \N$, we write $L^p(I;X)$, $W^{k,p}(I;X)$ and $H^{k}(I;X) = W^{k,2}(I;X)$ to denote the Lebesgue and Sobolev spaces of functions with values in $X$. The canonical norms are indicated by $\|\cdot\|_{L^p(I;X)}$, $\|\cdot\|_{W^{k,p}(I;X)}$ and $\|\cdot\|_{H^k(I;X)}$, respectively.
 	    We additionally define
\begin{align*}
    L^p_\mathrm{loc}(I;X)
    &:=
    \big\{
        u:I\to X \,\big\vert\, u \in L^p(J;X) \;\text{for every compact interval $J\subset I$}
    \big\}
    \\[1ex]
    L^p_\mathrm{uloc}(I;X)
    &:=
    \left\{ u:I\to X \,\middle|\,
    \begin{aligned}
    &u \in L^p_\mathrm{loc}(I;X) \;\text{and}\; \exists\, C>0\; \sup_{t\in I}\|u\|_{L^p(t,t+1;X)} \le C
    \end{aligned}
    \right\}.
\end{align*}
The spaces $W^{k,p}_\mathrm{loc}(I;X)$, $H^k_\mathrm{loc}(I;X)$, $W^{k,p}_\mathrm{uloc}(I;X)$, $H^k_\mathrm{uloc}(I;X)$ are defined analogously.
 %	For simplicity, we just write $(\cdot,\cdot) := (\cdot,\cdot)_{L^2(\Omega)}$, $\|\cdot\|:=\|\cdot\|_{L^2(\Omega)}$.
 	\item \textbf{Spaces of continuous functions.}
 	For any interval $I\subset\R$ and any Banach space $X$, $C(I;X)$ denotes the space of continuous functions mapping from $I$ to $X$ and $BC(I;X)$ denotes the space of bounded functions in $C(I;X)$. Moreover, $C_\mathrm{w}(I;X)$ denotes the space of functions mapping from $I$ to $X$, which are continuous on $I$ with respect to the weak topology of $X$, and $BC_\mathrm{w}(I;X)$ denotes the space of bounded functions in $C_\mathrm{w}(I;X)$. Then, we denote by $C^\gamma(I;X)$, $\gamma\in(0,1]$, the space of $\gamma$-H\"{o}lder (Lipschitz, if $\gamma=1$) continuous functions with values in $X$. In addition, $C_0^k(I;X)$ stands for the space of $k$-continuously differentiable functions with compact support mapping $I$ into $X$.

 	\item \textbf{Spaces of divergence-free functions.}
 We define the closed linear subspaces
 	\begin{align*}
 		\mathbf L^p_\sigma(\Omega)
 		&:=\overline{\{\mathbf{u}\in \mathbf{C}^\infty_0(\Omega) \,\big\vert\, \operatorname{div}\ \mathbf{u}=0\}}^{\mathbf{L}^p(\Omega)}
 		\subset \mathbf L^p(\Omega), \quad p\in[2,\infty),\\
 		\mathbf H^1_\sigma(\Omega)&:= \mathbf L^2_\sigma(\Omega) \cap \mathbf H^1(\Omega).
 	\end{align*}
In both cases, Korn's inequality yields
\begin{equation}
\Vert \mathbf{u}\Vert \leq \sqrt{2}\Vert D\mathbf{u}%
\Vert\leq \sqrt{2}\Vert \nabla \mathbf{u}\Vert
\quad \text{ for all } \mathbf{u}\in \mathbf H^1_\sigma (\Omega).
\label{korn}
\end{equation}
As a trivial consequence, $\|\nabla\cdot\|$ is a norm on $\mathbf H^1_\sigma(\Omega)$ that is equivalent to the standard norm $\|\cdot\|_{\mathbf{H}^1(\Omega)}$.
 \end{enumerate}

We now state the general assumptions that are supposed to hold throughout this paper unless otherwise specified.

\begin{enumerate}[label=\textnormal{(A\arabic*)},leftmargin=*]
    \item \label{ASS:1} Let $\Omega \subset \R^3$ be a bounded domain with $C^{2}$-boundary.
    \item \label{ASS:Viscosity} The viscosity $\nu\in W^{1,\infty}(\R)$ satisfies
    \begin{align*}
        0 < \nu_* \leq \nu(s) \leq \nu^*\qquad \text{ for all }s\in\R,
    \end{align*}
    for some positive constants $\nu_*,\nu^*\in \R$.
    %\item \label{ASS:3} The density $\rho$ and the mobility $m$ are positive constants. For convenience, we set $\rho=m=1$.
    \item \label{ASS:S1}
    The potential $f:[-1,1]\to \R$ exhibits the decomposition
    \begin{equation*}
        f(s)=F(s)-\frac{\theta_0}{2}s^2 \quad\text{for all $s\in [-1,1]$}
    \end{equation*}
    with a given constant $\theta_0>0$.
    Here, $F\in C([-1,1])\cap C^{2}(-1,1)$ has the properties
    \begin{equation*}
    \lim_{r\rightarrow -1}F^{\prime }(r)=-\infty ,
    \quad \lim_{r\rightarrow 1}F^{\prime }(r)=\infty ,
    \quad F^{\prime \prime }(s)\geq {\theta},
    \quad F'(0)=0
    \end{equation*}
    for all $s\in (-1,1)$ and a prescribed constant $\theta\in(0,\theta_0)$.
    Without loss of generality, we further assume $F(0)=0$.
    In particular, this means that $F(s)\geq 0$ for all $s\in [-1,1]$.

    For convenience, we extend $f$ and $F$ onto $\R\setminus[-1,1]$ by defining
    $f(s):=\infty $ and $F(s):=\infty $ for all $s\in\R\setminus [-1,1]$.
    %%%%%
    % \item \label{ASS:S2} In addition to \ref{ASS:S1}, there exists $\beta>\frac12$ such that
    % \begin{equation}
    % \frac{1}{F^{\prime }(1-2\delta )}=O\left( \frac{1}{|\ln (\delta )|^{\beta }}%
    % \right) ,\quad\text{ }\dfrac{1}{|F^{\prime }(-1+2\delta )|}=O\left( \frac{1}{%
    % |\ln (\delta )|^{\beta }}\right) .  \label{est}
    % \end{equation}
    % as $\delta\to 0^+$.
    %%%%%
    % \item\label{ASS:S3} In addition to \ref{ASS:S1}, it holds
    % \begin{alignat}{2}
    % 	\frac{1}{F^{\prime}(1-2\delta)}&=O\left(\frac{1}{\vert\ln(\delta)\vert}\right),
    %     &\quad\frac{1}{F^{\prime\prime}(1-2\delta)}&=O(\delta),
    % 	\label{F}
    %     \\
    % 	\dfrac{1}{\vert F^{\prime}(-1+2\delta)\vert }&=O\left(\frac{1}{\vert\ln(\delta)\vert}\right),
    %     &\quad\dfrac{1}{F^{\prime\prime}(-1+2\delta)}&=O\left(\delta\right).
    % 	\label{F2}
    % \end{alignat}
    % as $\delta\to 0^+$.
    % Moreover, there exists $\gamma_{0}>0$ such that $F^{\prime \prime }$ is monotonously increasing on $(-1,-1+\gamma_0]$ and on $[1-\gamma _{0},1)$.
\end{enumerate}
    \section{Main results}\label{MR}
    First of all, we state the existence of a weak solution. This result can be deduced from the proof of \cite[Thm.~3.1]{GGW2022} by assuming constant density.
    However, we report here an alternative proof (i.e. a full Galerkin scheme in place of a
    semi-Galerkin one) since it will be useful to prove the local existence of a strong solution (see Theorem \ref{THM:StrongWP} below).

    \begin{theorem}[Existence of global weak solutions]\label{THM:WeakExistence}
        Let the assumptions \ref{ASS:1}--\ref{ASS:S1} hold. Assume that $\bv_0\in ~\mathbf{L}^2_\sigma(\Omega)$ and $\varphi_0\in H^1(\Omega)$ with $\|\varphi_0\|_{L^\infty(\Omega)} \leq 1$ and $|\overline{\varphi_0}| < 1$. Then there exists a weak solution $(\bv,\varphi,\mu)$ to \eqref{Intro:Eq:1}--\eqref{Intro:Eq:6}, globally defined in $(0,\infty)$,  such that, for any $T>0$,
    \begin{enumerate}[label=\textnormal{(\roman*)}, topsep=0em, partopsep=0em, leftmargin=*]
            \item It holds
            \begin{align}
            \begin{cases}
                \bv \in L^\infty(0,T;\mathbf{L}^2_\sigma(\Omega))\cap L^2(0,T;\mathbf{H}^1_\sigma(\Omega)), \\
                \varphi \in L^\infty(0,T;H^1(\Omega))\cap L^2(0,T;H^2_N(\Omega))\text{ with }|\varphi |<1\ \text{a.e. in }\Omega_T, \\
                \partial_t\bv \in L^{4/3}(0,T;\mathbf{H}^1_\sigma(\Omega)^\prime), \\
                \partial_t\varphi \in L^2(0,T;H^1(\Omega)^\prime), \\
                \mu \in L^2(0,T;L^2(\Omega)).
            \end{cases}
            \end{align}
            \item The solution $(\bv,\varphi,\mu)$ fulfills:
            \begin{align}
                \label{w1}&\langle\partial_t\bv,\bw\rangle_{\mathbf{H}^1_\sigma(\Omega)} + ((\bv\cdot\nabla)\bv,\bw) + (\nu(\varphi)D\bv,D\bw) = (\mu\nabla\varphi,\bw), \\
                \label{w2}&\langle\partial_t\varphi,v\rangle_{H^1(\Omega)} + (\bv\cdot\nabla\varphi,v) + (\mu-\overline{\mu},v) = 0
            \end{align}
            for all $\bw\in \mathbf{H}^1_\sigma(\Omega)$, $v\in H^1(\Omega)$ and for almost all $t\in(0,T)$.
            \item The initial conditions $\bv(\cdot,0) = \bv_0$ and $\varphi(\cdot,0)=\varphi_0$ hold almost everywhere in $\Omega$.
            \item The energy inequality
            \begin{align}
                E(\bv(t),\varphi(t)) + \int_s^t\left[\|\sqrt{\nu(\varphi(\tau))}D\bv(\tau)\|_{\mathbf{L}^2(\Omega)}^2 + \|\mu(\tau)-\overline{\mu(\tau)}\|_{L^2(\Omega)}^2\right]\d\tau \leq E(\bv(s),\varphi(s))\label{ineq1}
            \end{align}
            holds for all $t\geq 0$ and almost any $s\geq0$, including $s=0$, where
     \begin{align}
                E(\bv,\varphi) := \frac{1}{2}\int_\Omega|\bv|^2\dx + \underbrace{\frac{1}{2}\int_\Omega|\nabla\varphi|^2\dx + \int_\Omega f(\varphi)\dx}_{:=E_{free}(\vphi)}.
            \end{align}
        \end{enumerate}
    \end{theorem}

 \begin{remark}
        	Notice that $\mathbf v\in C_{w}([0,T],\mathbf H^1_\sigma(\Omega))$, and $\varphi\in C([0,T];H^{\frac12}(\Omega))$, so that the initial conditions make sense, where the limit as $t\to0$ is intended in $\mathbf H^1_\sigma(\Omega)$-weak and in $H^\frac12(\Omega)$, respectively.
        Still about regularity, differently from the case of Navier-Stokes-Cahn-Hilliard systems (cf. \cite{GP}), here we cannot consider a split energy inequality as \cite[Eqs. (2.10)-(2.11)]{GP}. Indeed, the kinetic energy inequality would account for the term $(\bv\cdot\nabla \vphi,\mu),$ which does not belong to $L^1(0,T)$ given the regularity of a weak solution as in Theorem \ref{THM:WeakExistence}.
        \label{nosplit}\end{remark}

\begin{remark}
Going to the proof of Theorem \ref{THM:WeakExistence}, it is easy to realize that assumption \ref{ASS:Viscosity} can be weakened, that is, we can just assume that $\nu\in C^0([-1,1])$ and bounded from
below by a positive constant.
\label{visco}
\end{remark}

\begin{remark}
        As usual, the pressure term  is dropped in the weak formulation \eqref{w1}-\eqref{w2}.
The pressure $\pi$ can be recovered (up to a constant) thanks to the classical de Rham’s
theorem (see, for instance, \cite[Sec.V.1.5]{BoyerFabrie}).\label{pressure}
    \end{remark}

\begin{remark}
We recall that, in the case $\Omega\subset \R^2$, on account of the two-dimensional Lady\v{z}enskaja inequality, one can easily prove that (see, e.g., \cite[Theorem 3.2]{GGW2022}) $\partial _t\varphi\in L^2(0,T;L^2(\Omega))$ and $\partial_t\bv\in L^2(0,T;\bH_\sigma^1(\Omega)')$.
    \label{2d}
\end{remark}

We can now state our first main result.

    \begin{theorem}[Existence and uniqueness of local strong solutions]\label{THM:StrongWP}
        Let the assumptions \ref{ASS:1}--\ref{ASS:S1} hold. Assume that  $\bv_0\in ~\mathbf{H}^1_\sigma(\Omega)$ and $\varphi_0\in H^2_N(\Omega)$, $F^\prime(\varphi_0)\in L^2(\Omega)$, $\|\varphi_0\|_{L^\infty(\Omega)} \leq 1$ and $|\overline{\varphi_0}| < 1$. Then, there exists a strong solution $(\bv,\varphi,\mu)$ to \eqref{Intro:Eq:1}--\eqref{Intro:Eq:6} in a right-maximal time interval $[0,T_*)$ such that
        \begin{enumerate}[label=\textnormal{(\roman*)}, topsep=0em, partopsep=0em, leftmargin=*]
            \item It holds
            \begin{align}
                \begin{cases}
                    \bv \in L^\infty(0,T_*;\mathbf{H}^1_\sigma(\Omega))\cap L^2(0,T_*;\mathbf{H}^2(\Omega))\cap H^1(0,T_*;\mathbf{L}^2_\sigma(\Omega)),\\
                    \varphi\in L^\infty(0,T_*;H^2_N(\Omega))\cap H^1(0,T_*;H^1(\Omega)),\\
                    |\varphi(x,t)|<1\text{ for a.a. }x\in\Omega\text{ and all }t\in[0,T_*),\\
                    \mu\in L^\infty(0,T_*;L^2(\Omega))\cap L^2(0,T_*;H^1(\Omega)).
                \end{cases}
            \end{align}
            \item The triplet $(\bv,\varphi,\mu)$ fulfills equations \eqref{Intro:Eq:1}--\eqref{Intro:Eq:4} almost everywhere in $\Omega\times[0,T_*)$ and satisfies the boundary conditions \eqref{Intro:Eq:5} almost everywhere in $\partial\Omega\times(0,T_*)$.
            \item The initial conditions $\bv(\cdot,0) = \bv_0$ and $\varphi(\cdot,0)=\varphi_0$ hold in $\Omega$.
            \item For any $\tau>0$ the exists $\delta=\delta(\tau)\in(0,1)$ such that
\begin{align}
    \sup_{t\in(\tau,T_\star)}\norm{\varphi(t)}_{L^\infty(\Omega)}\leq 1-\delta.\label{sepa1}
\end{align}
If, additionally, $\varphi_0$ is strictly separated from pure phases, i.e., $  \exists \delta_0\in(0,1): \norm{\varphi_0}_{L^\infty(\Omega)}\leq 1-\delta_0$, then there exists $\delta\in(0,1)$, depending also on $\delta_0$, such that
\begin{align}
    \label{sepa2}\sup_{t\in[0,T_\star)}\norm{\varphi(t)}_{L^\infty(\Omega)}\leq 1-\delta.
\end{align}

            \item The solution $(\bv,\varphi)$ is unique and depends continuously on the initial data. More precisely, if $(\bv_i,\varphi_i,\mu_i)$ is the strong solution to \eqref{Intro:Eq:1}--\eqref{Intro:Eq:6} originating from the initial datum $(\bv_{0,i},\varphi_{0,i})$, $i=1,2$, then
            \begin{align}
            \begin{split}
                &\|\bv_1(t) - \bv_2(t)\|_{\mathbf{H}^1_\sigma(\Omega)^\prime}^2 + \|\varphi_1(t) - \varphi_2(t)\|_{L^2(\Omega)}^2 \\
                &\qquad \leq C\Big(\|\bv_{0,1} - \bv_{0,2}\|_{\mathbf{H}^1_\sigma(\Omega)^\prime}^2 + \|\varphi_{0,1} - \varphi_{0,2}\|_{L^2(\Omega)}^2\Big)e^{CT_*},
            \end{split}\label{contdep}
            \end{align}
           for some $C>0$ and for any $t\in(0,T_\star)$.
        \end{enumerate}
    \end{theorem}
    \begin{remark}
        The energy inequality \eqref{ineq1} is actually an identity for the strong solution.
    \end{remark}
    \begin{remark}
        For the strong solution, the pressure $\pi$ (see Remark \ref{pressure}) is more regular, namely, we have $\pi\in L^2(0,T_*;\Hu)$.
    \end{remark}
    \begin{remark}
    Adapting a technique devised in \cite{WW12} and then used in many subsequent papers, one can prove the existence of a global strong solution if the initial data are sufficiently close to an equilibrium (see \eqref{conv1t} below). This argument is based on a careful use of the {\L}ojasiewicz--Simon inequality (see Proposition \ref{LSineq}).
    \end{remark}
    The following result entails a conditional weak-strong uniqueness (cf. \cite{HM19} where the viscosity is constant, the potential is smooth and the Allen--Cahn equation is not the conserved one). We prove it by means of a relative energy method (see \cite{FPP2019} and references therein).

    \begin{theorem}[Conditional weak-strong continuous dependence]\label{THM:WSU}
        Let the assumptions \ref{ASS:1}--\ref{ASS:S1} hold. Let $(\bv_{0,1},\varphi_{0,1})$ be an initial datum as in Theorem \ref{THM:WeakExistence}. Then, we consider a weak solution $(\bv_1,\varphi_1,\mu_1)$ to \eqref{Intro:Eq:1}--\eqref{Intro:Eq:6} originating from $(\bv_{0,1},\varphi_{0,1})$ and  satisfying the additional regularity
        \begin{align}
        F'(\varphi_1)\in L^4(0,T_*;L^4(\Omega)).\label{additional}
        \end{align}
      Let $(\bv_{0,2},\varphi_{0,2})$ an initial datum as in Theorem \ref{THM:StrongWP}, such that $\overline{\vphi_{0,2}}=\overline{\vphi_{0,1}}$, and  $\varphi_{0,2}$ is strictly separated from pure phases, i.e.,   there exists $\delta_0\in(0,1): \norm{\varphi_{0,2}}_{L^\infty(\Omega)}\leq 1-\delta_0$. Then, we consider the strong solution $(\bv_2,\varphi_2,\mu_2)$ to \eqref{Intro:Eq:1}--\eqref{Intro:Eq:6} corresponding to $(\bv_{0,2},\varphi_{0,2})$. The following continuous dependence estimate holds
        \begin{align}
        \begin{split}
            &\|\bv_1(t)-\bv_2(t)\|_{\mathbf{L}^2_\sigma(\Omega)}^2 + \|\nabla(\varphi_1(t)-\varphi_2(t))\|_{\mathbf{L}^2(\Omega)}^2 \\
            &\qquad \leq C\Big(\|\bv_{0,1}-\bv_{0,2}\|_{\mathbf{L}^2_\sigma(\Omega)}^2 + \|\nabla(\varphi_{0,1}-\varphi_{0,2})\|_{\mathbf{L}^2(\Omega)}^2\Big)e^{\int_0^{T_*}\Lambda(t)\dt},
        \end{split}
        \label{contdep2}
        \end{align}
        for all $t\in(0,T_*)$, where $T_*>0$ is given as in Theorem \ref{THM:StrongWP} and $\Lambda\in L^1(0,T_*)$ is defined by
        \begin{align*}
            \Lambda(t) := C\Big(1 + \|\bv_2(t)\|_{\mathbf{H}^2(\Omega)}+\|\mu_2(t)\|_{H^1(\Omega)}^2+\norm{\varphi_1(t)}_{H^2(\Omega)}^2+ \|f^\prime(\varphi_2(t))\|_{L^4(\Omega)}^4 + \|f^\prime(\varphi_1(t))\|_{L^4(\Omega)}^4\Big),
        \end{align*}
        for almost any $t\in(0,T_*)$.
    \end{theorem}
    \begin{remark}
        The additional assumption \eqref{additional} on the weak solution is necessary to show that $\Lambda\in L^1(0,T_*)$. For instance, if the weak solution is strictly separated from pure phases, say
\begin{align}
   \sup_{t\in[0,T_*)} \norm{\varphi_1(t)}_{L^\infty(\Omega)}\leq 1-\delta,\label{str}
\end{align}
for some $\delta\in(0,1)$, then condition \eqref{additional} holds. On the other hand, notice that \eqref{additional} is actually much weaker than the strict separation assumption \eqref{str}.
    \end{remark}
The following result shows that weak solutions are unconditionally unique in dimension two, solving an open question left in \cite[Theorem 3.3]{GGW2022}, where the uniqueness of weak solutions is only conditional.

    \begin{theorem}[Continuous dependence on initial data of weak solutions in $2$D]\label{THM:Unique:2D}
    Let $\Omega\subset \R^2$ be a bounded domain with $C^2$-boundary, and let the assumptions \ref{ASS:1}--\ref{ASS:S1} hold. If $(\bv_{0,i},\varphi_{0,i})$, $i=1,2$, satisfy the assumptions of Theorem \ref{THM:WeakExistence}, with $\overline{\vphi_{0,1}}=\overline{\vphi_{0,2}}$, then, given two weak solutions $(\bv_i,\varphi_i)$ to \eqref{Intro:Eq:1}--\eqref{Intro:Eq:6} originating from $(\bv_{0,i},\varphi_{0,i})$, the following continuous dependence estimate holds for any fixed $T>0$
\begin{align}
   &\nonumber \sup_{t\in[0,T]}\left(\norm{\bv_1(t)-\bv_2(t)}_{\bHus'}+\norm{\varphi_1(t)-\vphi_2(t)}_{L^2(\Omega)}\right)\\&
   \leq C(T)(\norm{\bv_{0,1}-\bv_{0,2}}_{\bHus'}+\norm{\varphi_{0,1}-\vphi_{0,2}}_{L^2(\Omega)})e^{\int_{0}^T\Lambda_1(t)\dt}\label{estimate1},
\end{align}
where $C=C(T)>0$ depends on the structural constants of the system and $\Lambda_1\in L^1(0,T)$ is defined by
\begin{align*}
        \Lambda_1(t) := C\Big(\|\bv_1(t)\|_{\mathbf{L}^4(\Omega)}^4 + \|\bv_2(t)\|_{\mathbf{L}^4(\Omega)}^4 + \|D\bv_1(t)\|_{\mathbf{L}^2(\Omega)}^2 + \|\nabla\varphi_1(t)\|_{\mathbf{L}^4(\Omega)}^4 + \|\nabla\varphi_2(t)\|_{\mathbf{L}^4(\Omega)}^4\Big).
    \end{align*}
    \end{theorem}

We now discuss the longtime behavior of each (weak) trajectory, focusing on dimension three. This is nontrivial, as we first aim at only dealing with a weak solution, without exploiting any asymptotic regularization effects. Indeed, extending the novel technique introduced in \cite{GP,GP2}, we will prove that weak solutions do not need regularization to converge to a unique equilibrium. This holds especially under weaker conditions on the coefficients, namely we can even assume $\nu\in C^0([-1,1])$. If the viscosity $\nu$ is more regular, since we can further rely on some local regularization properties, we can even obtain convergence to the equilibrium in higher-order norms. 

In order to carry out this plan, it is convenient to first state a consequence of Theorem \ref{THM:WeakExistence}, which follows from the energy inequality \eqref{energyineq} (see the proof of Theorem \ref{THM:WeakExistence}).
\begin{corollary}[Global properties of weak solutions]
    Let the assumptions of Theorem \ref{THM:WeakExistence} hold. Then any weak solution $(\bv,\varphi,\mu)$ to \eqref{Intro:Eq:1}--\eqref{Intro:Eq:6} departing from an initial datum $(\bv_0,\varphi_0)$, given by Theorem \ref{THM:WeakExistence} and globally defined in time in $(0,\infty)$, enjoys the global properties
            \begin{align}
            \begin{cases}
               \label{reg} \bv \in L^\infty(0,\infty;\mathbf{L}^2_\sigma(\Omega))\cap L^2(0,\infty;\mathbf{H}^1_\sigma(\Omega)), \\
                \varphi \in L^\infty(0,T;H^1(\Omega))\cap L^2_{uloc}([0,\infty);H^2_N(\Omega))\text{ with }|\varphi |<1\ \text{a.e. in }\Omega\times(0,\infty), \\
                \partial_t\bv \in L^{4/3}_{uloc}([0,\infty);\mathbf{H}^1_\sigma(\Omega)^\prime), \\
                \partial_t\varphi \in L^2(0,\infty;H^1(\Omega)^\prime), \\
                \mu \in L^2_{uloc}([0,\infty);L^2(\Omega)).
            \end{cases}
            \end{align}

   \label{thmglobal}
\end{corollary}
% \begin{remark}
%     {\color{red}[ISSUE] Even with the embedding we have $\partial_t\varphi\in L^\frac43 \Ld\cap L^2H^1(\Omega)'$, and $\varphi \in L^2H^2$...at most it holds
%     $$
%     \varphi \in BUC([0,\infty);H^\frac12(\Omega)),
%     $$
%     not enough for LS inequality, where we need $H^1$ :(
%     }
% \end{remark}

Let us then consider the set of admissible initial data:
\begin{align*}
\mathcal{H}_m:=\left\{\varphi\in H^1(\Omega): \Vert\varphi\Vert_{L^\infty(\Omega)}\leq 1,\quad \vert\overline{\varphi}\vert= m \right\},
\end{align*}
with $m\in[0,1)$, and fix an initial datum $(\bv_0,\varphi_0)\in \bLds\times\mathcal{H}_m$. Let then $(\bv,\varphi)$ be a weak global-in-time solution departing from $(\bv_0,\varphi_0)$, whose existence is ensured by Corollary \ref{thmglobal}. We introduce the $\omega$-limit set associated to $(\varphi,\bv)$, i.e.,
\begin{align*}
\omega(\bv,\varphi)=\{(\widetilde\bv,\widetilde{\varphi})\in \bLds\times \mathcal{H}_m:\exists t_n\to \infty \text{ s.t. }\bv(t_n)\rightharpoonup \widetilde\bv\text{ in }\bLd,\quad \varphi(t_n)\rightharpoonup  \widetilde{\varphi}\text{ in } \Hu\},
\end{align*}
where both the convergences are in the weak sense.
By \eqref{reg} we have $\bv\in BC_w([0,\infty);\bLds)$, whereas $\varphi\in BC_w([0,\infty);H^1(\Omega))$. Therefore $\omega(\bv,\varphi)$ is non-empty. We can also characterize the set $\omega(\bv,\varphi)$, showing that it is composed by equilibrium points (i.e., stationary solutions), which are defined as follows
\begin{definition}
	\label{defequil}$(\bv_\infty,\varphi_\infty)$ is an equilibrium point to problem \eqref{Intro:Eq:1}-\eqref{Intro:Eq:6} if $\bv_\infty=\mathbf 0$, and $\varphi_\infty\in \mathcal{H}_m\cap H^2_N(\Omega)$ satisfies the stationary Cahn-Hilliard equation
	\begin{align}
		-\Delta \varphi_\infty+f^\prime(\varphi_\infty)=\mu_\infty,\quad \text{in }\Omega,
		\label{conv1t}
	\end{align}
	where $\mu_\infty\in \R$ is a real constant.
\end{definition}
If we introduce the set of all the stationary points of our problem
\begin{equation}
\mathcal{S}:=\left\{(\mathbf 0,\varphi_\infty)\in \bLds\times\mathcal{H}_m\cap H^2_N(\Omega): \varphi_\infty\text{ satisfies }\eqref{conv1t}\right\},
\end{equation}
then we can prove that $\omega(\bv,\varphi)\subset \mathcal{S}$, together with the validity of an asymptotic strict separation property and the precompactness of trajectories in suitable spaces.  Note that, as in \cite{GP2}, we can directly apply De Giorgi's iterations to the parabolic equation, and this allows to show the asymptotic uniform-in-time strict separation property for each weak solution, a completely new property for such a coupled model. 
% \textit{without} resorting to the more delicate analysis of bad and good times relatively to $\norm{\nabla \mu(\cdot)}_{\bLd}$.  
In particular, in Section \ref{proofequil} we will prove that the following holds
\begin{lemma}
				\label{convaaa}
				Let the assumptions of Corollary \ref{thmglobal} hold. Then, we have
				$$
				    \omega(\bv,\vphi)\subset \mathcal{S}.
				$$
				Moreover,
                \begin{itemize}
                    \item[(i)]  $\omega(\bv,\vphi)$ is bounded in $\{ \mathbf 0\} \times H^2(\Omega)$, and there exists $\delta_1>0$ such that
				\begin{align}
				\| \vphi_\infty\|_{L^\infty(\Omega)}\leq 1-2\delta_1,\quad \forall \: (\mathbf 0,\,\vphi_\infty)\in \omega(\bv,\,\vphi).
				\label{sepaglobal}
                \end{align}
                Moreover, there exists $E_\infty\in \R$ such that it holds 
                \begin{align}
                E(\mathbf 0,\vphi_\infty)=E_\infty,\quad \forall (\mathbf 0,\vphi_\infty)\in \omega(\bv,\vphi).
                    \label{E2}
                \end{align}
                    \item[(ii)] (Asymptotic strict separation property). There exists $T_S>0$ and  $\delta\in(0,\delta_1)$ such that
    \begin{align}
        \sup_{t\geq T_S}\norm{\varphi(t)}_{L^\infty(\Omega)}\leq 1-\delta, \label{asympt}
    \end{align}
    i.e., the weak solution $\varphi$ is asymptotically strictly separated from the pure phases.
    \item[(iii)] (Precompactness of trajectories).
                  Additionally, the $\omega$-limit can be characterized as
\begin{align*}
\omega(\bv,\varphi)=\{(\widetilde\bv,\widetilde{\varphi})\in \bLds\times \mathcal{H}_m:\exists t_n\to \infty \text{ s.t. }\bv(t_n)\to\mathbf 0\text{ in }\bLd,\; \varphi(t_n)\to  \widetilde{\varphi}\text{ in } \Hu\},
\end{align*}
where the convergences are now strong, and it holds

                \begin{align}
	\lim_{t\to\infty}\dist_{ \bH^s\times H^r(\Omega)}((\bu(t),\varphi(t)),\omega(\bu,\vphi))=0,\label{convergence1b}
\end{align}
for any $s\in[-1,0)$, and any $r\in(0,1)$.
        \end{itemize}

			\end{lemma}
Let us now introduce the set of good times, following \cite{GP}: for a fixed $M>0$ and $T>0$, we define
 $$
 A_{M}(T):=\left\{t\geq T:\ \nu_*\norm{D\bv(t)}_{\mathbf L^2(\Omega)}^2+\norm{\mu(t)-\overline{\mu(t)}}_{\Ld}^2\leq  M^2\right\},
 $$
 which is of course of finite measure due to the energy inequality \eqref{ineq1}.
 In order to complete the proof, convergence \eqref{convergence1b} is not enough, since we would need it to hold with $r=1$. Therefore, following \cite{GP2}, we introduce the set \textquotedblleft good equilibrium points\textquotedblright :
 \begin{align}
 	\omega_g(\bv,\vphi):=\{(\mathbf 0,\vphi_*)\in \omega(\bv,\vphi):\ \exists t_n\to \infty \text{ s.t.} \{t_n\}\subset A_{M}(T)\text{ and }\vphi(t_n)\rightharpoonup \vphi_*,\text{ in }H^1(\Omega)\},
 	\label{goodt}
 \end{align}
so that $\omega_g(\bv,\vphi)\subset \omega(\bv,\vphi)\subset \mathcal S$, and $\omega_g(\bv,\vphi)$ is nonempty and does not depend on  $T>0$ (see \cite{GP2}). We can prove the following 
\begin{lemma}\label{twoparts1}
 	Let the assumptions of Corollary \ref{thmglobal} hold.
 	%    Given the set
 	%    \begin{align}
 		%    	A_\delta(t) &:= \{x\in\Omega:\; |\varphi(x,t)|\geq 1-{\delta_1}\},\quad t\geq0 ,\label{Adelta}
 		%    \end{align}
 	%    it holds
 	%    \begin{align}
 		%    	\lim_{t\to\infty }\normmm{A_\delta(t)}\to 0,
 		%    	\label{conerg}
 		%    \end{align}
 	%    where $\delta_1>0$ is given in \eqref{sepaglobal}.
% 	Then, for any $M>0$ there exists $\delta\in(0,\delta_1)$ and $T_S>0$ such that
% 	\begin{align}
% 		\sup_{t\in A_M(T_*)}\norm{\varphi(t)}_{L^\infty(\Omega)}\leq 1-\delta. \label{asympt}
% 	\end{align}
 	Then, for $s\in(0,1)$, $r\in(0,2)$, for any $\varepsilon>0$ there exists $T=T(\varepsilon)\geq 0$ such that 
 	\begin{align}
 		\dist_{\mathbf H^s(\Omega)\times H^r(\Omega)}((\bv(t),\vphi(t)),\omega_g(\bv,\vphi))<\varepsilon,\quad \forall t\in A_{M}(T).
 		\label{precomp1a} 
 	\end{align} 
    Also, the set $\omega_g(\bv,\vphi)$ is compact in $\{\mathbf 0\}\times H^r(\Omega)$, for $r\in(0,2)$.
 \end{lemma}
 With this crucial lemma at hand, 
% \begin{lemma}[Asymptotic strict separation]\label{twoparts}
%     Under the same assumptions of Lemma \ref{thmglobal},
% \end{lemma}
we can demonstrate that the $\omega$-limit of each trajectory is formed by a unique element if $F$ is analytic in $(-1,1)$, namely,
\begin{theorem}\label{uniqueeq}
	Let the assumptions of Theorem \ref{THM:WeakExistence} hold. Suppose additionally that $F$ is real analytic in $(-1,1)$. Then, a weak solution $(\bv,\vphi)$, departing from $(\bv_0,\vphi_0)\in \bLds\times \mathcal{H}_m$, enjoys the property $\omega(\bv,\vphi)=\{(\mathbf 0,\vphi_\infty)\}$. In particular, we have
	\begin{align}
	\lim_{t\to\infty}\norm{\bv(t)}_{\mathbf H^s(\Omega)}=0,\quad \lim_{t\to \infty}\Vert\vphi(t)-\vphi_\infty\Vert_{H^r(\Omega)}=0,
		\label{equil}
	\end{align}
    for any $s\in[-1,0),\ r\in(0,1)$.
 \end{theorem}
 \begin{remark}
   The strength of this theorem is the fact that it is only based on the existence of a global weak solution satisfying an energy inequality. This means that, as already observed, we can assume minimal conditions on the coefficients. Recalling Remark \ref{visco}, the theorem above still holds if the viscosity is simply $\nu\in C^0([-1,1])$ and not necessarily Lipschitz continuous. Also, we do not need to exploit any regularization of the solution.
 \end{remark}
 \begin{remark}
     If we restrict ourselves to the good times, actually we can prove that the convergence to equilibrium holds in higher-order norms for $\bv$ and $\vphi$, by \eqref{precomp1a}.  
 \end{remark}
Furthermore, if $\nu \in W^{1,\infty}(\R)$, thanks to the combination of asymptotic strict separation property and conditional weak-strong uniqueness principle, we can actually prove more. Indeed, we can show that any weak solution asymptotically regularizes, so that,  as a trivial consequence  the convergences \eqref{equil} hold in  higher-order norms.
\begin{theorem}
\label{finn}
   Under the same assumptions of Theorem \ref{uniqueeq}, any weak solution $(\bv,\vphi)$, departing from $(\bv_0,\vphi_0)\in \bLds\times \mathcal{H}_m$, asymptotically regularizes, namely there exists $T_R>0$, depending only on $E(\bv_0,\vphi_0)$, such that  
   \begin{align}
                \begin{cases}
                    \bv \in L^\infty(T_R,\infty;\mathbf{H}^1_\sigma(\Omega))\cap L^2_{uloc}([T_R,\infty);\mathbf{H}^2(\Omega))\cap H^1_{uloc}([T_R,\infty);\mathbf{L}^2_\sigma(\Omega)),\\
                    \varphi\in L^\infty(T_R,\infty;H^2_N(\Omega))\cap H^1_{uloc}([T_R,\infty);H^1(\Omega)),\\
                    \mu\in L^\infty(T_R,\infty;L^2(\Omega))\cap L^2_{uloc}([T_R,\infty);H^1(\Omega)).
                \end{cases}
            \end{align}
   Additionally, $(\bv,\vphi)$ converges to a single equilibrium $(\mathbf 0,\vphi_\infty)\in \mathcal S$ in the sense that
	\begin{align}
\lim_{t\to\infty}\norm{\bv(t)}_{\bH^s(\Omega)}=0,\quad \lim_{t\to \infty}\Vert\vphi(t)-\vphi_\infty\Vert_{H^r(\Omega)}=0,\quad\forall s\in(0,1),\ r\in[1,2).
		\label{equil1}
	\end{align}
    % In particular, this also entails that the weak solution asymptotically regularizes, namely there exists $T_R>0$ such that
    % \begin{align*}
    % &\bv\in BC_w([T_R,\infty);\bH^s(\Omega)),\quad \forall s\in(0,1),\\&
    % \vphi\in BC_w([T_R,\infty);H^r(\Omega)), \quad \forall r\in[1,2).
    % \end{align*}
\end{theorem}
% \begin{remark}
%     Exploiting the weak-strong uniqueness principle of Theorem \ref{THM:WSU} together with the asymptotic strict separation of weak solutions, which are at the basis of the proof of Theorem \ref{finn} above, we could have avoided to resort to the notion of good times and directly have shown that the trajectories are precompact in $\bLds\times H^1(\Omega)$, immediately leading to prove the convergence of each trajectory to a unique equilibrium. Nevertheless, since weak-strong uniqueness is not usually available in similar problems (see the discussions in \cite{GP,GP2}), we preferred to present in the separate Theorem \ref{uniqueeq} which properties can be obtained only relying on weak solutions, without the help of regularization. For a prototypical example, recalling Remark \ref{visco}, Theorem \ref{uniqueeq} also holds in case $\nu\in C^0([-1,1])$, whereas for Theorem \ref{finn} to hold we need the existence of strong solutions, and thus we need to assume at least $\nu\in W^{1,\infty}(-1,1)$.
% \end{remark}

% \section{TO DO LIST}
% {\color{blue}Christoph:
% \begin{itemize}
%     \item Check the blue comments: minor issues like adding some modulus.
%     \item Major issue: check the proof of the existence of strong solutions, since in the Galerkin scheme some estimates are not allowed. Also, add the $\lambda$ approximation.
% \end{itemize}}
\begin{remark}
    We point out that, differently from Theorem \ref{uniqueeq}, to obtain the stronger results above, the assumption $\nu\in W^{1,\infty}(\R)$ cannot be relaxed to $\nu$ just continuous, since it is necessary, for instance, to apply Theorem \ref{THM:StrongWP}.
\end{remark}
    \section{Proof of Theorem \ref{THM:WeakExistence}}
    \label{weak}
    As mentioned in the Introduction, we report here a proof based on a full Galerkin scheme (cf. \cite{GGW2022}) which will be helpful in the subsequent section.
    \subsection{Regularization of the singular potential}
    First of all, due to the properties of the singular part $F$ (see \ref{ASS:S1}), it follows that $F$ is proper, convex and lower-semicontinuous on $[-1,1]$ and $F^\prime$ is continuous and non-decreasing on $(-1,1)$. Thus, $F$ can be identified as a maximal monotone graph on $\R$. Therefore, we can introduce its Yosida approximation $F_\lambda$ for $\lambda>0$.
    Recall that this approximation has the following useful properties (cf. \cite{barbu-monot}):
    \begin{enumerate}
    \item []
        \item\label{Yos:1} $F_\lambda$ is convex and converges to $F$ pointwise in $\R$ and is monotonically increasing as $\lambda\rightarrow0$;
        \item\label{Yos:2} $F_\lambda\in C^{1,1}(\R)$ and the Lipschitz constant of $F_\lambda$ is $\tfrac{1}{\lambda}$;
        \item\label{Yos:3} $|F_\lambda^\prime|$ converges to $|F^\prime|$ pointwise in $\R$ and is monotonically increasing as $\lambda\rightarrow0$;
        \item\label{Yos:4} $F_\lambda(0) = F_\lambda^\prime(0)=0$ for all $\lambda>0$;
        \item\label{Yos:5} for all $K>0$ there exists $\lambda_0>0$ such that, for all $\lambda\in(0,\lambda_0)$,
        \begin{align*}
            F_\lambda(s) \geq Ks^2-K \quad \text{ for all }s\in\R;
        \end{align*}
        \item\label{Yos:6} for all $\lambda\in(0,1]$, it holds
        \begin{align*}
            F_\lambda^{\prime\prime}(s) \geq \frac{\theta}{1+\theta} \quad \text{ for all }s\in\R.
        \end{align*}
        \item []
    \end{enumerate}
    Let us now introduce the regularized potential
    \begin{align*}
        f_\lambda(s) := F_\lambda(s) - \frac{\theta_0}{2}s^2\quad \text{ for }s\in\R,
    \end{align*}
    and consider the $\lambda$-approximated problem
    \begin{align}
        \partial_t\bv_\lambda + (\bv_\lambda\cdot\nabla)\bv_\lambda - \text{div}(\nu(\varphi_\lambda)D\bv_\lambda) + \nabla \pi_\lambda &= \mu_\lambda\nabla\varphi_\lambda &\quad&\text{ in }\Omega_T, \label{Eq:lam1}\\
        \nabla\cdot\bv_\lambda &= 0 &\quad&\text{ in }\Omega_T, \label{Eq:lam2}\\
        \partial_t\varphi_\lambda + \bv_\lambda\cdot\nabla\varphi_\lambda + \mu_\lambda - \overline{\mu_\lambda} &= 0 &\quad&\text{ in }\Omega_T, \label{Eq:lam3}\\
        \mu_\lambda &= - \Delta\varphi_\lambda + f_\lambda^\prime(\varphi_\lambda) &\quad&\text{ in }\Omega_T, \label{Eq:lam4}\\
        \bv_\lambda = \partial_\mathbf{n}\varphi_\lambda &= 0 &\quad&\text{ on }\partial\Omega_T, \label{Eq:lam5}\\
        \bv_\lambda(0) = \bv_0,\quad \varphi_\lambda(0) &= \varphi_0 &\quad&\text{ in }\Omega. \label{Eq:lam6}
    \end{align}
    \subsection{Discretization via Galerkin scheme}\label{Step:2}
    We intend employ a Galerkin scheme to solve problem \eqref{Eq:lam1}--\eqref{Eq:lam6}. It is well-known that the Stokes operator $A_S$ possesses a sequence of eigenvalues $(\alpha_j)_{j\in\mathbb{N}}$ along with a sequence of corresponding eigenfunctions $(\bw_j)_{j\in\mathbb{N}} \subset \bH^1_\sigma(\Omega)$, which can be chosen in such a way that they form an orthonormal basis of $\mathbf{L}^2_\sigma(\Omega)$ as well as an orthogonal Schauder's basis of $\bH^1_\sigma(\Omega)$. Using regularity theory for the Stokes operator, recalling that the domain is $C^2$, we further infer that $\bw_j \in \mathbf{H}^2(\Omega)$. Then, for any $m\in \mathbb{N}$, we introduce the finite-dimensional subspace
    \begin{equation*}
    	\mathbf{V}_m := \mathrm{span}\,\big\{ \bw_1, ..., \bw_m \big\} \subset \bH^1_\sigma(\Omega),
    \end{equation*}
    and $\mathbb{P}_m: \bH^1_\sigma(\Omega) \to \mathbf{V}_m$ denotes the $\mathbf{L}^2_\sigma(\Omega)$-orthogonal projection of $\bH^1_\sigma(\Omega)$ onto $\mathbf{V}_m$.
    \newline
    Analogously, the negative Laplacian with homogeneous Neumann boundary condition possesses a sequence of eigenvalues $(\beta_k)_{k\in\mathbb{N}}$ along with eigenvectors $(u_k)_{k\in\mathbb{N}}\subset H^2_N(\Omega)$. It is well-known that, up to renormalization, the set $(u_k)_{k\in\mathbb{N}}$ is orthonormal in $L^2(\Omega)$ and orthogonal in $H^1(\Omega)$. For any $m\in\mathbb{N}$, we introduce the finite-dimensional subspace
    \begin{align*}
        X_m := \mathrm{span}\,\{u_1,...,u_m\} \subset H^1(\Omega),
    \end{align*}
    and $P_m: H^1(\Omega)\rightarrow X_m$ denotes the $L^2(\Omega)$-orthogonal projection of $H^1(\Omega)$ onto $X_m$.
    Therefore, the discretized problem reads as follows
    \begin{align}
        \partial_t\bv_{\lambda,m} + (\bv_{\lambda,m}\cdot\nabla)\bv_{\lambda,m} - \text{div}(\nu(\varphi_{\lambda,m})D\bv_{\lambda,m}) + \nabla \pi_{\lambda,m} &= \mu_{\lambda,m}\nabla\varphi_{\lambda,m}, \label{Eq:m:lam1}\\
        \nabla\cdot\bv_{\lambda,m} &= 0, \label{Eq:m:lam2}\\
        \partial_t\varphi_{\lambda,m} + \bv_{\lambda,m}\cdot\nabla\varphi_{\lambda,m} + \mu_{\lambda,m} - \overline{\mu_{\lambda,m}} &= 0,  \label{Eq:m:lam3}\\
        \mu_{\lambda,m} &= - \Delta\varphi_{\lambda,m} + P_mf_\lambda^\prime(\varphi_{\lambda,m}), \label{Eq:m:lam4}
    \end{align}
    in $\Omega\times (0,T)$, with
    \begin{align}
        \bv_{\lambda,m} = \partial_\mathbf{n}\varphi_{\lambda,m} &= 0 &\quad&\text{ a.e. on }\partial\Omega_T, \label{Eq:m:lam5}\\
        \bv_{\lambda,m}(0) = \mathbb P_m\bv_0=:\bv_{0,m},\quad \varphi_{\lambda,m}(0) &= P_m\varphi_0=:\varphi_{0,m} &\quad&\text{ in }\Omega. \label{Eq:m:lam6}
    \end{align}
    Fixed any $\lambda\in(0,1)$ and $m\in\mathbb{N}$, we look for a solution to \eqref{Eq:m:lam1}--\eqref{Eq:m:lam6} of the form
    \begin{align}\label{Def:sol:m:lam}
        \bv_{\lambda,m} := \sum_{j=1}^m g_{\lambda,m}^j\bw_j,\quad \varphi_{\lambda,m} := \sum_{k=1}^md_{\lambda,m}^k u_k,\quad \mu_{\lambda,m} := \sum_{k=1}^mc_{\lambda,m}^k u_k,
    \end{align}
    where
    \begin{align*}
        \mathbf{g}_{\lambda,m}:[0,T]\rightarrow\R,\quad
        \mathbf{d}_{\lambda,m}:[0,T]\rightarrow\R, \quad
        \mathbf{c}_{\lambda,m}:[0,T]\rightarrow\R
    \end{align*}
    are suitable vectors of coefficients. Exploiting the relations in \eqref{Def:sol:m:lam} and testing \eqref{Eq:m:lam1}, \eqref{Eq:m:lam3}--\eqref{Eq:m:lam4} by $\bw_i$ and $u_i$, $i\in\{1,\ldots,m\}$, respectively, we infer that the coefficient vectors $\mathbf{g}_{\lambda,m}, \mathbf{d}_{\lambda,m}, \mathbf{c}_{\lambda,m}$ satisfy the following system of ordinary differential equations,
    \begin{align}
        \frac{\text{d}}{\text{d}t}g_{\lambda,m}^i = \sum_{l=1}^m\sum_{k=1}^m d_{\lambda,m}^l c_{\lambda,m}^k\int_\Omega u_l\nabla u_k\cdot\bw_i\dx
        - \sum_{j=1}^m\sum_{k=1}^mg_{\lambda,m}^jg_{\lambda,m}^k \,b(\bw_j,\bw_k,\bw_i) \label{ODE:1}\\
        -\int_\Omega\nu\Big(\sum_{k=1}^mc_{\lambda,m}^ku_k\Big)D\Big(\sum_{j=1}^mg_{\lambda,m}^j\bw_j\Big):D\bw_i\dx, \nonumber\\
        \frac{\text{d}}{\text{d}t}c_{\lambda,m}^i + d_{\lambda,m}^i = -\int_\Omega\Big(\sum_{j=1}^mg_{\lambda,m}^j\bw_j\Big)\cdot\nabla\Big(\sum_{k=1}^mc_{\lambda,m}^ku_k\Big)u_i\dx + \sum_{k=1}^md_{\lambda,m}^k\int_\Omega\overline{u_k}u_i\dx,\label{ODE:2}\\
        d_{\lambda_m}^i = \beta_ic_{\lambda,m}^i + \int_\Omega f_{\lambda}^\prime\big(\sum_{k=1}^mc_{\lambda,m}^ku_k\big)u_i\dx,\label{ODE:3}\\
        g_{\lambda,m}^i(0) = (\bv_0,\bv_i),\label{ODE:4}\\
        c_{\lambda,m}^i(0) = (c_0,u_i). \label{ODE:5}
    \end{align}
    Then, Peano's theorem ensures the existence of a $C^1$-solution to this initial value problem in a right-maximal time interval $[0,T_m^*)$, $T_m^*>0$.

    \subsection{Energy estimate} We test \eqref{Eq:m:lam1} with $\bv_{\lambda,m}$  and \eqref{Eq:m:lam3} with $\mu_{\lambda,m}-\overline{\mu_{\lambda,m}}$. Setting
    \begin{align*}
        E_{\lambda,m}(t) := \frac{1}{2}\|\bv_{\lambda,m}(t)\|_{\mathbf{L}^2_\sigma(\Omega)}^2 + \frac{1}{2}\|\nabla\varphi_{\lambda,m}(t)\|_{\mathbf{L}^2(\Omega)}^2 + \int_\Omega f_\lambda(\varphi_{\lambda,m}(t))\dx, \quad t\in[0,T],
    \end{align*}
    and using equation \eqref{Eq:m:lam4}, we infer the discretized energy inequality
    \begin{align*}
        \frac{\text{d}}{\text{d}t}E_{\lambda,m} + \nu_*\|D\bv_{\lambda,m}\|_{\mathbf{L}^2(\Omega)}^2 + \|\mu_{\lambda,m}-\overline{\mu_{\lambda,m}}\|_{L^2(\Omega)}^2 \leq 0.
    \end{align*}
    %Since $\varphi_0\in H^1(\Omega)$, then we have
    %\begin{align*}
    %   \varphi_0 = \sum_{k=1}^\infty (\varphi_0,u_i) u_i,
    %\end{align*}
    %and $\varphi_{0,m}\rightarrow\varphi_0$ in $H^1(\Omega)$ as $m\rightarrow\infty$. Thus, the sequence $(\varphi_{0,m})_{m\in\mathbb{N}}$ is bounded in $H^1(\Omega)$. Moreover, since the orthogonal projection $\mathbb{P}_m$ %is Lipschitz-continuous, we have
    %\begin{align*}
    %   \|\bv_{0,m}\|_{\mathbf{L}^2_\sigma(\Omega)} = \|\mathbb{P}_m\bv_0\|_{\mathbf{L}^2_\sigma(\Omega)} \leq \|\bv_0\|_{\mathbf{L}^2_\sigma(\Omega)}.
    %\end{align*}
    %Since $F_\lambda^\prime$ is linearly bounded, it follows that $F_\lambda$ is quadratically bounded and hence,
    %\begin{align*}
    %    \|f_\lambda(\varphi_{0,m})\|_{L^1(\Omega)} \leq C_\lambda(1+\|\varphi_0\|_{L^2(\Omega)}^2)
    %\end{align*}
    %by definition of $f_\lambda$,
    %where the constant $C_\lambda>0$ is independent of $m\in\mathbb{N}$. Altogether, this shows
    It is easy to show that
    \begin{align*}
        |E_{\lambda,m}(0)| \leq C,
    \end{align*}
    where the constant $C>0$ only depends on $\lambda, \|\varphi_0\|_{H^1(\Omega)}, \|\bv_0\|_{\mathbf{L}^2_\sigma(\Omega)}$. This yields
    \begin{align}\label{Inequ:EN}
        E_{\lambda,m}(t) + \nu_*\int_s^t\|D\bv_{\lambda,m}(t)\|_{\mathbf{L}^2(\Omega)}^2\dt + \int_0^s\|\mu_{\lambda,m}(t)-\overline{\mu_{\lambda,m}}(t)\|_{L^2(\Omega)}^2\dt\leq  E_{\lambda,m}(s) \leq C,
    \end{align}
    for all $m\in\mathbb{N}$, all $t\in [0,T_m^*]$ and almost all $s\in[0,T_m^*]$, including $s=0$.

    \subsection{A priori estimates} The energy inequality \eqref{Inequ:EN} implies that the solution is bounded uniformly with respect to $m$ in the following spaces
    \begin{align}
        \bv_{\lambda,m} &\in L^\infty(0,T_m^*;\mathbf{L}^2_\sigma(\Omega)), \label{Est:apr1}\\
        \nabla\varphi_{\lambda,m} &\in L^\infty(0,T_m^*;\mathbf{L}^2(\Omega)),\label{Est:apr2}\\
        D\bv_{\lambda,m}&\in L^2(0,T_m^*;\mathbf{L}^2(\Omega)),\label{Est:apr3}\\
        \mu_{\lambda,m}-\overline{\mu_{\lambda,m}}&\in L^2(0,T_m^*;L^2(\Omega)). \label{Est:apr4}
    \end{align}
    Also, Korn's inequality entails that
    \begin{align*}
        \bv_{\lambda,m}\in L^2(0,T_m^*;\mathbf{H}^1_\sigma(\Omega))
    \end{align*}
    and bounded uniformly with respect to $m$.

    Testing \eqref{Eq:m:lam3} with $1$, we get
    \begin{align*}
        \int_\Omega\varphi_{\lambda,m}(s)\dx = \int_\Omega\varphi_0\dx
    \end{align*}
    for all $s\in[0,T_m^*]$. Thus, on account of \eqref{Est:apr2} and the Poincaré inequality, it follows that
    \begin{align}
        \label{Est:apr5}
        \varphi_{\lambda,m}\in L^\infty(0,T_m^*;H^1(\Omega))
    \end{align}
    is bounded uniformly with respect to $m$. From equation \eqref{Eq:m:lam4}, we get
    \begin{align*}
        \overline{\mu_{\lambda,m}} = \overline{f_\lambda^\prime(\varphi_{\lambda,m})}.
    \end{align*}
    Using the fact that $F_\lambda^\prime$ is $\frac1\lambda$ Lipschitz-continuous, the definition of $f_\lambda$ yields
    \begin{align*}
        \|\overline{\mu_{\lambda,m}}\|_{L^\infty(0,T_m^*)}  \leq C_\lambda,
    \end{align*}
    where the constant $C_\lambda >0$ does not depend on $m$. Combining this estimate with \eqref{Est:apr4}, we conclude that
    \begin{align}
        \label{Est:apr6}
        \mu_{\lambda,m}\in L^2(0,T_m^*;L^2(\Omega))
    \end{align}
    is bounded uniformly with respect to $m$. Summing up, we have derived uniform estimates in the following spaces
    \begin{align}
        \bv_{\lambda,m} &\in L^2(0,T_m^*;\mathbf{H}^1_\sigma(\Omega)) \cap L^\infty(0,T_m^*;\mathbf{L}^2_\sigma(\Omega)), \label{EST:apr:v}\\
        \varphi_{\lambda,m}&\in L^\infty(0,T_m^*;H^1(\Omega)),\label{EST:apr:phi}\\
        \mu_{\lambda,m}&\in L^2(0,T_m^*;L^2(\Omega)).\label{EST:apr:mu}
    \end{align}
    This allows us to deduce that the solutions to problem \eqref{ODE:1}--\eqref{ODE:5} are uniformly bounded in $[0,T_m^*]$ with respect to $m$.
    Therefore, we can extend them to the whole time interval $[0,\infty)$.

    \subsection{Further estimates} Testing \eqref{Eq:m:lam3} with $\zeta\in L^2(0,T;H^1(\Omega))$ and using the bounds derived in the previous section, we arrive at
    \begin{align}
        \nonumber&\left|\int_{\Omega_T}\partial_t\varphi_{\lambda,m}\zeta\dx \right| = \left|\int_{\Omega_T}\partial_t\varphi_{\lambda,m}P_m\zeta\dx \right| \\
        &\leq \Big(\|\mu_{\lambda,m}-\overline{\mu_{\lambda,m}}\|_{L^2(0,T;L^2(\Omega))} + \|\bv_{\lambda}\|_{L^2(0,T;\mathbf{L}^4_\sigma(\Omega))}\|\varphi_{\lambda,m}\|_{L^\infty(0,T;L^4(\Omega))}\Big)\|\zeta\|_{L^2(0,T;H^1(\Omega))},\label{H1p}
    \end{align}
    which implies that $\partial_t\varphi_{\lambda,m}\in L^2(0,T;H^1(\Omega)^\prime)$. Hence
    $\varphi_{\lambda,m}\in H^1(0,T;H^1(\Omega)^\prime)$
    is uniformly bounded with respect to $m$.
    Using the Lipschitz continuity of $f_\lambda^\prime$ and the fact that $f_\lambda^\prime(0)=0$, we infer
    \begin{align}\label{Est:F:Lipschitz}
        |f_\lambda^\prime(\varphi_{\lambda,m})| \leq C_\lambda|\varphi_{\lambda,m}|.
    \end{align}
    Thus, we deduce (see \eqref{EST:apr:phi})
    \begin{align}
        \|f_\lambda^\prime(\varphi_{\lambda,m})\|_{L^2(0,T;L^2(\Omega))} \leq C_\lambda,
    \end{align}
    where $C_\lambda>0$ depends on $\lambda$ but is independent from $m$. In particular, a comparison in \eqref{Eq:m:lam4} implies
    \begin{align}
        \|\varphi_{\lambda,m}\|_{L^2(0,T;H^2(\Omega))} \leq C_\lambda
    \end{align}
    for some constant $C_\lambda>0$ similar to the one right above. Next, we test \eqref{Eq:m:lam1} with $\mathbf{\xi}\in L^4(0,T;\mathbf{H}^1_\sigma(\Omega))$.
    Then, on account of the uniform bounds derived so far and using the Gagliardo--Nirenberg inequality to handle $\mu_{\lambda,m}\nabla\varphi_{\lambda,m}$, we deduce
    \begin{align}
        \partial_t\bv_{\lambda,m}\in L^{4/3}(0,T;\mathbf{H}^1_\sigma(\Omega)^\prime).
    \end{align}
    \subsection{Compactness} Thanks to the uniform estimates derived in the previous steps, we can now pass to the limit $m\rightarrow\infty$, up to a convenient subsequence, keeping $\lambda\in(0,1)$ fixed. By standard compactness and diagonal arguments, we deduce that there exist $\varphi_\lambda,\bv_\lambda,\mu_\lambda$, globally defined on $(0,\infty)$, such that, up to a non-relabeled subsequence, the following convergences hold
    \begin{align*}
        \varphi_{\lambda,m} &\rightharpoonup \varphi_\lambda && \text{weakly-$*$ in }L^\infty(0,T;H^1(\Omega)),\\
        \varphi_{\lambda,m} &\rightharpoonup \varphi_\lambda && \text{weakly in }L^2(0,T;H^2(\Omega)),\\
        \partial_t\varphi_{\lambda,m}&\rightharpoonup \partial_t\varphi_\lambda && \text{weakly in }L^2(0,T;H^1(\Omega)^\prime),\\
        \varphi_{\lambda,m} &\rightarrow \varphi_\lambda && \text{strongly in }L^2(0,T;H^1(\Omega))\cap C^0([0,T];L^2(\Omega)),\\
        \varphi_{\lambda,m} &\rightarrow \varphi_\lambda && \text{a.e. in }\Omega_T,\\
        \mu_{\lambda,m} &\rightarrow \mu_\lambda && \text{weakly in }L^2(0,T;L^2(\Omega)),\\
        \bv_{\lambda,m} &\rightharpoonup \bv_\lambda && \text{weakly-$*$ in }L^\infty(0,T;\mathbf{L}^2_\sigma(\Omega)),\\
        \bv_{\lambda,m} &\rightharpoonup \bv_\lambda && \text{weakly in }L^2(0,T;\mathbf{H}^1_\sigma(\Omega)),\\
        \partial_t\bv_{\lambda,m} &\rightharpoonup \partial_t\bv_\lambda && \text{weakly in }L^{4/3}(0,T;\mathbf{H}^1_\sigma(\Omega)^\prime),\\
        \bv_{\lambda,m} &\rightarrow \bv_\lambda && \text{strongly in }L^2(0,T;\mathbf{L}^2_\sigma(\Omega)),\\
        \bv_{\lambda,m} &\rightarrow \bv_\lambda && \text{a.e. in }\Omega_T.
    \end{align*}
    \subsection{Passage to the limit $m\rightarrow\infty$}\label{Step:7}First of all, the Lipschitz-continuity of $f_\lambda^\prime$ implies
    \begin{align*}
        f_\lambda^\prime(\varphi_{\lambda,m}) \rightarrow f_\lambda^\prime(\varphi_\lambda) \quad \text{in }L^2(0,T;L^2(\Omega)).
    \end{align*}
    Next, we consider the Korteweg stress. Due to the weak convergence of $(\mu_{\lambda,m})_{m\in\mathbb{N}}$ in $L^2(0,T;L^2(\Omega))$ and the strong convergence of $(\nabla\varphi_{\lambda,m})_{m\in\mathbb{N}}$ in $L^2(0,T;\mathbf{L}^2(\Omega))$, we can pass to the limit $m\rightarrow\infty$ and obtain
    \begin{align*}
        (\mu_{\lambda,m},\nabla\varphi_{\lambda,m}\cdot\bw)  \rightarrow (\mu_\lambda,\nabla\varphi_\lambda\cdot \bw),
    \end{align*}
    for all the divergence-free functions $\bw \in \mathbf{C}^\infty_0(\Omega)$.
    Moreover, it is standard to check that for the nonlinear term in the Navier--Stokes equation, it holds
    \begin{align*}
        (\bv_{\lambda,m}\cdot\nabla)\bv_{\lambda,m}\rightarrow (\bv_{\lambda}\cdot\nabla)\bv_{\lambda}\quad \text{in }L^{4/3}(0,T;\mathbf{H}^1_\sigma(\Omega)^\prime)
    \end{align*}
    as $m\rightarrow\infty$.
    It remains to consider the convective term in the Cahn--Hilliard equation. There, we use a similar argument to conclude
    \begin{align*}
        \bv_{\lambda,m}\cdot\nabla\varphi_{\lambda,m} \rightharpoonup \bv_{\lambda}\cdot\nabla\varphi_{\lambda}\quad \text{in }L^1(0,T;L^{3/2}(\Omega))
    \end{align*}
    as $m\rightarrow\infty$.
    With these convergences at hand, we can now pass to the limit along a suitable non-relabeled sequence $m\rightarrow\infty$ in the weak formulation of \eqref{Eq:m:lam1}--\eqref{Eq:lam6} and then use a density argument to recover the weak formulation of problem \eqref{Eq:lam1}--\eqref{Eq:lam6}. Concerning the initial conditions, due to the convergence $\varphi_{0,m}\rightarrow\varphi_0$ in $L^2(\Omega)$ and the strong convergence $\varphi_{\lambda,m}\rightarrow\varphi_\lambda$ in $C^0([0,T];L^2(\Omega))$, it immediately follows that $\varphi_\lambda(0)=\varphi_0$ in $L^2(\Omega)$. Since the sequence $(\bv_{\lambda,m})_{m\in\mathbb{N}}\subset W^{1,4/3}(0,T;\mathbf{H}^1_\sigma(\Omega)^\prime)\hookrightarrow C^0([0,T];\mathbf{H}^1_\sigma(\Omega)^\prime)$ is bounded, there exists a subsequence such that $\bv_{\lambda,m}(0)\rightharpoonup\bv_\lambda(0)$ in $\mathbf{H}^1_\sigma(\Omega)^\prime$. On the other hand, since $\bv_{0,m}\rightarrow\bv_0$ in $\mathbf{L}^2_\sigma(\Omega)\hookrightarrow\mathbf{H}^1_\sigma(\Omega)^\prime$, we deduce $\bv_{\lambda}(0) = \bv_0$ in $\mathbf{H}^1_\sigma(\Omega)^\prime$.

    \subsection{Uniform estimates with respect to $\lambda$} First of all, recalling the definition of $f_\lambda$ we have that $\|f_\lambda(\varphi_0)\|_{L^1(\Omega)} \leq \|f(\varphi_0)\|_{L^1(\Omega)}$. Hence, exploiting the compactness arguments above we can pass to the limit in \eqref{Inequ:EN}, exploiting the lower-semicontinuity of the norms, to deduce
    \begin{align}
        E_\lambda(t) + \nu_*\int_s^t\|D\bv_\lambda(t)\|^2_{\mathbf{L}^2(\Omega)}\dt + \int_s^t\|\mu_\lambda(t)-\overline{\mu_\lambda}(t)\|_{L^2(\Omega)}^2\dt \leq E_\lambda(s)\leq C\label{energyineq}
    \end{align}
    for all $t\in[0,T]$ and almost all $s\in[0,T]$, including $s=0$,
    where the constant $C$ only depends on $\|\varphi_0\|_{H^1(\Omega)}, \|\bv_0\|_{\mathbf{L}^2_\sigma(\Omega)}$. This yields a priori estimates, which are uniform in $\lambda$, in the following spaces
    \begin{align}
        \bv_{\lambda} &\in L^\infty(0,T;\mathbf{L}^2_\sigma(\Omega))\cap L^2(0,T;\mathbf{H}^1_\sigma(\Omega)), \label{Est:apr:lam1}\\
        \nabla\varphi_{\lambda} &\in L^\infty(0,T;\mathbf{L}^2(\Omega)),\label{Est:apr:lam2}\\
        D\bv_{\lambda}&\in L^2(0,T;\mathbf{L}^2(\Omega)),\label{Est:apr:lam3}\\
        \mu_{\lambda}-\overline{\mu_{\lambda}}&\in L^2(0,T;L^2(\Omega)). \label{Est:apr:lam4}
    \end{align}
    Moreover, testing \eqref{Eq:lam3} with $1$, we infer that
    \begin{align*}
        \int_\Omega\varphi_\lambda(s)\dx = \int_\Omega\varphi_0\dx
    \end{align*}
    for any $s\in[0,T]$. Thus, by \eqref{Est:apr:lam2} and Poincaré's inequality, we obtain
    \begin{align}\label{Est:apr:lam5}
        \varphi_\lambda\in L^\infty(0,T;H^1(\Omega)).
    \end{align}
    From equation \eqref{Eq:lam4}, we obtain
    \begin{align*}
        \overline{\mu_\lambda} = \overline{f^\prime_\lambda(\varphi_\lambda)}.
    \end{align*}
    Testing \eqref{Eq:lam4} with $\varphi_\lambda-\overline{\varphi_\lambda}$, it follows that
    \begin{align*}
        \|\nabla\varphi_\lambda\|_{\mathbf{L}^2(\Omega)}^2 + \int_\Omega F_\lambda^\prime(\varphi_\lambda)(\varphi_\lambda-\overline{\varphi_\lambda})\dx = \int_\Omega (\mu_\lambda-\overline{\mu_\lambda})\varphi_\lambda\dx + \theta_0\int_\Omega \varphi_\lambda (\varphi_\lambda-\overline{\varphi_\lambda})\dx.
    \end{align*}
    Using the well-known inequality
    \begin{align*}
        \|F_\lambda^\prime(\varphi_\lambda)\|_{L^1(\Omega)}\leq C(\overline{\varphi_0})\Big(1+\normmm{\int_\Omega F_\lambda^\prime(\varphi_\lambda)(\varphi_\lambda-\overline{\varphi_\lambda})\dx}\Big),
    \end{align*}
    estimate \eqref{Est:apr:lam5} and exploiting Young's inequality, we then arrive at
    \begin{align*}
        \|F_\lambda^\prime(\varphi_\lambda)\|_{L^1(\Omega)} \leq C(1+ \|\mu_\lambda-\overline{\mu_\lambda}\|_{L^2(\Omega)}).
    \end{align*}
    In particular, by definition of $f_\lambda$ and \eqref{Est:apr:lam5} this yields
    \begin{align*}
        |\overline{\mu_\lambda}| \leq C(1+ \|\mu_\lambda-\overline{\mu_\lambda}\|_{L^2(\Omega)}),
    \end{align*}
    and hence, by \eqref{Est:apr:lam4},
    \begin{align*}
        \int_0^T|\overline{\mu_\lambda}(t)|^2\dt \leq C.
    \end{align*}
    Altogether, this shows
    \begin{align}\label{Est:apr:lam6}
        \mu_\lambda\in L^2(0,T;L^2(\Omega)).
    \end{align}
    Arguing as before, we can use \eqref{Est:apr:lam1},\eqref{Est:apr:lam4} and \eqref{Est:apr:lam5} to deduce
    \begin{align*}
        \varphi_\lambda \in H^1(0,T;H^1(\Omega)^\prime).
    \end{align*}
    We now test \eqref{Eq:lam4} with $F_\lambda^\prime(\varphi_\lambda)$. This gives
    \begin{align*}
        \|F_\lambda^\prime(\varphi_\lambda)\|_{L^2(\Omega)}^2 = -\int_\Omega\mu_\lambda F_\lambda^\prime(\varphi_\lambda)\dx - \int_\Omega\nabla\varphi_\lambda\cdot\nabla(F_\lambda^\prime(\varphi_\lambda))\dx + \theta_0\int_\Omega\varphi_\lambda F_\lambda^\prime(\varphi_\lambda)\dx.
    \end{align*}
    Owing to property \ref{Yos:6} of $F_\lambda^\prime$ and the chain rule, we observe
    \begin{align*}
        - \int_\Omega\nabla\varphi_\lambda\cdot\nabla(F_\lambda^\prime(\varphi_\lambda))\dx =- \int_\Omega F_\lambda^{\prime\prime}(\varphi_\lambda)|\nabla \varphi_\lambda|^2\dx \leq 0.
    \end{align*}
    Thus, by means of Hölder's and Young's inequality, we infer that
    \begin{align*}
        \|F_\lambda^\prime(\varphi_\lambda)\|_{L^2(\Omega)}^2 \leq C\big(\|\mu_\lambda\|_{L^2(\Omega)}^2+\|\varphi_\lambda\|_{L^2(\Omega)}^2\big).
    \end{align*}
    Hence, by \eqref{Est:apr:lam5} and \eqref{Est:apr:lam6}, we have
    \begin{align}\label{Est:apr:lam7}
        F_\lambda^\prime(\varphi_\lambda) \in L^2(0,T;L^2(\Omega)).
    \end{align}
    Therefore, by elliptic regularity theory and \eqref{Est:apr:lam5}--\eqref{Est:apr:lam7}, we conclude
    \begin{align}\label{Est:apr:lam8}
        \varphi_\lambda \in L^2(0,T;H^2(\Omega)).
    \end{align}
    Moreover, on account of \eqref{Est:apr:lam1}, \eqref{Est:apr:lam5}, \eqref{Est:apr:lam6}, a comparison argument in \eqref{Eq:lam1} yield as before
    \begin{align*}
        \bv_\lambda\in W^{1,4/3}(0,T;\mathbf{H}^1_\sigma(\Omega)^\prime).
    \end{align*}
    We have thus obtained uniform bounds with respect to $\lambda$ in all the above spaces.

    \subsection{Passage to the limit $\lambda\rightarrow0$} We can now pass to the limit along a suitable subsequence $\lambda_k\rightarrow0$, which still denote by $\lambda$ for ease of notation.
    By compactness, we immediately deduce that there exist $\varphi,\bv,\mu$, globally defined on $(0,\infty)$, such that, up to a subsequence, the following convergences hold
    \begin{align*}
        \varphi_{\lambda} &\rightharpoonup \varphi && \text{weakly-$*$ in }L^\infty(0,T;H^1(\Omega)),\\
        \varphi_{\lambda} &\rightharpoonup \varphi && \text{weakly in }L^2(0,T;H^2(\Omega)),\\
        \partial_t\varphi_{\lambda}&\rightharpoonup \partial_t\varphi && \text{weakly in }L^2(0,T;H^1(\Omega)^\prime),\\
        \varphi_{\lambda} &\rightarrow \varphi && \text{strongly in }L^2(0,T;H^1(\Omega))\cap C^0([0,T];L^2(\Omega)),\\
        \varphi_{\lambda} &\rightarrow \varphi && \text{a.e. in }\Omega_T,\\
        \mu_{\lambda} &\rightarrow \mu && \text{weakly in }L^2(0,T;L^2(\Omega)),\\
        \bv_{\lambda} &\rightharpoonup \bv && \text{weakly-$*$ in }L^\infty(0,T;\mathbf{L}^2_\sigma(\Omega)),\\
        \bv_{\lambda} &\rightharpoonup \bv && \text{weakly in }L^2(0,T;\mathbf{H}^1_\sigma(\Omega)),\\
        \partial_t\bv_{\lambda} &\rightharpoonup \partial_t\bv && \text{weakly in }L^{4/3}(0,T;\mathbf{H}^1_\sigma(\Omega)^\prime),\\
        \bv_{\lambda} &\rightarrow \bv && \text{strongly in }L^2(0,T;\mathbf{L}^2_\sigma(\Omega)),\\
        \bv_{\lambda} &\rightarrow \bv && \text{a.e. in }\Omega_T.
    \end{align*}

    %\begin{align}\label{Conv:f_lam}
    %    f_\lambda^\prime(\varphi_\lambda) \rightarrow f^\prime(\varphi) \quad \text{ in }L^2(0,T;L^2(\Omega)).
    %\end{align}
    By definition of $f_\lambda^\prime$, due to the weak-strong closure of maximal monotone operators (see, for instance, \cite[Lemma 2.3]{barbu-monot}) and the strong convergence for $(\varphi_\lambda)_{\lambda>0}$, we deduce
    \begin{align*}
        F_\lambda^\prime(\varphi_\lambda) \rightharpoonup F^\prime(\varphi)\quad \text{ in }L^2(0,T;L^2(\Omega)).
    \end{align*}
    Next, arguing as before (see Step \ref{Step:7}),
    %\begin{align*}
    %    \mu_\lambda\nabla\varphi_\lambda &\rightarrow \mu\nabla\varphi && \text{ in }L^1(\Omega_T), \\
    %  B(\bv_\lambda,\bv_\lambda) &\rightarrow B(\bv,\bv) && \text{ in }L^{4/3}(0,T;\mathbf{H}^1_\sigma(\Omega)^\prime),\\
    %   \bv_\lambda\cdot\nabla\varphi_\lambda &\rightarrow \bv\cdot\nabla\varphi  &&\text{ in }L^1(0,T;L^{3/2}(\Omega))
    %\end{align*}
    we can pass to the limit $\lambda\rightarrow 0^+$, along a suitable subsequence, in the weak formulation of \eqref{Eq:lam1}--\eqref{Eq:lam6}. Regarding the initial conditions, we first observe that $\varphi_\lambda(0) = \varphi_0$ almost everywhere  in $\Omega$ for all $\lambda>0$ and $\varphi_\lambda\rightarrow\varphi$ in $C^0([0,T];L^2(\Omega))$. Thus, using $\varphi_0 - \varphi(0)\in L^2(\Omega)$, we obtain
    \begin{align*}
        \int_\Omega\varphi_0(\varphi_0 - \varphi(0))\dx =  \int_\Omega\varphi_\lambda(0)(\varphi_0 - \varphi(0))\dx \rightarrow \int_\Omega\varphi(0)(\varphi_0 - \varphi(0))\dx,
    \end{align*}
    which implies $\varphi(0) = \varphi_0$ almost everywhere in $\Omega$. We know that there exists a subsequence such that $\bv_{\lambda}(0)\rightharpoonup\bv(0)$ in $\mathbf{H}^1_\sigma(\Omega)^\prime$. On the other hand, since $\bv_\lambda(0)=\bv_0$ in $\mathbf{H}^1_\sigma(\Omega)^\prime$, we deduce $\bv(0) = \bv_0$ in $\mathbf{H}^1_\sigma(\Omega)^\prime$. Thanks to the result in \cite[Lemma 4.1]{Abels09}, it holds
    \begin{align*}
        L^\infty(0,T;\mathbf{L}^2_\sigma(\Omega))\cap BUC([0,T];\mathbf{H}^1_\sigma(\Omega)^\prime) \hookrightarrow BC_{w}([0,T];\mathbf{L}^2_\sigma(\Omega)).
    \end{align*}
    Therefore, we have
    \begin{align*}
        \int_\Omega\bv(0)\cdot \bw \dx = \int_\Omega\bv_0\cdot \bw \dx
    \end{align*}
    for all $\bw\in\mathbf{H}^1_\sigma(\Omega)$. Since $\bv(0)-\bv_0 \in \mathbf{L}^2_\sigma(\Omega)$, a denseness argument allows us to test this equation by $\bw := \bv(0)-\bv_0$. In particular, this implies $\bv(0) = \bv_0$ almost everywhere in $\Omega$.

    Finally, we prove the energy inequality. Owing to the convergences above, we can use the lower-semicontinuity of norms and let $\lambda\rightarrow0^+$, along a suitable subsequence, in \eqref{energyineq}.
    Observe, in particular, that $F_\lambda(\varphi_\lambda) \geq F(J_\lambda\varphi_\lambda)$, where $J_\lambda := (I+\lambda F)^{-1}:\R\rightarrow\R$ denotes the resolvent of $F$. Then, from
    \begin{align*}
        |J_\lambda\varphi_\lambda - \varphi| \leq |J_\lambda\varphi_\lambda - \varphi_\lambda| + |\varphi_\lambda-\varphi| \leq |\lambda F_\lambda(\varphi_\lambda)|  + |\varphi_\lambda-\varphi|,
    \end{align*}
    we deduce that $J_\lambda\varphi_\lambda \rightarrow \varphi$ almost everywhere in $\Omega_T$. Thus, by lower-semicontinuity and Fatou's lemma, it follows that
    \begin{align*}
        \|F(\varphi(t))\|_{L^1(\Omega)} \leq \liminf_{\lambda\rightarrow0}\|F_\lambda(\varphi_\lambda(t))\|_{L^1(\Omega)},
    \end{align*}
    for almost every $t>0$.
    Finally, concerning the right-hand side of \eqref{energyineq}, note that $F_\lambda(\varphi(s)) \leq F(\varphi(s))$ for any $s>0$. Hence, using the convergence above, we can deduce from \eqref{energyineq} the energy inequality \eqref{ineq1} that is satisfied for almost any $t>0$. Nevertheless, by lower-semicontinuity of the norms involved it is immediate to see that it also holds for any $t>0$. This concludes the proof.
    
    \section{Proof of Theorem \ref{THM:StrongWP}}
    \label{strong}
    This section is split into three main parts. The first is devoted to the proof of local existence of a strong solution. The second to its uniqueness, whereas the third is dedicated to the validity of the strict separation property.
    \subsection{Existence of a local strong solution}  This will be performed again by means of a Galerkin \textit{ansatz}. To this end, we can take advantage of the Galerkin scheme as introduced in Step \ref{Step:2}. In particular, we consider the following discretized problem
    \begin{align}
        \partial_t\bv_{\lambda,m} + (\bv_{\lambda,m}\cdot\nabla)\bv_{\lambda,m} - \text{div}(\nu(\varphi_{\lambda,m})D\bv_{\lambda,m}) + \nabla \pi_{\lambda,m} &= \mu_{\lambda,m}\nabla\varphi_{\lambda,m}, \label{Eq:m1}\\
        \nabla\cdot\bv_{\lambda,m} &= 0, \label{Eq:m2}\\
        \partial_t\varphi_{\lambda,m} + \bv_{\lambda,m}\cdot\nabla\varphi_{\lambda,m} + \mu_{\lambda,m} - \overline{\mu_{\lambda,m}} &= 0, \label{Eq:m3}\\
        \mu_{\lambda,m} &= - \Delta\varphi_{\lambda,m} +P_m f^\prime_\lambda(\varphi_{\lambda,m}), \label{Eq:m4}
    \end{align}
    in $\Omega \times (0,T)$, subject to
    \begin{align}
        \bv_{\lambda,m} = \partial_\mathbf{n}\varphi_{\lambda,m} &= 0 &\quad&\text{ on }\partial\Omega \times (0,T), \label{Eq:m5}\\
        \bv_{\lambda,m}(0) = \bv_{0,m},\quad \varphi_{\lambda,m}(0) &= \varphi_{0,m} &\quad&\text{ in }\Omega. \label{Eq:m6}
    \end{align}
    All uniform estimates derived in the previous sections for the scheme above with $\lambda>0$, $m\in \N$, are clearly still valid, namely \eqref{EST:apr:v}-\eqref{H1p}.
    \subsubsection{Higher-order estimates for the Allen--Cahn equation}
    We test equation \eqref{Eq:m3} with $\partial_t(\mu_{\lambda,m}-\overline{\mu_{\lambda,m}})$. This yields
    \begin{align}
        \int_\Omega \partial_t\varphi_{\lambda,m}\partial_t\mu_{\lambda,m}\dx + \int_\Omega\bv_{\lambda,m}\cdot\nabla\varphi_{\lambda,m}\,\partial_t\mu_{\lambda,m}\dx + \frac{1}{2}\ddt\|\mu_{\lambda,m}-\overline{\mu_{\lambda,m}}\|_{L^2(\Omega)}^2 = 0.\label{basicA}
    \end{align}
    For the first term on the left-hand side, the definition of $\mu_{\lambda,m}$ (see \eqref{Eq:m4}) implies
    \begin{align*}
        \int_\Omega \partial_t\varphi_{\lambda,m}\partial_t\mu_{\lambda,m}\dx = \|\nabla\partial_t\varphi_{\lambda,m}\|_{\mathbf{L}^2(\Omega)}^2 + \int_\Omega F^{\prime\prime}_\lambda(\varphi_{\lambda,m})|\partial_t\varphi_{\lambda,m}|^2\dx-\int_\Omega \theta_0|\partial_t\varphi_{\lambda,m}|^2\dx.
    \end{align*}
    Note that (see property \eqref{Yos:6})
    \begin{align*}
        \int_\Omega F^{\prime\prime}_\lambda(\varphi_{\lambda,m})|\partial_t\varphi_{\lambda,m}|^2\dx \geq 0.
    \end{align*}
    Also, using standard interpolations and Young's inequality, we get
   \begin{align*}
        \int_\Omega \theta_0|\partial_t\varphi_{\lambda,m}|^2\dx\leq C\norm{\partial_t\varphi_{\lambda,m}}_{H^1(\Omega)'}\norm{\nabla\partial_t\varphi_{\lambda,m}}_{\mathbf{L}^2(\Omega)}\leq C\norm{\partial_t\varphi_{\lambda,m}}_{H^1(\Omega)'}^2+\frac{ 1}{12}\norm{\nabla \partial_t\varphi_{\lambda,m}}_{\mathbf{L}^2(\Omega)}^2.
     \end{align*}
    Arguing similarly as for \eqref{H1p}, we then observe by comparison and standard embeddings that
    \begin{align*}
        \norm{\partial_t\varphi_{\lambda,m}}_{H^1(\Omega)'}^2&\leq C(\norm{\mu_{\lambda,m}-\overline{\mu_{\lambda,m}}}_{\Ld}^2+\norm{\bv_{\lambda,m}}_{\mathbf L^4(\Omega)}^2\norm{\vphi_{\lambda,m}}_{L^4(\Omega)}^2)\nonumber\\& \leq C(\norm{\mu_{\lambda,m}-\overline{\mu_{\lambda,m}}}_{\Ld}^2+\norm{D\bv_{\lambda,m}}_{\mathbf L^2(\Omega)}^2),
    \end{align*}
   so that we get
     \begin{align*}
        \int_\Omega \theta_0|\partial_t\varphi_{\lambda,m}|^2\dx\leq C(\norm{\mu_{\lambda,m}-\overline{\mu_{\lambda,m}}}_{\Ld}^2+\norm{D\bv_{\lambda,m}}_{\mathbf L^2(\Omega)}^2)+\frac{ 1}{12}\norm{\nabla \partial_t\varphi_{\lambda,m}}_{\mathbf{L}^2(\Omega)}^2.
    \end{align*}
    Invoking the definition of $\mu_{\lambda,m}$ once again, we then infer that
    \begin{align*}
        \int_\Omega\bv_{\lambda,m}\cdot\nabla\varphi_{\lambda,m}\,\partial_t\mu_{\lambda,m}\dx &= \int_\Omega\bv_{\lambda,m}\cdot\nabla\varphi_{\lambda,m}\,\partial_t\Delta\varphi_{\lambda,m}\dx + \int_\Omega \bv_{\lambda,m}\cdot\nabla\varphi_{\lambda,m}\,(F^{\prime\prime}_\lambda(\varphi_{\lambda,m})-\theta_0)\partial_t\varphi_{\lambda,m}\dx \\
        &= -\int_\Omega \nabla\varphi_{\lambda,m}\cdot\nabla\bv_{\lambda,m}\nabla\partial_t\varphi_{\lambda,m}\dx - \int_\Omega \bv_{\lambda,m}\cdot D^2\varphi_{\lambda,m}\nabla\partial_t\varphi_{\lambda,m}\dx \\
        &\quad +\int_\Omega \bv_{\lambda,m}\cdot\nabla\varphi_{\lambda,m}\,f''_\lambda(\varphi_{\lambda,m})\partial_t\varphi_{\lambda,m}\dx \\
        &=: I_1 + I_2 + I_3.
    \end{align*}
    Observe first that, by Gagliardo--Nirenberg's inequalities,
    \begin{align}
        \vert I_1\vert = \left|\int_\Omega \nabla\varphi_{\lambda,m}\cdot \nabla\bv_{\lambda,m}\nabla\partial_t\varphi_{\lambda,m}\dx\right|
        &\nonumber\leq C\|\nabla \varphi_{\lambda,m}\|_{\mathbf L^6(\Omega)}\|\nabla\bv_{\lambda,m}\|_{\mathbf{L}^3(\Omega)}\|\nabla\partial_t\varphi_{\lambda,m}\|_{\mathbf{L}^2(\Omega)}\\&
        \leq C\|\varphi_{\lambda,m}\|_{ H^2(\Omega)}\|\nabla\bv_{\lambda,m}\|_{\mathbf{L}^2(\Omega)}^\frac12\|\bv_{\lambda,m}\|_{\mathbf{H}^2(\Omega)}^\frac12\|\nabla\partial_t\varphi_{\lambda,m}\|_{\mathbf{L}^2(\Omega)}\nonumber\\&
        \leq C(\|\varphi_{\lambda,m}\|_{ H^2(\Omega)}^4\|\nabla\bv_{\lambda,m}\|_{\mathbf{L}^2(\Omega)}^2)\nonumber\\
        &+\frac{\nu_*}{36}\|\bv_{\lambda,m}\|_{\mathbf{H}^2(\Omega)}^2+\frac1{12}\|\nabla\partial_t\varphi_{\lambda,m}\|_{\mathbf{L}^2(\Omega)}^2.
    \label{Est:Strong:1}
    \end{align}
    % By comparison in \eqref{Eq:m4}, we obtain {\color{blue}[Pay attention! This estimate can be obtained only after the passage to the limit as $m\to\infty$! In Galerkin there is the projector in front of $f'$...the estimates must be done in a different way in Galerkin, using the Lipschitz properties of $F_\lambda$, and making estimates depending on $\lambda$. Then after passin to the limit in $m$ all what follows is correct...]}
    Testing \eqref{Eq:m4} with $-\Delta\varphi_{\lambda,m}$, it holds
    \begin{align*}
        -\int_\Omega\mu_{\lambda,m}\Delta\varphi_{\lambda,m}\dx = \|\Delta\varphi_{\lambda,m}\|_{L^2(\Omega)}^2 - \int_\Omega P_m f_\lambda^\prime(\varphi_{\lambda,m})\Delta\varphi_{\lambda,m}\dx.
    \end{align*}
    Recalling \eqref{Est:F:Lipschitz} and using the H\"{o}lder and the Young inequalities, we infer that
    % \begin{align*}
    %     \|\varphi_m\|_{H^2(\Omega)} \leq C\big(\|f^\prime_\lambda(\varphi_m)\|_{L^2(\Omega)} + \|\mu_m\|_{L^2(\Omega)}\big),
    % \end{align*}
    % as well as
    % \begin{align*}
    %     \|f^\prime_\lambda(\varphi_m)\|_{L^2(\Omega)} \leq C\big(1+\|\mu_m\|_{L^2(\Omega)}\big)
    % \end{align*}
    % and thus,
    \begin{align}\label{Est:Comp:1}
        \|\varphi_{\lambda,m}\|_{H^2(\Omega)} \leq C_\lambda\big(1+\|\mu_{\lambda,m}\|_{L^2(\Omega)}\big),
    \end{align}
    where the constant $C_\lambda>0$ depends on $\lambda$ but is independent from $m$.  From now on, any constant appearing in the estimates will be independent of $m$.
    % Then, by comparison in \eqref{Eq:m2}, using Gagliardo--Nirenberg's inequalities, we have
    % \begin{align}
    %     \|\partial_t\varphi_{\lambda,m}\|_{L^2(\Omega)}
    %     &\leq \|\mu_{\lambda,m}-\overline{\mu_{\lambda,m}}\|_{L^2(\Omega)} + \|\bv_{\lambda,m}\cdot\nabla\varphi_{\lambda,m}\|_{L^2(\Omega)} \nonumber\\&\nonumber \|\mu_{\lambda,m}-\overline{\mu_{\lambda,m}}\|_{L^2(\Omega)} + \|\bv_{\lambda,m}\|_{\mathbf L^6(\Omega)}\|\nabla\varphi_{\lambda,m}\|_{\mathbf L^3(\Omega)}\\
    %     &\leq \|\mu_{\lambda,m}-\overline{\mu_{\lambda,m}}\|_{L^2(\Omega)} + C\|D\bv_{\lambda,m}\|_{\mathbf{L}^2(\Omega)}\|\varphi_{\lambda,m}\|_{H^2(\Omega)}^{\frac{1}{2}}\|\varphi_{\lambda,m}\|_{H^1(\Omega)}^{\frac{1}{2}} \nonumber\\
    %     &\leq \|\mu_{\lambda,m}-\overline{\mu_{\lambda,m}}\|_{L^2(\Omega)} + C_\lambda\|D\bv_{\lambda,m}\|_{\mathbf{L}^2(\Omega)}\big(1+\|\mu_{\lambda,m}\|_{L^2(\Omega)}\big)^{\frac{1}{2}}. \label{Est:Comp:2}
    % \end{align}
    Inserting \eqref{Est:Comp:1} into \eqref{Est:Strong:1}, we arrive at
    \begin{align*}
        |I_1|
        &\leq C_\lambda(1+\|\mu_{\lambda,m}\|_{ L^2(\Omega)}^4)\|D\bv_{\lambda,m}\|_{\mathbf{L}^2(\Omega)}^2+\frac{\nu_*}{36}\|A_S\bv_{\lambda,m}\|_{\mathbf{L}^2(\Omega)}^2+\frac1{12}\|\nabla\partial_t\varphi_{\lambda,m}\|_{\mathbf{L}^2(\Omega)}^2.
    \end{align*}
    Observe now that
    \begin{align*}
        |I_2|
        &\leq \|\varphi_{\lambda,m}\|_{H^2(\Omega)}\|\bv_{\lambda,m}\|_{\mathbf{L}^\infty(\Omega)}\|\nabla\partial_t\varphi_{\lambda,m}\|_{\mathbf{L}^2(\Omega)} \\
        &\leq C\|\varphi_{\lambda,m}\|_{H^2(\Omega)}\|\bv_{\lambda,m}\|_{\mathbf{H}^2(\Omega)}^{1/2}\|D\bv_{\lambda,m}\|_{\mathbf{L}^2(\Omega)}^{1/2}\|\nabla\partial_t\varphi_{\lambda,m}\|_{\mathbf{L}^2(\Omega)} \\
        &\leq \frac{1}{12}\|\nabla\partial_t\varphi_{\lambda,m}\|_{\mathbf{L}^2(\Omega)}^2 + \frac{\nu_*}{36}\|A_S\bv_{\lambda,m}\|_{\mathbf{L}^2(\Omega)}^2 +C_\lambda\|D\bv_{\lambda,m}\|_{\mathbf{L}^2(\Omega)}^2(1+\|\mu_{\lambda,m}\|_{L^2(\Omega)}^4),
    \end{align*}
    where we again used \eqref{Est:Comp:1} and Agmon's inequality.

    Consider $I_3$. An integration by parts gives
    \begin{align*}
        I_3=\int_\Omega\bv_{\lambda,m}\cdot\nabla\varphi_{\lambda,m} f_\lambda^{\prime\prime}(\varphi_{\lambda,m})\partial_t\varphi_{\lambda,m}\dx
        &= \int_\Omega\bv_{\lambda,m}\cdot\nabla f_\lambda^{\prime}(\varphi_{\lambda,m})\partial_t\varphi_{\lambda,m}\dx
       \\& = - \int_\Omega\bv_{\lambda,m}\cdot \nabla\partial_t\varphi_{\lambda,m} f^\prime_\lambda(\varphi_{\lambda,m})\dx.
    \end{align*}
    Thus, we infer that
    \begin{align}
        &\left|\int_\Omega\bv_{\lambda,m}\cdot\nabla\varphi_{\lambda,m} f_\lambda^{\prime\prime}(\varphi_{\lambda,m})\partial_t\varphi_{\lambda,m}\dx\right|
        \leq C\|\bv_{\lambda,m}\|_{\mathbf{L}^6(\Omega)}\|\nabla\partial_t\varphi_{\lambda,m}\|_{\mathbf{L}^2(\Omega)}\|f_\lambda^\prime(\varphi_{\lambda,m})\|_{L^3(\Omega)}\nonumber\\&\leq C\|D\bv_{\lambda,m}\|_{\mathbf{L}^2(\Omega)}\|\nabla\partial_t\varphi_{\lambda,m}\|_{\mathbf{L}^2(\Omega)}\|f^\prime_\lambda(\varphi_{\lambda,m})\|_{L^3(\Omega)}.\label{Est:Comp:3}
    \end{align}
    Using now the Lipschitz continuity of $f_\lambda^\prime$, we arrive at
    \begin{align}\label{Est:f_lambda}
        \|f^\prime_\lambda(\varphi_{\lambda,m})\|_{L^3(\Omega)} \leq C_\lambda(1+\|\varphi_\lambda\|_{L^3(\Omega)})\leq C_\lambda.
    \end{align}
    Hence, we conclude
    % \begin{align*}
    %     &\left|\int_\Omega\bv_{\lambda,m}\cdot\nabla\varphi_{\lambda,m} f_\lambda^{\prime\prime}(\varphi_{\lambda,m})\partial_t\varphi_{\lambda,m}\dx\right| \\
    %     &\leq C\|D\bv_{\lambda,m}\|_{\mathbf{L}^2(\Omega)}\|\nabla\partial_t\varphi_{\lambda,m}\|_{\mathbf{L}^2(\Omega)}\Big(1
    %     + \|\mu_{\lambda,m}\|_{L^2(\Omega)}^{1/2}C\Big(\|\nabla\partial_t\varphi_{\lambda,m}\|_{\mathbf{L}^2(\Omega)} +\|\bv_{\lambda,m}\|_{\mathbf{H}^2(\Omega)}\big(1+\|\mu_{\lambda,m}\|_{L^2(\Omega)}\big)^{1/2} \\
    %     &\qquad+C\|D\bv_{\lambda,m}\|_{\mathbf{L}^2(\Omega)}^{1/2}\|\bv_{\lambda,m}\|_{\mathbf{H}^2(\Omega)}^{1/2}\big(1+\|\mu_{\lambda,m}\|_{L^2(\Omega)}\big)\Big)^{1/2}\Big) \\
    %     &\leq \frac{\nu_*}{36}\|A_S\bv_{\lambda,m}\|_{\mathbf{L}^2(\Omega)}^2 + \frac{1}{12}\|\nabla\partial_t\varphi_{\lambda,m}\|_{\mathbf{L}^2(\Omega)}^2 + C\Big(1 + \|D\bv_{\lambda,m}\|_{\mathbf{L}^2(\Omega)}^6 + \|\mu_{\lambda,m}\|_{L^2(\Omega)}^6\Big).
    % \end{align*}
    \begin{align}\label{Est:I:3}
        \normmm{I_3}\leq \left|\int_\Omega\bv_{\lambda,m}\cdot\nabla\varphi_{\lambda,m} f_\lambda^{\prime\prime}(\varphi_{\lambda,m})\partial_t\varphi_{\lambda,m}\dx\right| \leq \frac{1}{12}\|\nabla\partial_t\varphi_{\lambda,m}\|_{\mathbf{L}^2(\Omega)}^2 + C\|D\bv_{\lambda,m}\|_{\mathbf{L}^2(\Omega)}^2.
    \end{align}
    % To estimate $I_3$ we are left to control one more term. By Gagliardo--Nirenberg's inequalities together with \eqref{Est:Comp:2}, we find
    % \begin{align*}
    %     &\normmm{\int_\Omega\bv_{\lambda,m}\cdot\nabla\varphi_{\lambda,m}\partial_t\varphi_{\lambda,m}\dx} \\
    %     &\leq C\|\bv_{\lambda,m}\|_{\mathbf{L}^6(\Omega)}\|\nabla\varphi_{\lambda,m}\|_{\mathbf{L}^2(\Omega)}\|\partial_t\varphi_{\lambda,m}\|_{L^3(\Omega)}\\&\leq C\|D\bv_{\lambda,m}\|_{\mathbf{L}^2(\Omega)}\|\nabla\varphi_{\lambda,m}\|_{\mathbf{L}^2(\Omega)}\|\partial_t\varphi_{\lambda,m}\|_{L^2(\Omega)}^{1/2}\|\nabla\partial_t\varphi_{\lambda,m}\|_{\mathbf{L}^2(\Omega)}^{1/2} \\
    %     &\leq C\|D\bv_{\lambda,m}\|_{\mathbf{L}^2(\Omega)}\Big(\|\mu_{\lambda,m}\|_{L^2(\Omega)}+\|D\bv_{\lambda,m}\|_{\mathbf{L}^2(\Omega)}\big(1+\|\mu_{\lambda,m}\|_{L^2(\Omega)}\big)^{\frac{1}{2}}\Big)^{1/2}\|\nabla\partial_t\varphi_{\lambda,m}\|_{\mathbf{L}^2(\Omega)}^{1/2} \\
    %     &\leq \frac{1}{12}\|\nabla\partial_t\varphi_{\lambda,m}\|_{\mathbf{L}^2(\Omega)}^2 + C\|D\bv_{\lambda,m}\|_{\mathbf{L}^2(\Omega)}^6 + C(1+\|\mu_{\lambda,m}\|_{L^2(\Omega)}^2).
    % \end{align*}
    Summing up, we have
    \begin{align}\label{Est:AC}
    \begin{split}
        &\frac{1}{2}\ddt\|\mu_{\lambda,m}-\overline{\mu_{\lambda,m}}\|_{L^2(\Omega)}^2  + \frac{3}{4}\|\nabla\partial_t\varphi_{\lambda,m}\|_{\mathbf{L}^2(\Omega)}^2 \\
        & \leq \frac{\nu_*}{18}\|A_S\bv_{\lambda,m}\|_{\mathbf{L}^2(\Omega)}^2 +C_\lambda\big(1+ \|D\bv_{\lambda,m}\|_{\mathbf{L}^2(\Omega)}^6 +\|\mu_{\lambda,m}\|_{L^2(\Omega)}^6\big).
    \end{split}
    \end{align}
    \subsubsection{Higher-order estimates for the Navier--Stokes equation} We test equation \eqref{Eq:m1} with $A_S\bv_{\lambda,m}$. This yields
    \begin{align}\label{Est:NS:1}
    \begin{split}
        &\frac{1}{2}\ddt\|D\bv_{\lambda,m}\|_{\mathbf{L}^2(\Omega)}^2 = \int_\Omega\mu_{\lambda,m}\nabla\varphi_{\lambda,m}\cdot A_S\bv_{\lambda,m}\dx \\
        &- \int_\Omega(\bv_{\lambda,m}\cdot\nabla)\bv_{\lambda,m}\cdot A_S\bv_{\lambda,m}\dx - \int_\Omega\nu(\varphi_{\lambda,m})D\bv_{\lambda,m}:DA_S\bv_{\lambda,m}\dx\\
        &= :I_4 + I_5 +I_6.
    \end{split}
    \end{align}
    To obtain an estimate on $\nabla \mu_{\lambda,m}$, we test \eqref{Eq:m3} by $-\Delta\mu_{\lambda,m}$ and integrate by parts. This yields
    \begin{align*}
       \int_\Omega\nabla\partial_t\varphi_{\lambda,m}\cdot\nabla\mu_{\lambda,m}\dx + \int_\Omega\nabla(\bv_{\lambda,m}\cdot\nabla\varphi_{\lambda,m})\nabla\mu_{\lambda,m}\dx = \|\nabla\mu_{\lambda,m}\|_{\mathbf{L}^2(\Omega)}^2.
    \end{align*}
    For the first term on the left-hand side, Hölder's inequality gives
    \begin{align*}
        \left|\int_\Omega\nabla\partial_t\varphi_{\lambda,m}\cdot\nabla\mu_{\lambda,m}\dx\right| \leq \|\nabla\partial_t\varphi_{\lambda,m}\|_{\mathbf{L}^2(\Omega)}\|\nabla\mu_{\lambda,m}\|_{\mathbf{L}^2(\Omega)}.
    \end{align*}
   Moreover, using \eqref{Est:Comp:1}, we have by Agmon's inequality and Sobolev embeddings
    \begin{align*}
        &\left|\int_\Omega\nabla(\bv_{\lambda,m}\cdot\nabla\varphi_{\lambda,m})\nabla\mu_{\lambda,m}\dx\right| \\
        &\leq \left|\int_\Omega \nabla\varphi_{\lambda,m}\cdot\nabla\bv_{\lambda,m}\nabla\mu_{\lambda,m}\dx\right| + \left|\int_\Omega \bv_{\lambda,m}\cdot D^2\varphi_{\lambda,m} \nabla\mu_{\lambda,m}\dx\right| \\
        &\leq C\|D\bv_{\lambda,m}\|_{\mathbf{L}^6(\Omega)}\|\nabla\varphi_{\lambda,m}\|_{\mathbf{L}^3(\Omega)}\|\nabla\mu_{\lambda,m}\|_{\mathbf{L}^2(\Omega)} + \|\bv_{\lambda,m}\|_{\mathbf{L}^\infty(\Omega)}\|D^2\varphi_{\lambda,m}\|_{\mathbf{L}^2(\Omega)}\|\nabla\mu_{\lambda,m}\|_{\mathbf{L}^2(\Omega)} \\
        &\leq C\|\bv_{\lambda,m}\|_{\mathbf{H}^2(\Omega)}\|\nabla\varphi_{\lambda,m}\|_{\mathbf{L}^2(\Omega)}^{1/2}\|\varphi_{\lambda,m}\|_{H^2(\Omega)}^{1/2}\|\nabla\mu_{\lambda,m}\|_{\mathbf{L}^2(\Omega)} \\
        &\qquad + C\|D\bv_{\lambda,m}\|_{\mathbf{L}^2(\Omega)}^{1/2}\|\bv_{\lambda,m}\|_{\mathbf{H}^2(\Omega)}^{1/2}\|\varphi_{\lambda,m}\|_{H^2(\Omega)}\|\nabla\mu_{\lambda,m}\|_{\mathbf{L}^2(\Omega)} \\
        &\leq C_\lambda\|\bv_{\lambda,m}\|_{\mathbf{H}^2(\Omega)}\Big(1+\|\mu_{\lambda,m}\|_{L^2(\Omega)}\Big)^{1/2}\|\nabla\mu_{\lambda,m}\|_{\mathbf{L}^2(\Omega)} \\
        &\qquad + C_\lambda\|D\bv_{\lambda,m}\|_{\mathbf{L}^2(\Omega)}^{1/2}\|\bv_{\lambda,m}\|_{\mathbf{H}^2(\Omega)}^{1/2}\Big(1+\|\mu_{\lambda,m}\|_{L^2(\Omega)}\Big)\|\nabla\mu_{\lambda,m}\|_{\mathbf{L}^2(\Omega)}.
    \end{align*}
    Combining these estimates, we then arrive at
    \begin{align}\label{Est:Nabla_mu}
    \begin{split}
        \|\nabla\mu_{\lambda,m}\|_{\mathbf{L}^2(\Omega)}
        &\leq C\|\nabla\partial_t\varphi_{\lambda,m}\|_{\mathbf{L}^2(\Omega)} + C_\lambda\|\bv_{\lambda,m}\|_{\mathbf{H}^2(\Omega)}\Big(1+\|\mu_{\lambda,m}\|_{L^2(\Omega)}\Big)^{1/2} \\
        &\qquad + C_\lambda\|D\bv_{\lambda,m}\|_{\mathbf{L}^2(\Omega)}^{1/2}\|\bv_{\lambda,m}\|_{\mathbf{H}^2(\Omega)}^{1/2}\Big(1+\|\mu_{\lambda,m}\|_{L^2(\Omega)}\Big).
    \end{split}
    \end{align}
Now, thanks to \eqref{Est:Comp:1} and \eqref{Est:Nabla_mu}, we obtain, by Gagliardo-Nirenberg's inequality,
    \begin{align}\label{Est:NS:2}
    \begin{split}
        \vert I_4\vert = \left|\int_\Omega\mu_{\lambda,m}\nabla\varphi_{\lambda,m}\cdot A_S\bv_{\lambda,m}\dx\right|
        &\leq \|A_S\bv_{\lambda,m}\|_{\mathbf{L}^2(\Omega)}\|\mu_{\lambda,m}\|_{L^3(\Omega)}\|\nabla\varphi_{\lambda,m}\|_{\mathbf{L}^6(\Omega)} \\
        &\leq C\|A_S\bv_{\lambda,m}\|_{\mathbf{L}^2(\Omega)}\|\mu_{\lambda,m}\|_{L^2(\Omega)}^{1/2}\|\nabla\mu_{\lambda,m}\|_{\mathbf{L}^2(\Omega)}^{1/2}\|\varphi_{\lambda,m}\|_{H^2(\Omega)} \\
        &\leq C_\lambda\|A_S\bv_{\lambda,m}\|_{\mathbf{L}^2(\Omega)}\|\mu_{\lambda,m}\|_{L^2(\Omega)}^{1/2}\big(1+\|\mu_{\lambda,m}\|_{L^2(\Omega)}\big) \\
        &\qquad \times\Big(\|\nabla\partial_t\varphi_{\lambda,m}\|_{\mathbf{L}^2(\Omega)} + \|\bv_{\lambda,m}\|_{\mathbf{H}^2(\Omega)}(1+\|\mu_{\lambda,m}\|_{L^2(\Omega)})^{1/2} \\
        &\qquad
        + \|D\bv_{\lambda,m}\|_{\mathbf{L}^2(\Omega)}^{1/2}\|\bv_{\lambda,m}\|_{\mathbf{H}^2(\Omega)}^{1/2}(1+\|\mu_{\lambda,m}\|_{L^2(\Omega)})\Big)^{1/2} \\
        &\leq \frac{\nu_*}{36}\|A_S\bv_{\lambda,m}\|_{\mathbf{L}^2(\Omega)}^2 + C_\lambda(1+\|\mu_{\lambda,m}\|_{L^2(\Omega)}^7) \\&\quad+ \omega_2\|\nabla\partial_t\varphi_{\lambda,m}\|_{\mathbf{L}^2(\Omega)}^2 + C_\lambda\|D\bv_{\lambda,m}\|_{\mathbf{L}^2(\Omega)}^6,
        \end{split}
    \end{align}
    for some $\omega_2>0$ to be chosen later on.
    Then, by Agmon's and Young's inequalities, we get
    \begin{align}\label{Est:NS:3}
    \begin{split}
        \vert I_5\vert = &\left|\int_\Omega(\bv_{\lambda,m}\cdot\nabla)\bv_{\lambda,m}\cdot A_S\bv_{\lambda,m}\dx\right| \\
        &\leq \norm{\bv_{\lambda,m}}_{\mathbf{L}^\infty(\Omega)}\norm{\nabla \bv_{\lambda,m}}_{\mathbf{L}^2(\Omega)}\norm{A_S\bv_{\lambda,m}}_{\mathbf{L}^2(\Omega)}\leq \frac{\nu_*}{36}\|A_S\bv_{\lambda,m}\|_{\mathbf{L}^2(\Omega)}^2 + C\|D\bv_{\lambda,m}\|_{\mathbf{L}^2(\Omega)}^6.
    \end{split}
    \end{align}
    Moreover, using integration by parts, we find
    \begin{align*}
        \vert I_6\vert &= \int_\Omega\nu(\varphi_{\lambda,m})D\bv_{\lambda,m}:DA_S\bv_{\lambda,m}\dx\\
        &=\int_\Omega\nu(\varphi_{\lambda,m})D\bv_{\lambda,m}:\nabla A_S\bv_{\lambda,m}\dx \\&= -\int_\Omega D\bv_{\lambda,m}\nabla\nu(\varphi_{\lambda,m})\cdot A_S\bv_{\lambda,m}\dx - \frac12\int_\Omega\nu(\varphi_{\lambda,m})\Delta\bv_{\lambda,m}\cdot A_S\bv_{\lambda,m}\dx \\
        &=: I_7 + I_8.
    \end{align*}
    By means of the chain rule, recalling that $\nu\in W^{1,\infty}(\R)$ and \eqref{Est:Comp:1}, $I_7$ can be estimated as follows
    \begin{align*}
        \left|\int_\Omega D\bv_{\lambda,m}\nabla\nu(\varphi_{\lambda,m})\cdot A_S\bv_{\lambda,m}\dx\right|
        &\leq
C\|D\bv_{\lambda,m}\|_{\mathbf{L}^3(\Omega)}\|A_S\bv_{\lambda,m}\|_{\mathbf{L}^2(\Omega)}\|\nabla\varphi_{\lambda,m}\|_{\mathbf L^6(\Omega)}
        \\&\leq C\|D\bv_{\lambda,m}\|_{\mathbf{L}^2(\Omega)}^{1/2}\|A_S\bv_{\lambda,m}\|_{\mathbf{L}^2(\Omega)}^{3/2}\|\varphi_{\lambda,m}\|_{H^2(\Omega)} \\
        &\leq C_\lambda\|D\bv_{\lambda,m}\|_{\mathbf{L}^2(\Omega)}^{1/2}\|A_S\bv_{\lambda,m}\|_{\mathbf{L}^2(\Omega)}^{3/2}\big(1+\|\mu_{\lambda,m}\|_{L^2(\Omega)}\big)\\
        &\leq \frac{\nu_*}{36}\|A_S\bv_{\lambda,m}\|_{\mathbf{L}^2(\Omega)}^2 + C_\lambda\big(1 + \|D\bv_{\lambda,m}\|_{\mathbf{L}^2(\Omega)}^6 + \|\mu_{\lambda,m}\|_{L^2(\Omega)}^6\big).
    \end{align*}
    The regularity theory for the Stokes operator entails that there exists some $\pi_{\lambda,m}\in C^0([0,T];H^1(\Omega))$ such that $-\Delta\bv_{\lambda,m} + \nabla \pi_{\lambda,m} = A_S\bv_{\lambda,m}$ holds almost everywhere in $\Omega$. Inserting this equation into $I_8$, we obtain
    \begin{align*}
        I_8 = \frac12\int_\Omega\nu(\varphi_{\lambda,m})\nabla \pi_{\lambda,m}\cdot A_S\bv_{\lambda,m}\dx - \frac12\int_\Omega\nu(\varphi_{\lambda,m})A_S\bv_{\lambda,m}\cdot A_S\bv_{\lambda,m}\dx =: I_9 + I_{10}.
    \end{align*}
    Recalling the interpolation result (see \cite[Lemma 3.1, Remark 3.2]{GGGP})
    \begin{align*}
        \norm{\pi_{\lambda,m}}_{L^3(\Omega)}\leq C\norm{A_S\bv_{\lambda,m}}^\frac12_{\mathbf{L}^2(\Omega)}\norm{D\bv_{\lambda,m}}^\frac12_{\mathbf{L}^2(\Omega)},
    \end{align*}
    integrating by parts, and using \eqref{Est:Comp:1}, we infer that
    \begin{align}\label{Est:NS:4}
    \begin{split}
       \vert I_9\vert&= \left|\frac12\int_\Omega\nu(\varphi_{\lambda,m})\nabla \pi_{\lambda,m}\cdot A_S\bv_{\lambda,m}\dx\right|\\
        &= \left|-\frac12\int_\Omega\nu^\prime(\varphi_{\lambda,m})\pi_{\lambda,m}\nabla \varphi_{\lambda,m} \cdot A_S\bv_{\lambda,m}\dx\right| \\
        &\leq C\|\nabla\varphi_{\lambda,m}\|_{\mathbf L^6(\Omega)}\|\pi_{\lambda,m}\|_{L^3(\Omega)}\|A_S\bv_{\lambda,m}\|_{\mathbf{L}^2(\Omega)}\\&\leq C\|\varphi_{\lambda,m}\|_{H^2(\Omega)}\|\pi_{\lambda,m}\|_{L^3(\Omega)}\|A_S\bv_{\lambda,m}\|_{\mathbf{L}^2(\Omega)} \\
        &\leq C_\lambda\big(1+\|\mu_{\lambda,m}\|_{L^2(\Omega)}\big)\|A_S\bv_{\lambda,m}\|_{\mathbf{L}^2(\Omega)}^{3/2}\|D\bv_{\lambda,m}\|_{\mathbf{L}^2(\Omega)}^{1/2} \\
        &\leq \frac{\nu_*}{36}\|A_S\bv_{\lambda,m}\|_{\mathbf{L}^2(\Omega)}^2 + C_\lambda\big(1+ \|D\bv_{\lambda,m}\|_{\mathbf{L}^2(\Omega)}^6 + \|\mu_{\lambda,m}\|_{L^2(\Omega)}^6\big).
    \end{split}
    \end{align}
    The term $I_{10}$ can be controlled from below. Indeed, we have (see \ref{ASS:Viscosity})
    \begin{align}\label{Est:NS:5}
        I_{10} \geq \frac{\nu_*}{2}\|A_S\bv_{\lambda,m}\|_{\mathbf{L}^2(\Omega)}^2.
    \end{align}
    In the next step, we test equation \eqref{Eq:m1} with $\partial_t\bv_{\lambda,m}$. This yields
    \begin{align}\label{Est:NS:6}
        \nonumber\|\partial_t\bv_{\lambda,m}\|_{\mathbf{L}^2(\Omega)}^2 &= -\int_\Omega(\bv_{\lambda,m}\cdot\nabla)\bv_{\lambda,m}\cdot\partial_t\bv_{\lambda,m}\dx \\
        &+ \int_\Omega\mu_{\lambda,m}\nabla\varphi_{\lambda,m}\cdot\partial_t\bv_{\lambda,m}\dx - \int_\Omega\nu(\varphi_{\lambda,m})D\bv_{\lambda,m}:D\partial_t\bv_{\lambda,m}\dx \nonumber\\
        &=: I_{11} + I_{12} + I_{13}.
    \end{align}
    Arguing as in \eqref{Est:NS:3}, we get
    \begin{align}\label{Est:NS:7}
        \vert I_{11}\vert&=\left|\int_\Omega(\bv_{\lambda,m}\cdot\nabla)\bv_{\lambda,m}\cdot\partial_t\bv_{\lambda,m}\dx\right| \\
        &\leq \frac{1}{12}\|\partial_t\bv_{\lambda,m}\|_{\mathbf{L}^2(\Omega)}^2 + \omega\|A_S\bv_{\lambda,m}\|_{\mathbf{L}^2(\Omega)}^2 + C(\omega)\|D\bv_{\lambda,m}\|_{\mathbf{L}^2(\Omega)}^6,
    \end{align}
    for some $\omega>0$ to be chosen later on.
    Consider now $I_{13}$. We first integrate by parts and then, recalling  \eqref{Est:Comp:1} and using Gagliardo--Nirenberg's as well as H\"{o}lder's and Young's inequalities, we deduce
    \begin{align}\label{Est:NS:8}
    \begin{split}
    &\left|-\int_\Omega\nu(\varphi_{\lambda,m})D\bv_{\lambda,m}:D\partial_t\bv_{\lambda,m}\dx\right| \\
    &\leq  \left|\int_\Omega\nu^\prime(\varphi_{\lambda,m})D\bv_{\lambda,m}\nabla\varphi_{\lambda,m}\cdot\partial_t\bv_{\lambda,m}\dx\right| + \left|\int_\Omega\nu(\varphi_{\lambda,m})\Delta\bv_{\lambda,m}\cdot \partial_t\bv_{\lambda,m}\dx\right|\\&
    \leq C \norm{D\bv_{\lambda,m}}_{\mathbf{L}^3(\Omega)}\norm{\nabla \varphi_{\lambda,m}}_{\mathbf{L}^6(\Omega)}\norm{\partial_t\bv_{\lambda,m}}_{\mathbf{L}^2(\Omega)}+C\norm{A_S\bv_{\lambda,m}}\norm{\partial_t\bv_{\lambda,m}}_{\mathbf{L}^2(\Omega)}\\
    &\leq C\|D\bv_{\lambda,m}\|_{\mathbf{L}^2(\Omega)}^{1/2}\|A_S\bv_{\lambda,m}\|_{\mathbf{L}^2(\Omega)}^{1/2}\|\varphi_{\lambda,m}\|_{H^2(\Omega)}\|\partial_t\bv_{\lambda,m}\|_{\mathbf{L}^2(\Omega)}\\&\quad+ C_S\|A_S\bv_{\lambda,m}\|_{\mathbf{L}^2(\Omega)}^2 + \frac{1}{12}\|\partial_t\bv_{\lambda,m}\|_{\mathbf{L}^2(\Omega)}^2 \\
    &\leq C_\lambda\|D\bv_{\lambda,m}\|_{\mathbf{L}^2(\Omega)}^{1/2}\|A_S\bv_{\lambda,m}\|_{\mathbf{L}^2(\Omega)}^{1/2}\big(1+\|\mu_{\lambda,m}\|_{L^2(\Omega)}\big)\|\partial_t\bv_{\lambda,m}\|_{\mathbf{L}^2(\Omega)}\\&\quad  + C_S\|A_S\bv_{\lambda,m}\|_{\mathbf{L}^2(\Omega)}^2 + \frac{1}{12}\|\partial_t\bv_{\lambda,m}\|_{\mathbf{L}^2(\Omega)}^2 \\
    &\leq \frac{1}{6}\|\partial_t\bv_{\lambda,m}\|_{\mathbf{L}^2(\Omega)}^2 + (\omega+C_S)\|A_S\bv_{\lambda,m}\|_{\mathbf{L}^2(\Omega)}^2 + C_\lambda(\omega)\big(1 + \|D\bv_{\lambda,m}\|_{\mathbf{L}^2(\Omega)}^6 +\|\mu_{\lambda,m}\|_{L^2(\Omega)}^6\big),
    \end{split}
    \end{align}
    where $C_S>0$ is a constant, and we can assume $C_S>1$.
    Moreover, we have, using again Gagliardo--Nirenberg's inequality combined with \eqref{Est:Comp:1} and \eqref{Est:Nabla_mu},
    \begin{align}\label{Est:NS:9}
    \begin{split}
        \vert I_{12}\vert = \left|\int_\Omega\mu_{\lambda,m}\nabla\varphi_{\lambda,m}\cdot\partial_t\bv_{\lambda,m}\dx\right|
        &\leq C\|\mu_{\lambda,m}\|_{L^3(\Omega)}\|\varphi_{\lambda,m}\|_{ L^6(\Omega)}\|\partial_t\bv_{\lambda,m}\|_{\mathbf{L}^2(\Omega)}\\&\leq C\|\mu_{\lambda,m}\|_{L^2(\Omega)}^{1/2}\|\nabla\mu_{\lambda,m}\|_{\mathbf{L}^2(\Omega)}^{1/2}\|\varphi_{\lambda,m}\|_{H^2(\Omega)}\|\partial_t\bv_{\lambda,m}\|_{\mathbf{L}^2(\Omega)} \\
        &\leq C_\lambda\|\mu_{\lambda,m}\|_{L^2(\Omega)}^{1/2}\|\partial_t\bv_{\lambda,m}\|_{\mathbf{L}^2(\Omega)}\big(1+\|\mu_{\lambda,m}\|_{L^2(\Omega)}\big) \\
        &\qquad \times\Big(\|\nabla\partial_t\varphi_{\lambda,m}\|_{\mathbf{L}^2(\Omega)} + \|\bv_{\lambda,m}\|_{\mathbf{H}^2(\Omega)}(1+\|\mu_{\lambda,m}\|_{L^2(\Omega)})^{1/2} \\
        &\qquad
        + \|D\bv_{\lambda,m}\|_{\mathbf{L}^2(\Omega)}^{1/2}\|\bv_{\lambda,m}\|_{\mathbf{H}^2(\Omega)}^{1/2}(1+\|\mu_{\lambda,m}\|_{L^2(\Omega)})\Big)^{1/2} \\
        &\leq \frac{1}{12}\|\partial_t\bv_{\lambda,m}\|_{\mathbf{L}^2(\Omega)}^2 +\omega_2\|\nabla\partial_t\varphi_{\lambda,m}\|_{\mathbf{L}^2(\Omega)}^2 + \omega\|A_S\bv_{\lambda,m}\|_{\mathbf{L}^2(\Omega)}^2 \\
        &\qquad + C_\lambda\big(1 + \|D\bv_{\lambda,m}\|_{\mathbf{L}^2(\Omega)}^6 +\|\mu_{\lambda,m}\|_{L^2(\Omega)}^7\big).
    \end{split}
    \end{align}
    Thanks to the above bounds, we deduce from \eqref{Est:NS:6} that
    \begin{align}
        &\nonumber\frac23\norm{\partial_t\bv_{\lambda,m}}_{\mathbf{L}^2(\Omega)}^2\\&\leq (3\omega+C_S)\norm{A_S\bv_{\lambda,m}}^2_{\mathbf{L}^2(\Omega)}+C_\lambda(\omega,\omega_2)\big(1 + \|D\bv_{\lambda,m}\|_{\mathbf{L}^2(\Omega)}^6 +\|\mu_{\lambda,m}\|_{L^2(\Omega)}^7\big)+{\omega_2}\norm{\nabla \partial_t\vphi_{\lambda,m}}_{\mathbf{L}^2(\Omega)}^2.\label{ptv}
    \end{align}
    Thus, adding \eqref{Est:NS:1} to \eqref{ptv} multiplied by $\frac32\frac{\omega}{ C_S+3\omega}$, choosing $\omega=\frac{\nu_*}{18C_S}$, and using \eqref{Est:NS:2}--\eqref{Est:NS:9}, we end up with
\begin{align}\label{Est:NS:10}
    \begin{split}
          \frac{1}{2}\ddt\|D\bv_{\lambda,m}\|_{\mathbf{L}^2(\Omega)}^2 + \frac{\omega}{ C_S+3\omega}\|\partial_t\bv_{\lambda,m}\|_{\mathbf{L}^2(\Omega)}^2 + \frac{11\nu_*}{36 }\|A_S\bv_{\lambda,m}\|_{\mathbf{L}^2(\Omega)}^2 \\
          \leq C_\lambda\big(1 + \|D\bv_{\lambda,m}\|_{\mathbf{L}^2(\Omega)}^6 +\|\mu_{\lambda,m}\|_{L^2(\Omega)}^7\big) + \frac14\|\nabla\partial_t\varphi_{\lambda,m}\|_{\mathbf{L}^2(\Omega)}^2,
    \end{split}
    \end{align}
    where we chose $\omega_2>0$ such that $\omega_2\big(1+\frac{3\omega}{2(C_S+3\omega)}\big)=\frac14$.
    \subsubsection{Completion of the proof}
    Combining estimates \eqref{Est:AC} and \eqref{Est:NS:10}, we find
    \begin{align}
    \begin{split}
        &\frac{1}{2}\ddt\Big(\|\mu_{\lambda,m}-\overline{\mu_{\lambda,m}}\|_{L^2(\Omega)}^2 + \|D\bv_{\lambda,m}\|_{\mathbf{L}^2(\Omega)}^2\Big) \\
        &\qquad+ \frac{\omega}{4 (C_S+3\omega)}\|\partial_t\bv_{\lambda,m}\|_{\mathbf{L}^2(\Omega)}^2 + \frac{\nu_*}{4}\|A_S\bv_{\lambda,m}\|_{\mathbf{L}^2(\Omega)}^2 + \frac{1}{2}\|\nabla\partial_t\varphi_{\lambda,m}\|_{\mathbf{L}^2(\Omega)}^2 \\
        &\leq C_\lambda\big(1 +\|D\bv_{\lambda,m}\|_{\mathbf{L}^2(\Omega)}^6 +\|\mu_{\lambda,m}\|_{L^2(\Omega)}^7\big)\\&\leq C_\lambda\big(1 + \|D\bv_{\lambda,m}\|_{\mathbf{L}^2(\Omega)}^2 +\|\mu_{\lambda,m}-\overline{\mu_{\lambda,m}}\|_{L^2(\Omega)}^2\big)^\frac72,
    \end{split}
    \label{analog}
    \end{align}
    where we used the fact that $\overline{\mu_{\lambda,m}}\leq C_\lambda(1+\norm{\mu_{\lambda,m}}_{\Ld})$, thanks to the regularity of $F_\lambda$.
    The inequality can be written as
    \begin{align*}
        u(t) \leq u(0) + \int_0^t C_\lambda w(u(s))\d s,
    \end{align*}
    where
    $$
        u := \|\mu_{\lambda,m}-\overline{\mu_{\lambda,m}}\|_{L^2(\Omega)}^2 + \|D\bv_{\lambda,m}\|_{\mathbf{L}^2(\Omega)}^2, \quad w(s) := (1+s)^{7/2}.
    $$
    Notice that, for $\lambda>0$ sufficiently small (see \eqref{Eq:m4}), we have
    \begin{align}\label{Est:u_00}
        u(0) = \|\mu_{\lambda,m}(0)-\overline{\mu_{\lambda,m}(0)}\|_{L^2(\Omega)}^2 + \|D\bv_{\lambda,m}(0)\|_{\mathbf{L}^2(\Omega)}^2 \leq C_{A,\lambda},
    \end{align}
    where the constant $C_{A,\lambda}$ does not depend on $m$ but it may depend on $\lambda$. Indeed, one can argue as in the proof of \cite[Theorem 3.1]{GPCAC} and control the quantity $\|\mu_{\lambda,m}(0)-\overline{\mu_{\lambda,m}(0)}\|_{L^2(\Omega)}$ uniformly in $m$, for $\lambda$ sufficiently small.

    Moreover, we note that
    \begin{align*}
        \int_0^T C_\lambda\;\d s:= C_{B,\lambda}(T),
    \end{align*}
    where $C_{B,\lambda}$ does not depend on $m$ but still depends on $\lambda$. Observe that $C_{B,\lambda}$ fulfills $C_{B,\lambda}(T)\rightarrow 0$ as $T\rightarrow0$ and is increasing with respect to $T$.
    Let us set
    \begin{align}
        G(x) := \int_1^x \frac{1}{(1+y)^{7/2}}\;\d y = \frac{1}{5\cdot 2^{3/2}} - \frac{2}{5(x+1)^{5/2}},\quad x>1,\label{G}
    \end{align}
    and observe that
    \begin{align*}
        G^{-1}(y) = \Big(2^{-5/2} - \frac{5}{2}y\Big)^{-2/5} - 1,\quad 0 < y < \frac{1}{5\cdot2^{3/2}}.
    \end{align*}
    As $G^{-1}$ is increasing, thanks to Bihari's inequality (see, e.g., \cite[Lemma II.4.12]{BoyerFabrie}), we obtain
    \begin{align}
        \max_{t\in[0,T]}u(t) \leq G^{-1}(G(C_{A,\lambda}) + C_{B,\lambda}) \leq C_{M,\lambda}(T),\label{max0}
    \end{align}
    uniformly in $m$, where $T<{T^*_\lambda}$, with $T^*_\lambda>0$ such that
    \begin{align}
        1 + C_{A,\lambda} < \Big(\frac{2}{5C_{B,\lambda}(T^*_\lambda)}\Big)^{2/5}.\label{maxtime0}
    \end{align}
    Note that $T^*_\lambda$ does not depend on $m$ (see \eqref{Est:u_00}). Also, being $C_{B,\lambda}(T)\rightarrow 0$ as $T\rightarrow0$, there exists some $T^*_\lambda$ satisfying the inequality.
    In particular, we have thus proved that the following estimates
    \begin{align}
        &\|\bv_{\lambda,m}\|_{L^\infty(0,T;\mathbf{H}^1_\sigma(\Omega))} + \|\bv_{\lambda,m}\|_{H^1(0,T;\mathbf{L}^2_\sigma(\Omega))} + \|\bv_{\lambda,m}\|_{L^2(0,T;\mathbf{H}^2(\Omega))} \leq C_{1,\lambda}, \label{Est:m1}\\
        &\|\varphi_{\lambda,m}\|_{H^1(0,T;H^1(\Omega))} + \|\varphi_{\lambda,m}\|_{L^\infty(0,T;H^1(\Omega))} + \|\varphi_{\lambda,m}\|_{L^2(0,T;H^2(\Omega))} \leq C_{2,\lambda}, \label{Est:m2}\\
        &\|\mu_{\lambda,m}\|_{L^\infty(0,T;L^2(\Omega))} + \|\mu_{\lambda,m}\|_{L^2(0,T;H^1(\Omega))} \leq C_{3,\lam} \label{Est:m3},
    \end{align}
    hold for all $T\in (0,T^*_\lambda)$, where the constants $C_{1,\lam}, C_{2,\lam}, C_{3,\lam}$ are independent of $m$.
    Observe that \eqref{Est:Comp:1} and \eqref{Est:m3} imply that $\|\varphi_{\lambda,m}\|_{L^\infty(0,T;H^2(\Omega))}\leq {C_\lam}$.
    Thus, passing to the limit $m\rightarrow\infty$ along a suitable subsequence, we conclude that there exists $( \bv_\lambda, \vphi_\lambda, \mu_\lambda)$ satisfying  problem \eqref{Eq:lam1}--\eqref{Eq:lam6}, and enjoying the regularity
    \begin{align}
        &\|\bv_\lambda\|_{L^\infty(0,T;\mathbf{H}^1_\sigma(\Omega))} + \|\bv_\lambda\|_{H^1(0,T;\mathbf{L}^2_\sigma(\Omega))} + \|\bv_\lambda\|_{L^2(0,T;\mathbf{H}^2(\Omega))} \leq C_{1,\lam}, \\
        &\label{H1H1}\|\varphi_\lambda\|_{H^1(0,T;H^1(\Omega))} + \|\varphi_\lambda\|_{L^\infty(0,T;H^2(\Omega))} \leq C_{2,\lam}, \\
        &\|\mu_\lambda\|_{L^\infty(0,T;L^2(\Omega))} + \|\mu_\lambda\|_{L^2(0,T;H^1(\Omega))} \leq C_{3,\lam},\label{L2H2}
    \end{align}
    for all $T\in (0,T^*)$.
Notice also that, thanks to \eqref{H1H1}, by comparison in \eqref{Eq:lam4}, we infer (using for instance time increment quotients) that
\begin{align*}
    \norm{\pt \mu_\lambda}_{L^2(0,T;H^1(\Omega)')}\leq C_\lambda\norm{\pt\vphi_\lambda}_{L^2(0,T;H^1(\Omega))},
\end{align*}
for any $T\in(0,T^*)$, so that, recalling \eqref{L2H2}, we get (this can be obtained first using smooth functions and then arguing by density) 
\begin{align}
\int_0^T\normmm{\frac12\ddt \norm{\mu_\lambda-\overline{\mu_\lambda}}_{L^2(\Omega)}^2}\dt\leq \norm{\mu_\lambda}_{L^2(0,T;H^1(\Omega))}\norm{\pt\mu_\lambda}_{L^2(0,T;H^1(\Omega)')}\leq C_\lambda,
    \label{mul}
\end{align}
entailing that 
\begin{align}
\norm{\mu_\lambda-\overline{\mu_\lambda}}_{L^2(\Omega)}\in AC([0,T]),
    \label{mul1}
\end{align} for any $T>0$. Similarly, we notice that, since $\pt \bv_\lambda\in L^2(0,T;\mathbf L_\sigma^2(\Omega))$ and $\bv_\lambda\in L^2(0,T;\bH^2(\Omega))$, 
\begin{align}
\norm{\bv_\lambda}_{\mathbf H^1(\Omega)}\in AC([0,T]),\quad \forall T\in (0,T^*).
    \label{bv1}
\end{align}
We now need to prove estimates which are independent from $\lambda$. Thanks to \eqref{mul1}, we can choose $\partial_t(\mu_\lambda-\overline\mu_\lambda)$ as a test function in \eqref{Eq:lam3} to obtain again \eqref{basicA}. Analogously, thanks to \eqref{bv1}, we can also use $A_S\bv_\lambda$ and $\pt \bv_\lambda$  as test functions in \eqref{Eq:lam1}, obtaining \eqref{Est:NS:1} as well as \eqref{Est:NS:6}. Therefore, we can follow the same estimates above to obtain an analog of \eqref{analog}, up to some differences which we now highlight. In particular we first need to revise estimate \eqref{Est:Comp:1}. Observe that, using a comparison argument, we get
    \begin{align*}
\|\varphi_\lambda\|_{H^2(\Omega)} \leq C(\|f_\lambda^\prime(\varphi_\lambda)\|_{L^2(\Omega)} + \|\mu_\lambda\|_{L^2(\Omega)})
    \end{align*}
    as well as
    \begin{align}
        \|f_\lambda^\prime(\varphi_\lambda)\|_{L^2(\Omega)} \leq C(1+\|\mu_\lambda\|_{L^2(\Omega)}).\label{essential}
    \end{align}
    In particular, this implies
    \begin{align}
        \|\varphi_\lambda\|_{H^2(\Omega)} \leq C(1+\|\mu_\lambda\|_{L^2(\Omega)}),\label{cpp}
    \end{align}
    where the constant $C>0$ does not depend on $\lambda$. In order to perform the same arguments as before, it only remains to modify the estimate in \eqref{Est:f_lambda}, since the right-hand side depends on $\lambda$.
    {\color{black} Therefore, we need to control $\|F_\lambda^\prime(\varphi_{\lambda})\|_{L^3(\Omega)}$ uniformly in $\lambda$. To this end, we test equation \eqref{Eq:m4} with $|F^\prime_\lambda(\varphi_{\lambda})|F_\lambda^\prime(\varphi_{\lambda})$. This yields
    $$
    \int_\Omega\mu_{\lambda} |F^\prime_\lambda(\varphi_{\lambda})|F^\prime_\lambda(\varphi_{\lambda})\dx= \int_\Omega\nabla\varphi_{\lambda}\cdot\nabla(|F^\prime_\lambda(\varphi_{\lambda})|F_\lambda^\prime(\varphi_{\lambda}))\dx + \int_\Omega f_\lambda^\prime(\varphi_{\lambda})|F_\lambda^\prime(\varphi_{\lambda})|F_\lambda^\prime(\varphi_{\lambda})\dx.
    $$
    Using the chain rule and invoking the properties of $F_\lambda$ (see property \eqref{Yos:6}), we infer that
    \begin{align*}
        \int_\Omega\nabla\varphi_{\lambda}\cdot\nabla(|F_\lambda^\prime(\varphi_{\lambda})|F_\lambda^\prime(\varphi_{\lambda}))\dx \geq 0.
    \end{align*}
    Moreover, it follows that (see \ref{ASS:S1})
    $$
        \int_\Omega f^\prime_\lambda(\varphi_{\lambda})|F_\lambda^\prime(\varphi_{\lambda})|F_\lambda^\prime(\varphi_{\lambda})\dx
    = \|F_\lambda^\prime(\varphi_{\lambda})\|^3_{L^3(\Omega)} - \theta_0\int_\Omega \varphi_{\lambda}|F^\prime_\lambda(\varphi_{\lambda})|F^\prime_\lambda(\varphi_{\lambda})\dx.
    $$
    Hence, we have
    $$
        \int_\Omega f^\prime_\lambda(\varphi_{\lambda})|F_\lambda^\prime(\varphi_{\lambda})|F_\lambda^\prime(\varphi_{\lambda})\dx
    \geq \frac{1}{2} \|F^\prime_\lambda(\varphi_{\lambda})\|^3_{L^3(\Omega)} - C,
    $$
    where we used Young's inequality and the fact that $\varphi_{\lambda}$ is bounded in $L^\infty(0,T;H^1(\Omega))$ uniformly with respect to $\lambda$. On the other hand, applying Hölder's and Young's inequalities, we find
    \begin{align*}
        \left|\int_\Omega\mu_{\lambda} |F_\lambda^\prime(\varphi_{\lambda})|F_\lambda^\prime(\varphi_{\lambda})\dx\right| \leq C\|\mu_{\lambda}\|_{L^3(\Omega)}^3 + \frac{1}{4}\|F_\lambda^\prime(\varphi_{\lambda})\|_{L^3(\Omega)}^3.
    \end{align*}
    Combining all these estimates, we then conclude, by Gagliardo-Nirenberg's inequalities, that
    \begin{align}\label{Est:F_prime1}
        \|F^\prime_\lambda
(\varphi_{\lambda})\|_{L^3(\Omega)} \leq C\Big(1+\|\mu_{\lambda}\|_{L^3(\Omega)}\Big) \leq C\Big(1+\|\mu_{\lambda}\|_{L^2(\Omega)}^{1/2}\|\mu_{\lambda}\|_{H^1(\Omega)}^{1/2}\Big).
    \end{align}
    }
    Then we need to control $\|\mu_{\lambda}\|_{H^1(\Omega)}$ in the next step. To this end, arguing as to obtain \eqref{Est:Nabla_mu}, we get
    \begin{align*}
\int_\Omega\nabla\partial_t\varphi_{\lambda}\cdot\nabla\mu_{\lambda}\dx + \int_\Omega\nabla(\bv_{\lambda}\cdot\nabla\varphi_{\lambda})\nabla\mu_{\lambda}\dx = \|\nabla\mu_{\lambda}\|_{\mathbf{L}^2(\Omega)}^2.
    \end{align*}
   Again, for the first term on the left-hand side, Hölder's inequality gives
    \begin{align*}
        \left|\int_\Omega\nabla\partial_t\varphi_{\lambda}\cdot\nabla\mu_{\lambda}\dx\right| \leq \|\nabla\partial_t\varphi_{\lambda}\|_{\mathbf{L}^2(\Omega)}\|\nabla\mu_{\lambda}\|_{\mathbf{L}^2(\Omega)}.
    \end{align*}
    For the second term, we use the chain rule to obtain from \eqref{cpp} (in place of \eqref{Est:Comp:1}),
    \begin{align*}
        &\left|\int_\Omega\nabla(\bv_{\lambda}\cdot\nabla\varphi_{\lambda})\nabla\mu_{\lambda}\dx\right| \\
        &\leq \left|\int_\Omega \nabla\varphi_{\lambda}\cdot\nabla\bv_{\lambda}\nabla\mu_{\lambda}\dx\right| + \left|\int_\Omega \bv_{\lambda}\cdot D^2\varphi_{\lambda} \nabla\mu_{\lambda}\dx\right| \\
        &\leq C\|D\bv_{\lambda}\|_{\mathbf{L}^6(\Omega)}\|\nabla\varphi_{\lambda}\|_{\mathbf{L}^3(\Omega)}\|\nabla\mu_{\lambda}\|_{\mathbf{L}^2(\Omega)} + \|\bv_{\lambda}\|_{\mathbf{L}^\infty(\Omega)}\|D^2\varphi_{\lambda}\|_{\mathbf{L}^2(\Omega)}\|\nabla\mu_{\lambda}\|_{\mathbf{L}^2(\Omega)} \\
        &\leq C\|\bv_{\lambda}\|_{\mathbf{H}^2(\Omega)}\|\nabla\varphi_{\lambda}\|_{\mathbf{L}^2(\Omega)}^{1/2}\|\varphi_{\lambda}\|_{H^2(\Omega)}^{1/2}\|\nabla\mu_{\lambda}\|_{\mathbf{L}^2(\Omega)} \\
        &\qquad + C\|D\bv_{\lambda}\|_{\mathbf{L}^2(\Omega)}^{1/2}\|\bv_{\lambda}\|_{\mathbf{H}^2(\Omega)}^{1/2}\|\varphi_{\lambda}\|_{H^2(\Omega)}\|\nabla\mu_{\lambda}\|_{\mathbf{L}^2(\Omega)} \\
        &\leq C\|\bv_{\lambda}\|_{\mathbf{H}^2(\Omega)}\Big(1+\|\mu_{\lambda}\|_{L^2(\Omega)}\Big)^{1/2}\|\nabla\mu_{\lambda}\|_{\mathbf{L}^2(\Omega)} \\
        &\qquad + C\|D\bv_{\lambda}\|_{\mathbf{L}^2(\Omega)}^{1/2}\|\bv_{\lambda}\|_{\mathbf{H}^2(\Omega)}^{1/2}\Big(1+\|\mu_{\lambda}\|_{L^2(\Omega)}\Big)\|\nabla\mu_{\lambda}\|_{\mathbf{L}^2(\Omega)},
    \end{align*}
    where this time all the constants do not depend on $\lambda$.
    Combining these estimates, we finally obtain 
\begin{align}\label{Est:Nabla_mu1}
    \begin{split}
        \|\nabla\mu_{\lambda}\|_{\mathbf{L}^2(\Omega)}
        &\leq C\|\nabla\partial_t\varphi_{\lambda}\|_{\mathbf{L}^2(\Omega)} + C\|\bv_{\lambda}\|_{\mathbf{H}^2(\Omega)}\Big(1+\|\mu_{\lambda}\|_{L^2(\Omega)}\Big)^{1/2} \\
        &\qquad + C\|D\bv_{\lambda}\|_{\mathbf{L}^2(\Omega)}^{1/2}\|\bv_{\lambda}\|_{\mathbf{H}^2(\Omega)}^{1/2}\Big(1+\|\mu_{\lambda}\|_{L^2(\Omega)}\Big).
    \end{split}
    \end{align}
    Hence, inserting \eqref{Est:Nabla_mu1} into \eqref{Est:F_prime1}, we obtain
    \begin{align}\label{Est:F_prime2}
    \begin{split}
        &\|f_\lambda^\prime(\varphi_{\lambda})\|_{L^3(\Omega)} \leq C(1+\|F^\prime_\lambda(\varphi_{\lambda})\|_{L^3(\Omega)})\\& \leq C\Big(1
        + C\|\mu_{\lambda}\|_{L^2(\Omega)}^{1/2}\Big(\|\nabla\partial_t\varphi_{\lambda}\|_{\mathbf{L}^2(\Omega)} +\|\bv_{\lambda}\|_{\mathbf{H}^2(\Omega)}\big(1+\|\mu_{\lambda}\|_{L^2(\Omega)}\big)^{1/2} \\
        &+C\|D\bv_{\lambda}\|_{\mathbf{L}^2(\Omega)}^{1/2}\|\bv_{\lambda}\|_{\mathbf{H}^2(\Omega)}^{1/2}\big(1+\|\mu_{\lambda}\|_{L^2(\Omega)}\big)\Big)^{1/2}\Big).
        \end{split}
    \end{align}
    Now, analogously to \eqref{Est:I:3}, we have
    \begin{align*}
        &\left|\int_\Omega\bv_{\lambda}\cdot\nabla\varphi_{\lambda} f_\lambda^{\prime\prime}(\varphi_{\lambda})\partial_t\varphi_{\lambda}\dx\right| \\
        &\leq C\|D\bv_{\lambda}\|_{\mathbf{L}^2(\Omega)}\|\nabla\partial_t\varphi_{\lambda}\|_{\mathbf{L}^2(\Omega)}\Big(1
        + \|\mu_{\lambda}\|_{L^2(\Omega)}^{1/2}C\Big(\|\nabla\partial_t\varphi_{\lambda}\|_{\mathbf{L}^2(\Omega)} +\|\bv_{\lambda}\|_{\mathbf{H}^2(\Omega)}\big(1+\|\mu_{\lambda}\|_{L^2(\Omega)}\big)^{1/2} \\
        &\qquad+C\|D\bv_{\lambda}\|_{\mathbf{L}^2(\Omega)}^{1/2}\|\bv_{\lambda}\|_{\mathbf{H}^2(\Omega)}^{1/2}\big(1+\|\mu_{\lambda}\|_{L^2(\Omega)}\big)\Big)^{1/2}\Big) \\
        &\leq \frac{\nu_*}{36}\|A_S\bv_{\lambda}\|_{\mathbf{L}^2(\Omega)}^2 + \frac{1}{12}\|\nabla\partial_t\varphi_{\lambda}\|_{\mathbf{L}^2(\Omega)}^2 + C\Big(1 + \|D\bv_{\lambda}\|_{\mathbf{L}^2(\Omega)}^6 + \|\mu_{\lambda}\|_{L^2(\Omega)}^6\Big).
    \end{align*}
    Using these results, we can finally argue as to obtain \eqref{analog}, but this time all bounds are independent of $\lambda$. This gives
    \begin{align}
    \begin{split}
        &\frac{1}{2}\ddt\Big(\|\mu_{\lambda}-\overline{\mu_{\lambda}}\|_{L^2(\Omega)}^2 + \|D\bv_{\lambda}\|_{\mathbf{L}^2(\Omega)}^2\Big) \\
        &\qquad+ \frac{\omega}{4 (C_S+3\omega)}\|\partial_t\bv_{\lambda}\|_{\mathbf{L}^2(\Omega)}^2 + \frac{\nu_*}{4}\|A_S\bv_{\lambda}\|_{\mathbf{L}^2(\Omega)}^2 + \frac{1}{4}\|\nabla\partial_t\varphi_{\lambda}\|_{\mathbf{L}^2(\Omega)}^2\\
        &\leq C\big(1 + \|D\bv_{\lambda}\|_{\mathbf{L}^2(\Omega)}^6 +\|\mu_{\lambda}\|_{L^2(\Omega)}^7\big)\\&
        \leq 
        C\big(1 + \|D\bv_{\lambda}\|_{\mathbf{L}^2(\Omega)} ^2+\|\mu_{\lambda}-\overline{\mu_\lambda}\|_{L^2(\Omega)}^2\big)^\frac72,
    \end{split}
    \end{align}
    where in the last inequality we used \eqref{essential} to control $\overline{\mu_\lambda}$.
    In particular, we also have
    \begin{align*}
        u(t) \leq u(0) + \int_0^t C_Rw(u(s))\d s,
    \end{align*}
   for $C_R>0$, where
    \begin{align*}
        &u := \|\mu_{\lambda}-\overline{\mu_{\lambda}}\|_{L^2(\Omega)}^2 + \|D\bv_{\lambda}\|_{\mathbf{L}^2(\Omega)}^2, \\
        &w(s) := (1+s)^{7/2}.
    \end{align*}
    Notice that by the assumptions on the initial data as well as \eqref{Eq:m4}, it holds, for $\lambda>0$ sufficiently small,
    \begin{align}\label{Est:u_0}
        u(0) = \|\mu_{\lambda}(0)-\overline{\mu_{\lambda}(0)}\|_{L^2(\Omega)}^2 + \|D\bv_\lambda(0)\|_{\mathbf{L}^2(\Omega)}^2 \leq C_{A},
    \end{align}
    where the constant $C_{A}$ is independent of $\lambda$ thanks to the properties of the approximation $f_\lambda$.
    Moreover, we note that
    \begin{align*}
        \int_0^T C_R\;\d s:= C_{B}(T)
    \end{align*}
   does not depend on $\lambda$. Observe that $C_B$ fulfills $C_B(T)\rightarrow 0$ as $T\rightarrow0$, but is also increasing with $T$. A careful inspection of the estimates above reveals that $C_B$ only depends on $E(\bv_0,\varphi_0)$, i.e., on the initial energy as well as on the parameters of the problem. Recalling the definition of  $G$ (see \eqref{G}), we deduce again by Bihari's inequality, cf. \cite[Lemma II.4.12]{BoyerFabrie}, that
    \begin{align}
        \max_{t\in[0,T]}u(t) \leq G^{-1}(G(C_{A}) + C_{B}) \leq C_{M}(T),\label{max}
    \end{align}
    uniformly in $\lambda$, where $T<{T^*}$, with $T^*>0$ such that
    \begin{align}
        1 + C_{A} < \Big(\frac{2}{5C_{B}(T^*)}\Big)^{2/5}.\label{maxtime}
    \end{align}
    Observe also that $T^*$ does not depend on $\lambda$ thanks to \eqref{Est:u_0}, and, thanks to $C_{B}(T)\rightarrow 0$ as $T\rightarrow0$, there exists some $T^*$ satisfying the constraint. Without loss of generality, we can assume that there exists a sequence $(\bv_\lambda, \vphi_\lambda)$ which is defined on $[0,T^*)$. Indeed, in case $T^*_\lambda<T^*$ for some $\lambda>0$, i.e., the maximal time of existence of each solution for fixed $\lambda>0$ is less than $T^*$, since $T^*$ is uniform in $\lambda$ and since each solution is unique (see the proof in Section \ref{unique}, which also holds with $f_\lambda$), one can extend any solution on $[0,T^*)$ by means of a continuation argument. We thus have obtained the following uniform estimates

    \begin{align*}
        &\|\bv_\lambda\|_{L^\infty(0,T;\mathbf{H}^1_\sigma(\Omega))} + \|\bv_\lambda\|_{H^1(0,T;\mathbf{L}^2_\sigma(\Omega))} + \|\bv_\lambda\|_{L^2(0,T;\mathbf{H}^2(\Omega))} \leq C_1, \\
        &\|\varphi_\lambda\|_{H^1(0,T;H^1(\Omega))} + \|\varphi_\lambda\|_{L^\infty(0,T;H^2(\Omega))} \leq C_2, \\
        &\|\mu_\lambda\|_{L^\infty(0,T;L^2(\Omega))} + \|\mu_\lambda\|_{L^2(0,T;H^1(\Omega))} \leq C_3,
    \end{align*}
    for all $T\in (0,T^*)$, where the constants $C_1(T), C_2(T), C_3(T)$ are independent of $\lambda$. Thus, passing to the limit $\lambda\rightarrow0^+$, we obtain
    \begin{align*}
        &\|\bv\|_{L^\infty(0,T;\mathbf{H}^1_\sigma(\Omega))} + \|\bv\|_{H^1(0,T;\mathbf{L}^2_\sigma(\Omega))} + \|\bv\|_{L^2(0,T;\mathbf{H}^2(\Omega))} \leq C_1, \\
        &\|\varphi\|_{H^1(0,T;H^1(\Omega))} + \|\varphi\|_{L^\infty(0,T;H^2(\Omega))} \leq C_2, \\
        &\|\mu\|_{L^\infty(0,T;L^2(\Omega))} + \|\mu\|_{L^2(0,T;H^1(\Omega))} \leq C_3,
    \end{align*}
    for all $T\in (0,T^*)$. The argument to show that $(\bv,\varphi,\mu)$ is a local strong solution is now standard. This concludes the proof of the existence of a strong solution to our problem.

    \subsection{Uniqueness of strong solutions}\label{unique}
    Consider two sets of initial data $(\bv_{0,1},\varphi_{0,1})$ and $(\bv_{0,2},\varphi_{0,2})$ satisfying the assumptions of Theorem \ref{THM:StrongWP}. For $i=1,2$, let $(\bv_i,\varphi_i,\mu_i)$ denote a strong solution corresponding to the initial data $(\bv_{0,i},\varphi_{0,i})$, respectively. We further write
    \begin{align*}
        (\bv_0,\varphi_0) &:= (\bv_{0,1},\varphi_{0,1}) -(\bv_{0,2},\varphi_{0,2}), \\
        (\bv,\varphi,\mu) &:= (\bv_1,\varphi_1,\mu_1)-(\bv_2,\varphi_2,\mu_2).
    \end{align*}
    This means that the triple $(\bv,\varphi,\mu)$ fulfills the following system of equations in a strong sense:
    \begin{align}
        \partial_t\bv + (\bv_1\cdot\nabla)\bv_1 - (\bv_2\cdot\nabla)\bv_2 - \text{div}(\nu(\varphi_1)D\bv_1) \nonumber\\
        +  \text{div}(\nu(\varphi_2)D\bv_2) + \nabla \pi &= \mu_1\nabla\varphi_1 -\mu_2\nabla\varphi_2 &\quad&\text{ in }\Omega_T, \label{Eq:1}\\
        \nabla\cdot\bv &= 0 &\quad&\text{ in }\Omega_T, \label{Eq:2}\\
        \partial_t\varphi + \bv_1\cdot\nabla\varphi_1 - \bv_2\cdot\nabla\varphi_2 + \mu - \overline{\mu} &= 0 &\quad&\text{ in }\Omega_T, \label{Eq:3}\\
        \bv = \partial_\mathbf{n}\varphi &= 0 &\quad&\text{ on }\partial\Omega_T, \label{Eq:4}\\
        \bv(0) = \bv_0,\quad \varphi(0) &= \varphi_0 &\quad&\text{ in }\Omega. \label{Eq:5}
    \end{align}
    We now test \eqref{Eq:1} by $A_S^{-1}\bv$ and \eqref{Eq:3} by $\varphi$, and add the resulting equations. Integrating by parts and using the identity
    \begin{align*}
        \frac{1}{2}\ddt\|\bv\|_{\mathbf{H}^1_\sigma(\Omega)^\prime}^2 = \frac{1}{2}\ddt\|\nabla A_S^{-1}\bv\|_{\mathbf{L}^2(\Omega)}^2 = \langle\partial_t\bv, A_S^{-1}\bv\rangle_{\mathbf{H}^1_\sigma(\Omega)},
    \end{align*}
    we then infer that
    \begin{align}\label{Est:U:1}
    \begin{split}
        &\frac{1}{2}\ddt\Big(\|\bv\|_{\mathbf{H}^1_\sigma(\Omega)^\prime}^2 + \|\varphi\|_{L^2(\Omega)}^2\Big) = \int_\Omega\big(\mu_1\nabla\varphi_1 - \mu_2\nabla\varphi_2\big)\cdot A_S^{-1}\bv \dx\\
        &- \int_\Omega\big((\bv_1\cdot\nabla)\bv_1 - (\bv_2\cdot\nabla)\bv_2\big)\cdot A_S^{-1}\bv\dx -\int_\Omega\big(\nu(\varphi_1)D\bv_1 - \nu(\varphi_2)D\bv_2\big):DA_S^{-1}\bv\dx\\
        &- \int_\Omega(\bv_1\cdot\nabla\varphi_1 - \bv_2\cdot\nabla\varphi_2)\varphi\dx
        - \int_\Omega(\mu-\overline{\mu})\varphi\dx .
    \end{split}
    \end{align}
    Recalling that $\bv_i$ are divergence-free, integrating by parts, using Gagliardo-–Nirenberg's inequality and Young's inequality, we deduce, since $\bv_i\in L^\infty(0,T;\mathbf H^1_\sigma(\Omega))$,
    \begin{align}\label{Est:U:2}
    \begin{split}
        &\left|\int_\Omega\big((\bv_1\cdot\nabla)\bv_1 - (\bv_2\cdot\nabla)\bv_2\big)\cdot A_S^{-1}\bv\dx\right| \\
        &\leq  \int_\Omega|(\bv_1\otimes\bv):\nabla A_S^{-1}\bv|\dx +\int_\Omega|(\bv\otimes\bv_2):\nabla A_S^{-1}\bv|\dx \\
        &\leq \|\bv_1\|_{\mathbf{L}^6(\Omega)}\|\nabla A_S^{-1}\bv\|_{\mathbf{L}^3(\Omega)}\|\bv\|_{\mathbf{L}_\sigma^2(\Omega)} + \|\bv\|_{\mathbf{L}_\sigma^2(\Omega)}\|\nabla A_S^{-1}\bv\|_{\mathbf{L}^3(\Omega)}\|\bv_2\|_{\mathbf{L}^6(\Omega)} \\
        &
        \leq C\norm{\nabla A_S^{-1}\bv}_{\mathbf{L}^2(\Omega)}^\frac12\norm{\bv}^\frac32_{\mathbf{L}^6(\Omega)} \leq \frac{\nu_*}{28} \|\bv\|_{\mathbf{L}_\sigma^2(\Omega)}^2 + C\|\bv\|_{\mathbf{H}^1_\sigma(\Omega)^\prime}^2.
        \end{split}
        \end{align}
        Moreover, we observe that
        \begin{align*}
            \int_\Omega\big(\nu(\varphi_1)D\bv_1 - \nu(\varphi_2)D\bv_2\big):DA_S^{-1}\bv\dx
            &= \int_\Omega\big(\nu(\varphi_1)-\nu(\varphi_2)\big)D\bv_1:DA_S^{-1}\bv\dx \\
            &+ \int_\Omega\nu(\varphi_2)D\bv:DA_S^{-1}\bv\dx.
        \end{align*}
        Invoking the Lipschitz continuity of $\nu$, we infer
        \begin{align}\label{Est:U:3}
        \begin{split}
            &\left|\int_\Omega\big(\nu(\varphi_1)-\nu(\varphi_2)\big)D\bv_1:DA_S^{-1}\bv\dx\right|\\
            &\leq C\|\varphi\|_{L^2(\Omega)}\|D\bv_1\|_{\mathbf{L}^6(\Omega)}\|DA_S^{-1}\bv\|_{\mathbf{L}^3(\Omega)} \\
            &\leq C\|\varphi\|_{L^2(\Omega)}\|\bv_1\|_{\mathbf{H}^2(\Omega)}\|DA_S^{-1}\bv\|_{\mathbf{L}^2(\Omega)}^{1/2}\|\bv\|_{\mathbf{L}_\sigma^2(\Omega)}^{1/2} \\
            &\leq \frac{\nu_*}{28}\|\bv\|_{\mathbf{L}_\sigma^2(\Omega)}^2 +C\big(\|\bv_1\|_{\mathbf{H}^2(\Omega)}^2\|\varphi\|_{L^2(\Omega)}^2 + \|\bv\|_{\mathbf{H}^1_\sigma(\Omega)^\prime}^2\big).
        \end{split}
        \end{align}
        In the remaining term, an integration by parts gives
        \begin{align*}
            \int_\Omega\nu(\varphi_2)D\bv:DA_S^{-1}\bv\dx = -\int_\Omega\nu^\prime(\varphi_2)D A_S^{-1}\bv\nabla\varphi_2\cdot \,\bv\dx - \frac12\int_\Omega\nu(\varphi_2)\bv\cdot\Delta A_S^{-1}\bv\dx.
        \end{align*}
        Since $\nu\in W^{1,\infty}(\R)$, we infer by Gagliardo-Nirenberg's inequalities, that
        \begin{align}\label{Est:U:4}
        \begin{split}
            \normmm{\int_\Omega\nu^\prime(\varphi_2)D A_S^{-1}\bv\nabla\varphi_2\cdot \,\bv\dx }
            &\leq C\|\nabla\varphi_2\|_{\mathbf{L}^6(\Omega)}\|\bv\|_{\mathbf{L}^2_\sigma(\Omega)}\|DA_S^{-1}\bv\|_{\mathbf{L}^3(\Omega)} \\&
            \leq C\|\bv\|_{\mathbf{L}^2_\sigma(\Omega)}^\frac32\|DA_S^{-1}\bv\|_{\mathbf{L}^2(\Omega)}^\frac12 \\
            &\leq \frac{\nu_*}{28}\|\bv\|_{\mathbf{L}_\sigma^2(\Omega)}^2 + C\|\bv\|_{\mathbf{H}^1_\sigma(\Omega)^\prime}^2,
        \end{split}
        \end{align}
        where we used $\vphi_i\in L^\infty(0,T;H^2(\Omega)),$ $i=1,2$.
        Now, by regularity theory for the Stokes operator, there exists some $\pi\in C^0([0,T];H^1(\Omega))$ such that $-\Delta A_S^{-1}\bv + \nabla \pi = \bv$ holds almost everywhere in $\Omega$. This yields
        \begin{align*}
            -\frac12\int_\Omega\nu(\varphi_2)\bv\cdot\Delta A_S^{-1}\bv\dx = \frac12\int_\Omega\nu(\varphi_2)|\bv|^2\dx - \frac12\int_\Omega\nu(\varphi_2)\bv\cdot\nabla \pi\dx.
        \end{align*}
    Using the properties of $\nu$, it follows that
    \begin{align}\label{Est:U:5}
        \frac12\int_\Omega\nu(\varphi_2)|\bv|^2\dx \geq \frac{\nu_*}{2}\|\bv\|_{\mathbf{L}_\sigma^2(\Omega)}^2.
    \end{align}
    Moreover, from \cite[Lemma 3.1, Remark 3.2]{GGGP}, we deduce
\begin{align*}
    \norm{\pi}_{L^3(\Omega)}\leq C\norm{\nabla A_S\bv}_{\mathbf L^2(\Omega)}^\frac12\norm{\bv}_{\mathbf L^2(\Omega)}^\frac12,
\end{align*}
so that, after an integration by parts, we obtain
    \begin{align}\label{Est:U:6}
    \begin{split}
        \left|-\frac12\int_\Omega\nu(\varphi_2)\bv\cdot\nabla \pi\dx\right| = \left|\frac12\int_\Omega\nu^\prime(\varphi_2)\nabla\varphi_2\cdot\bv \pi\dx\right|
        &\leq C\|\nabla\varphi_2\|_{\mathbf{L}^6(\Omega)}\|\bv\|_{\mathbf{L}^2_\sigma(\Omega)}\|\pi\|_{L^3(\Omega)} \\
        &\leq C\|\nabla\varphi_2\|_{\mathbf{L}^6(\Omega)}\|\bv\|_{\mathbf{L}^2_\sigma(\Omega)}^{3/2}\|\nabla A_S^{-1}\bv\|_{\mathbf{L}^2(\Omega)}^{1/2} \\
        &\leq \frac{\nu_*}{28}\|\bv\|_{\mathbf{L}_\sigma^2(\Omega)}^2 + C\|\bv\|_{\mathbf{H}^1_\sigma(\Omega)^\prime}^2.
    \end{split}
    \end{align}
    Regarding the Korteweg stress, we use the definition of $\mu_i$ to deduce
    \begin{align*}
        \int_\Omega\big(\mu_1\nabla\varphi_1 - \mu_2\nabla\varphi_2\big)\cdot A_S^{-1}\bv \dx &= -\int_\Omega\Delta\varphi_1\nabla\varphi_1\cdot A_S^{-1}\bv\dx + \int_\Omega f^\prime(\varphi_1)\nabla\varphi_1\cdot A_S^{-1}\bv\dx \\
        &+
        \int_\Omega\Delta\varphi_2\nabla\varphi_2\cdot A_S^{-1}\bv\dx - \int_\Omega f^\prime(\varphi_1)\nabla\varphi_1\cdot A_S^{-1}\bv\dx.
    \end{align*}
    Since $\bv_i$ is divergence-free, integration by parts entails
    \begin{align*}
        \int_\Omega f^\prime(\varphi_i)\nabla\varphi_i\cdot A_S^{-1}\bv\dx = \int_\Omega \nabla f(\varphi_i)\cdot A_S^{-1}\bv\dx = 0
    \end{align*}
    for $i=1,2$. Thus, it holds
    \begin{align*}
        \int_\Omega\big(\mu_1\nabla\varphi_1 - \mu_2\nabla\varphi_2\big)\cdot A_S^{-1}\bv \dx
        &= -\int_\Omega\Delta\varphi_1\nabla\varphi_1\cdot A_S^{-1}\bv\dx + \int_\Omega\Delta\varphi_2\nabla\varphi_2\cdot A_S^{-1}\bv\dx \\
        &= -\int_\Omega\Delta\varphi_1\nabla\varphi\cdot A_S^{-1}\bv\dx - \int_\Omega\Delta\varphi\nabla\varphi_2\cdot A_S^{-1}\bv\dx.
    \end{align*}
    For the first term on the right-hand side, we infer, by Agmon's inequality,
    \begin{align}\label{Est:U:7}
    \begin{split}
        \left|\int_\Omega\Delta\varphi_1\nabla\varphi\cdot A_S^{-1}\bv\dx\right|
        &\leq \|\varphi_1\|_{H^2(\Omega)}\|\nabla\varphi\|_{\mathbf{L}^2(\Omega)}\|A_S^{-1}\bv\|_{\mathbf{L}^\infty(\Omega)} \\
        &\leq C\|\nabla\varphi\|_{\mathbf{L}^2(\Omega)}\|\bv\|_{\mathbf{L}^2_\sigma(\Omega)}^{1/2}\|\bv\|_{\mathbf{H}^1_\sigma(\Omega)^\prime}^{1/2} \\
        &\leq \frac{1}{8}\|\nabla\varphi\|_{\mathbf{L}^2(\Omega)}^2 + \frac{\nu_*}{28}\|\bv\|_{\mathbf{L}_\sigma^2(\Omega)}^2 + C\|\bv\|_{\mathbf{H}^1_\sigma(\Omega)^\prime}^2.
    \end{split}
    \end{align}
    Moreover, integrating by parts, we can control the second term by means of Agmon's inequality as follows
    \begin{align}\label{Est:U:8}
    \begin{split}
        \left|-\int_\Omega\Delta\varphi\nabla\varphi_2\cdot A_S^{-1}\bv\dx\right|
        &= \int_\Omega\nabla\varphi\cdot D^2\varphi_2 A_S^{-1}\bv\dx + \int_\Omega\nabla\varphi\cdot DA_S^{-1}\bv\,\nabla\varphi_2\dx \\
        &\leq C\|\varphi_2\|_{H^2(\Omega)}\|\nabla\varphi\|_{\mathbf{L}^2(\Omega)}\|A_S^{-1}\bv\|_{\mathbf{L}^\infty(\Omega)} \\
        &\qquad + C\|\nabla\varphi\|_{\mathbf{L}^2(\Omega)}\|DA_S^{-1}\bv\|_{\mathbf{L}^3(\Omega)}\|\nabla\varphi_2\|_{\mathbf L^6(\Omega)}
        \\
        &\leq C\|\nabla\varphi\|_{\mathbf{L}^2(\Omega)}\|\nabla A_S^{-1}\bv\|_{\mathbf{L}^2(\Omega)}^\frac12\|\bv\|_{\mathbf{L}^2(\Omega)} ^\frac12
        \\
        &\leq \frac{1}{8}\|\nabla\varphi\|_{\mathbf{L}^2(\Omega)}^2 + \frac{\nu_*}{28}\|\bv\|_{\mathbf{L}_\sigma^2(\Omega)}^2 + C\|\bv\|_{\mathbf{H}^1_\sigma(\Omega)^\prime}^2.
    \end{split}
    \end{align}
    Next, we consider the convective term in the Allen--Cahn equation. There, it holds
    \begin{align*}
        \int_\Omega(\bv_1\cdot\nabla\varphi_1 - \bv_2\cdot\nabla\varphi_2)\varphi\dx = \int_\Omega(\bv\cdot\nabla\varphi_1)\varphi\dx + \int_\Omega(\bv_2\cdot\nabla\varphi)\varphi\dx.
    \end{align*}
    Since $\bv_2$ is divergence-free, integration by parts yields
    \begin{align*}
        \int_\Omega(\bv_2\cdot\nabla\varphi)\varphi\dx = 0.
    \end{align*}
    For the remaining term, we obtain
    \begin{align}\label{Est:U:9}
    \begin{split}
        \left|\int_\Omega(\bv\cdot\nabla\varphi_1)\varphi\dx\right|
        &\leq \|\bv\|_{\mathbf{L}^2_\sigma(\Omega)}\|\nabla\varphi_1\|_{\mathbf{L}^6(\Omega)}\|\varphi\|_{L^3(\Omega)} \\
        &\leq C\|\bv\|_{\mathbf{L}^2_\sigma(\Omega)}\|\varphi_1\|_{H^2(\Omega)}\|\varphi\|_{L^2(\Omega)}^{1/2}\|\varphi\|_{H^1(\Omega)}^{1/2} \\
        &\leq \frac{\nu_*}{28}\|\bv\|_{\mathbf{L}^2_\sigma(\Omega)}^2 + \frac{1}{8}\|\nabla\varphi\|_{\mathbf{L}^2(\Omega)}^2 + C\|\varphi\|_{L^2(\Omega)}^2.
    \end{split}
    \end{align}
    For the last term, we use that $\overline{\varphi}=0$. Thus, on account of Fubini's theorem and the definition of $\mu$, we infer that
    \begin{align*}
        \int_\Omega(\mu-\overline{\mu})\varphi\dx = \int_\Omega\mu\,\varphi\dx = \|\nabla\varphi\|_{\mathbf{L}^2(\Omega)}^2 + \int_\Omega\big(f^\prime(\varphi_1)-f^\prime(\varphi_2)\big)\varphi\dx.
    \end{align*}
    Recalling \ref{ASS:S1}, it follows that
    \begin{align}\int_\Omega\big(F^\prime(\varphi_1)-F^\prime(\varphi_2)\big)\varphi\dx \geq \theta\|\varphi\|_{L^2(\Omega)}^2,
    \end{align}
    as well as
    \begin{align}
        \label{Est:U:10}
        \normmm{\int_\Omega \frac{\theta_0}2\varphi^2\dx}\leq C\norm{\varphi}_{L^2(\Omega)}^2.
    \end{align}
    Eventually, combining the estimates \eqref{Est:U:1}--\eqref{Est:U:10}, we arrive at
    \begin{align*}
        \frac{1}{2}\ddt\Big(\|\bv\|_{\mathbf{H}^1_\sigma(\Omega)^\prime}^2 + \|\varphi\|_{L^2(\Omega)}^2\Big) + \frac{1}{2}\|\nabla\varphi\|_{\mathbf{L}^2(\Omega)}^2 + \frac{\nu_*}{2}\|\bv\|_{\mathbf{L}^2_\sigma(\Omega)}^2 \leq C\Big(\|\bv\|_{\mathbf{H}^1_\sigma(\Omega)^\prime}^2 +\norm{\bv_1}_{\mathbf H^2(\Omega)}^2 \|\varphi\|_{L^2(\Omega)}^2\Big).
    \end{align*}
     Thus, Gronwall's lemma gives \eqref{contdep}, recalling $\bv_1\in L^2(0,T;\bH^2(\Omega))$, and the proof is finished.
     \subsection{Strict separation property}
We start from the instantaneous strict separation property for strong solutions. The proof is the same as in \cite[Theorem 3.10]{GPCAC} and takes inspiration from the ideas developed in \cite{GGPC} concerning the control over $\overline\mu$. Nevertheless, since here we consider only two phases, for the sake of clarity we present the main steps of the proof by means of De Giorgi iterations. Let us fix $\delta\in(0,1)$ to be chosen later on. Set then $T\in(0,T_*]$, where $T_*\in(0,\infty]$ is the maximal time of existence of the strong solution, and then $\widetilde{\tau}>0$ such that $T-3\widetilde{\tau}\geq \frac{\tau}{2}$. We also fix $\widetilde{\tau}$ such that
\begin{equation}
2\widetilde{\tau}+\frac{\tau}{2}\leq \tau.
\label{tt}
\end{equation}
Let us then fix $\delta>0$ so that \eqref{delt1} and \eqref{delta} below hold. The choice of $\delta$  \textit{does not} depend on the specific $T$, but clearly depends on $\tau$ through \eqref{tt}. We now define, as usual in this argument, the sequence
\begin{align}
k_n=1-\delta-\frac{\delta}{2^n}, \quad \forall n\geq 0,
\label{kn}
\end{align}
where
\begin{align}
1-2\delta< k_n<k_{n+1}<1-\delta,\qquad \forall n\geq 1,\qquad k_n\to 1-\delta\qquad \text{as }n\to \infty,
\label{kn1}
\end{align}
and the sequence of times
\begin{align}
\begin{cases}
t_{-1}=T-3\widetilde{\tau},\\
t_n=t_{n-1}+\frac{\widetilde{\tau}}{2^n},\qquad n\geq 0,
\end{cases}
\end{align}
satisfying
$$
t_{-1}<t_n<t_{n+1}< T-\widetilde{\tau},\qquad \forall n\geq 0.
$$
We also introduce a cutoff function $\eta_n\in C^1(\R)$ by setting
\begin{align}
\eta_n(t):=\begin{cases}
0,\quad t\leq t_{n-1},\\
1,\quad t\geq t_{n},
\end{cases}\text{ and }\quad \vert \eta^\prime_n(t)\vert\leq \frac{2^{n+1}}{\widetilde{\tau}},
\label{cutoff}
\end{align}
on account of the above definition of the sequence $\{t_n\}_n$. We then set
\begin{align}
\varphi_n(x,t):=(\varphi-k_n)^+,
\label{phik0}
\end{align}
and, for any $n\geq 0$, we introduce the interval $I_n=[t_{n-1},T]$ and the set
$$
A_n(t):=\{x\in \Omega: \varphi(x,t)-k_n\geq 0\},\quad \forall t\in I_n.
$$
Clearly, we have
$$
I_{n+1}\subseteq I_n,\qquad \forall n\geq 0,$$
$$A_{n+1}(t)\subseteq A_n(t),\qquad \forall n\geq 0,\qquad \forall t\in I_{n+1}.
$$
In conclusion, we set
$$
y_n=\int_{I_n}\int_{A_n(s)}1dxds,\qquad \forall n\geq0.
$$
 Now, for any $n\geq 0$, we consider the test function $v=\varphi_n\eta_n^2$, and integrate over $[t_{n-1},t]$, $t_n\leq t\leq T$. Then we have
\begin{align}
&\nonumber\int_{t_{n-1}}^t\langle \partial_t\varphi,\varphi_n\eta_n^2\rangle_{\Hu',\Hu}ds+\int_{t_{n-1}}^t\int_{A_n(s)}\nabla\varphi\cdot \nabla\varphi_n \eta_n^2 dx ds+\int_{t_{n-1}}^t\int_{\Omega}F^\prime(\varphi)(\varphi_n-\overline{\varphi}_n)\eta_n^2dxds\\&=\int_{t_{n-1}}^t\int_{\Omega}\eta_n^2\theta_0\varphi(\varphi_n-\overline\varphi_n)dx ds,
\label{phin}
\end{align}
where we used
\begin{align*}
    \int_{t_{n-1}}^t\int_{\Omega}\eta_n^2(\mu-\overline{\mu})\varphi_ndxds=     \int_{t_{n-1}}^t\int_{\Omega}\eta_n^2\mu(\varphi_n-\overline \varphi_n)dxds.
\end{align*}
Furthermore, recall that, as also observed in \cite[Remark 4.7]{P}, since $\bv$ is divergence-free and $\bv=\mathbf 0$ on $\partial\Omega$,
 for almost any $t\geq 0$ we have
			\begin{align}
			\int_\Omega (\bv\cdot \nabla\varphi)\varphi_n\eta_n^2\: \d x=\int_\Omega \bv\cdot \frac 1 2\nabla\vert\varphi_n\vert^2\eta_n^2 \:\d x=0.\label{vanish}
			\end{align}
Now we obtain
\begin{align}
\int_{t_{n-1}}^t\eta_n^2\int_{A_n(s)}\nabla\varphi\cdot \nabla\varphi_n dx ds=\int_{t_{n-1}}^t\eta_n^2\Vert \nabla\varphi_n\Vert^2_{\bLd}ds,\label{nabla1}
\end{align}
but also, since, on $A_n(t)$, $F^\prime(\varphi(t))\geq F^\prime(1-2\delta)$,
$$
\int_{t_{n-1}}^t\int_{\Omega}F^\prime(\varphi)\varphi_n\eta_n^2dxds\geq F^\prime(1-\delta)\int_{t_{n-1}}^t\int_{A_n(s)}\varphi_n\eta_n^2dxds
$$
and, for the remaining terms in \eqref{phin}, recalling that $\vert\varphi\vert< 1$ a.e. in $\Omega\times (0,\infty)$ and that $\overline{F'(\varphi)}\in L^\infty(\tau,\infty)$, we find
\begin{align*}
&\int_{t_{n-1}}^t\int_{\Omega}\eta_n^2\theta_0\varphi(\varphi_n-\overline\varphi_n)dx ds+\int_{t_{n-1}}^t\int_{\Omega}\eta_n^2F'(\varphi)\overline\varphi_ndxds\\& =\int_{t_{n-1}}^t\int_{A_n(s)}\eta_n^2\theta_0\varphi\varphi_ndx ds-\int_{t_{n-1}}^t\int_{\Omega}\eta_n^2\theta_0\varphi\overline\varphi_ndx ds+\int_{t_{n-1}}^t\int_{\Omega}\eta_n^2\overline{F'(\varphi)}\varphi_ndxds
\\&
\leq C(\tau)\int_{t_{n-1}}^t\int_{A_n(s)}\eta_n^2\varphi_ndx ds,
\end{align*}
where we used the fact that $\varphi_n\geq 0$ and
\begin{align}
\overline\varphi_n(s)=\frac{\int_{A_n(s)}\varphi_n\dx}{\normmm{\Omega}}.
    \label{mean}
\end{align}
Moreover, we have
\begin{align}
\int_{t_{n-1}}^t<\partial_t\varphi,\varphi_n\eta_n^2>ds=\frac{1}{2}\Vert\varphi_n(t)\Vert^2_{\Ld}-\int_{t_{n-1}}^t\Vert\varphi_n(s)\Vert^2_{\Ld}\eta_n\partial_t\eta_nds.\label{timeder}
\end{align}
Note that, as crucially first observed in \cite{P}, since $ \vert\varphi\vert<1\text{ a.e. in }\Omega$, for any $t\geq \frac{\tau}{2}$,
\begin{align}
\label{delti}
0\leq \varphi_n\leq 2\delta\quad\text{ a.e. in }\Omega,\quad \forall t\geq \frac{\tau}{2}.
\end{align}
Then we have
\begin{align}
\nonumber
&\int_{t_{n-1}}^t\Vert\varphi_n(s)\Vert_{\Ld}^2\eta_n\partial_t\eta_nds=\int_{t_{n-1}}^t\int_\Omega \varphi_n^2(s)\eta_n\partial_t\eta_ndxds=\int_{t_{n-1}}^t\int_{A_n(s)} \varphi_n^2(s)\eta_n\partial_t\eta_ndxds\\&
\leq
\int_{t_{n-1}}^t\int_{A_n(s)} (2\delta)^2\frac{2^{n+1}}{\widetilde{\tau}}dxds\leq \frac{2^{n+3}\delta^2}{\widetilde{\tau}}y_n.
\label{del}
\end{align}
This allows us to deduce the estimate
\begin{align*}
&\frac{1}{2}\Vert\varphi_n(t)\Vert^2_{\Ld}+\left(F^\prime(1-2\delta)-C(\tau)\right)\int_{t_{n-1}}^t \eta_n^2\Vert\varphi_n(s)\Vert _{L^1(\Omega)}ds+\int_{t_{n-1}}^t \eta_n^2\Vert \nabla\varphi_n(s)\Vert^2_{\bLd}ds\\&\leq\frac{2^{n+3}\delta^2}{\widetilde{\tau}}y_n,
\end{align*}
for any $t\in[t_n,T]$. Notice that, choosing $\delta>0$ sufficiently small, since $F^\prime(1-2\delta)\to \infty$, we can find $C_0>0$ such that
\begin{align}
F^\prime(1-2\delta)-C(\tau)\geq {C_0},
\label{delt1}
\end{align}
so that we get
\begin{align}
\max_{t\in I_{n+1}}\Vert\varphi_n(t)\Vert^2_{\Ld}\leq X_n,\qquad  {\color{black}}2\int_{I_{n+1}}\Vert \nabla\varphi_n\Vert^2_{\bLd}ds \leq X_n,
\label{est}
\end{align}
where
$$
X_n:= \frac{2^{n+4}\delta^2}{\widetilde{\tau}}y_n.
$$
On the other hand, for any $t\in I_{n+1}$ and for almost any $x\in A_{n+1}(t)$, we get
\begin{align}
&\nonumber\varphi_n(x,t)=\varphi(x,t)-\left[1-\delta-\frac{\delta}{2^n}\right]\\&=
\underbrace{\varphi(x,t)-\left[1-\delta-\frac{\delta}{2^{n+1}}\right]}_{\varphi_{n+1}(x,t)\geq 0}+\delta\left[\frac{1}{2^{n}}-\frac{1}{2^{n+1}}\right]\geq \frac{\delta}{2^{n+1}},\label{basic}
\end{align}
which implies
\begin{align*}
\int_{I_{n+1}}\int_{\Omega}\vert\varphi_n\vert^3dxds\geq \int_{I_{n+1}}\int_{A_{n+1}(s)}\vert\varphi_n\vert^3dxds\geq \left(\frac{\delta}{2^{n+1}}\right)^3\int_{I_{n+1}}\int_{A_{n+1}(s)}dxds=\left(\frac{\delta}{2^{n+1}}\right)^3y_{n+1}.
\end{align*}
Then we have
\begin{align}
&\nonumber\left(\frac{\delta}{2^{n+1}}\right)^3y_{n+1}\leq \int_{I_{n+1}}\int_{\Omega}\vert\varphi_n\vert^3dxds\\&=\int_{I_{n+1}}\int_{A_n(s)}\vert\varphi_n\vert^3dxds \leq \left(\int_{I_{n+1}}\int_{\Omega}\vert\varphi_n\vert^{\frac{10}{3}}dxds\right)^{\frac{9}{10}}\left(\int_{I_{n+1}}\int_{A_n(s)}1\ dxds\right)^{\frac{1}{10}}.
\label{est2}
\end{align}
Notice that, by Gagliardo-Nirenberg's inequalities, we obtain
\begin{align*}
&\int_{I_{n+1}}\int_{\Omega}\vert\varphi_n\vert^{\frac{10}{3}}dxds\leq \hat{C} \int_{I_{n+1}}\Vert\varphi_n\Vert_{H^1(\Omega)}^{2}\Vert\varphi_n\Vert^{\frac{4}{3}}_{\Ld}ds\\
&\leq \hat{C}\int_{I_{n+1}} \left(\Vert\varphi_n\Vert^2_{\Ld}+\Vert\nabla\varphi_n\Vert^2_{\bLd}\right)\Vert\varphi_n\Vert^{\frac{4}{3}}_{\Ld}ds.
\end{align*}
By \eqref{est} we then have
\begin{align*}
&\int_{I_{n+1}}\int_{\Omega}\vert\varphi_n\vert^{\frac{10}{3}}dxds\leq \hat{C}\max_{t\in I_{n+1}}\Vert\varphi_n(t)\Vert_{\Ld}^{\frac{4}{3}}\left(3\widetilde \tau\max_{t\in I_{n+1}}\Vert\varphi_n(t)\Vert^2_{\Ld}+\int_{I_{n+1}}\Vert\nabla\varphi_n\Vert^2_{\bLd}ds\right)\\&
\leq \frac{{\color{black}}\hat{C}}{2}X_n^{\frac{2}{3}}\left(6\widetilde\tau X_n+2\int_{I_{n+1}}\Vert\nabla\varphi_n\Vert^2_{\bLd}ds\right)\leq \frac{{\color{black}}\hat{C}}{2}X_n^{\frac{5}{3}}(6\widetilde\tau +1)\leq
\frac{{\color{black}}\hat{C}(6\widetilde\tau +1)}{2}\dfrac{2^{\frac 5 3n +\frac{20}{3}}\delta^{\frac {10} 3}}{\widetilde{\tau}^{\frac  53}}y_n^{\frac{5}{3}}.
\end{align*}
Coming back to \eqref{est2}, we immediately infer
\begin{align}
&\nonumber\left(\frac{\delta}{2^{n+1}}\right)^3y_{n+1}\leq \left(\int_{I_{n+1}}\int_{\Omega}\vert\varphi_n\vert^{\frac{10}{3}}dxds\right)^{\frac{9}{10}}\left(\int_{I_{n+1}}\int_{A_n(s)}1\ dxds\right)^{\frac{1}{10}}\\&\leq \delta^3\frac{2^{{\frac{3}{2}n+{6}}}\hat{C}^{\frac{9}{10}}(6\widetilde\tau +1)^{\frac{9}{10}}}{2^{\frac{9}{10}}\widetilde{\tau}^{\frac 3{2}}}y_n^{\frac{8}{5}}.
\label{est3}
\end{align}
In conclusion, we end up with
\begin{align}
y_{n+1}\leq \frac{2^{{\frac{9}{2}n+\frac{81}{10}}}\hat{C}^{\frac{9}{10}}(6\widetilde\tau +1)^{\frac{9}{10}}}{\widetilde{\tau}^{\frac 3{2}}}y_n^{\frac{8}{5}},\qquad \forall n\geq 0.
\label{last01}
\end{align}
Thus we can apply the well known geometric lemma \cite[Lemma 3.8]{P}. In particular, in the notation of the lemma, we have $b=2^\frac{9}{2}>1$, $C=\frac{2^{{\frac{81}{10}}}\hat{C}^{\frac{9}{10}}(6\widetilde\tau +1)^{\frac{3}{2}}}{\widetilde{\tau}^{\frac 9{10}}}>0$, $\varepsilon=\frac{3}{5}$, to get that ${y}_n\to 0$, as long as
$$
{y}_0\leq C^{-\frac{5}{3}}b^{-\frac{25}{9}},
$$
i.e.,
\begin{align}
y_0\leq \dfrac{2^{-26}\widetilde{\tau}^\frac 32}{\hat{C}^{\frac{3}{2}}(6\widetilde\tau +1)^{\frac{5}{2}}}.
\label{last}
\end{align}
Then, recalling that $\Vert F^\prime(\varphi)\Vert_{L^\infty(\frac{\tau}{2},T;L^1(\Omega))}\leq C(\tau)$  and $F^\prime$ is monotone in a neighborhood of $+1$, we immediately infer
\begin{align*}
&y_0=\int_{I_0}\int_{A_0(s)}1dxds\leq\int_{I_0}\int_{\{x\in\Omega:\ \varphi(x,t) \geq 1-2\delta\}}1dxds\\&\leq \int_{I_0}\int_{A_0(s)}\frac{\vert F^\prime(\varphi)\vert}{F^\prime(1-2\delta)} dxds\leq \frac{3C(\tau)\widetilde{\tau}}{{F^\prime(1-2\delta)}}.
\end{align*}
Therefore, if we ensure that
$$
\frac{3C(\tau)\widetilde{\tau}}{{F^\prime(1-2\delta)}}\leq \dfrac{2^{-26}\widetilde{\tau}^\frac 3 2}{\hat{C}^{\frac{3}{2}}(6\widetilde\tau +1)^{\frac{5}{2}}},
$$
then \eqref{last} holds. This means that we need
\begin{align}
\label{delta}
\frac{3C(\tau)2^{26}\hat{C}^{\frac{3}{2}}(6\widetilde\tau +1)^{\frac{5}{2}}}{\widetilde{\tau}^{\frac 1 2}}\leq {F^\prime(1-2\delta)}.
\end{align}
Having fixed $\widetilde{\tau}$ such that \eqref{tt} holds,
we obtain the result by choosing $\delta$ sufficiently small, since $F^\prime(1-2\delta)\to \infty$ as $\delta\to0$.
In the end, passing to the limit in $y_n$ as $n\to\infty$, we have obtained that
$$
\Vert(\varphi-(1-\delta))^+\Vert_{L^\infty(\Omega\times({T}-\widetilde{\tau},{T}))}=0.
$$
 We now repeat ($\widetilde{\tau}$ can be kept the same) exactly the same argument for the case $(\varphi-(-1+\delta))^-$ (using $\varphi_n(t)=(\varphi(t)+k_n)^-$). Moreover, the argument is exactly the same since $F^\prime(-1+2\delta)\to-\infty$ as $\delta\to0$. We can then choose the minima between the $\delta$s  obtained in the two cases, to get in the end that there exists $\delta=\delta(\tau)>0$ such that
 \begin{align}
-1+\delta\leq \varphi(x,t)\leq 1-\delta,\quad\text{a.e. in }\Omega\times ({T}-\widetilde{\tau},{T}).
\label{end}
\end{align}
Finally, notice that, due to the choice of $\widetilde \tau$, we have $T-\widetilde{\tau}=2\widetilde{\tau}+\frac{\tau}{2}\leq \tau$, therefore, if $T_*<\infty$ we are done by choosing $T=T_*$, otherwise, we can repeat the same procedure on the interval $({T},T+\widetilde{\tau})$ (this means that the new starting time will be $t_{-1}={T}-2\widetilde{\tau}\geq \frac{\tau}{2}$) and so on, reaching eventually the entire interval $[\tau,T_*)$. The proof of \eqref{sepa1} is thus concluded.

In case the initial datum $\varphi_0$ is already strictly separated, we can repeat a similar proof without using De Giorgi iterations (say without using any cutoff function), but simply using a modification of the maximum principle (see, e.g., \cite[Theorem 3.7]{GGPC}, \cite[Theorem 3.1]{GPCAC}) to obtain that, given $0<T<T_*$, it holds

\begin{align}
-1+\delta\leq \varphi(x,t)\leq 1-\delta,\quad\text{a.e. in }\Omega\times [0,{T}).
\label{end1}
\end{align}
As a consequence, putting together this result with  \eqref{end}, we immediately infer \eqref{sepa2}, concluding the proof of the theorem.\label{sec:sep1}

    \section{Proof of Theorem \ref{THM:WSU}: Conditional weak-strong uniqueness}
    \label{weak-strong}
    
    We consider two sets of initial data $(\bv_{0,1},\varphi_{0,1})$ and $(\bv_{0,2},\varphi_{0,2})$ satisfying the assumptions in Theorems \ref{THM:WeakExistence} and \ref{THM:StrongWP}, respectively. We assume that $\varphi_{0,2}$ is strictly separated from pure phases, namely, there exists $\delta_0\in(0,1): \norm{\varphi_{0,2}}_{L^\infty(\Omega)}\leq 1-\delta_0$. Also, we suppose that $\overline{\vphi_{0,1}}=\overline{\vphi_{0,2}}$. Let $(\bv_1,\varphi_1,\mu_1)$ be a weak solution corresponding to $(\bv_{0,1},\varphi_{0,1})$ and let $(\bv_2,\varphi_2,\mu_2)$ be the (unique) strong solution corresponding to $(\bv_{0,2},\varphi_{0,2})$. Our conditional assumption is
    \begin{align}
\norm{F'(\varphi_1)}_{L^4(0,T_*;L^4(\Omega))}\leq C(T_*).
        \label{L4}
    \end{align}
    Let us set
    \begin{align*}
        (\bv_0,\varphi_0) &:= (\bv_{0,1},\varphi_{0,1}) -(\bv_{0,2},\varphi_{0,2}), \\
        (\bv,\varphi,\mu) &:= (\bv_1,\varphi_1,\mu_1)-(\bv_2,\varphi_2,\mu_2).
    \end{align*}
    Note that, by mass conservation, $\overline{\vphi}(t)=0$ for any $t\geq0$.
    This means that $(\bv_i,\varphi_i,\mu_i)$, $i=1,2$, fulfill the following systems of equations at least in a weak sense:
    \begin{align}
        \partial_t\bv_1 + (\bv_1\cdot\nabla)\bv_1 - \text{div}(\nu(\varphi_1)D\bv_1) 
         + \nabla \pi_1 &= \mu_1\nabla\varphi_1 &\quad&\text{ in }\Omega_T, \label{Eq:1'0}\\
        \nabla\cdot\bv_1 &= 0 &\quad&\text{ in }\Omega_T, \label{Eq:2'0}\\
        \partial_t\varphi_1 + \bv_1\cdot\nabla\varphi_1  + \mu_1 - \overline{\mu}_1&= 0 &\quad&\text{ in }\Omega_T, \label{Eq:3'0}\\
         \mu_1 &= -\Delta \vphi_1+f'(\vphi_1)&\quad&\text{ in }\Omega_T, \label{Eq:3'01}\\
        \bv_1 = \partial_\mathbf{n}\varphi_1 &= 0 &\quad&\text{ on }\partial\Omega_T, \label{Eq:4'0}\\
        \bv_1(0) = \bv_{0,1},\quad \varphi_1(0) &= \varphi_{0,1} &\quad&\text{ in }\Omega,
        \label{Eq:5'0}
    \end{align}
    \begin{align}
        \partial_t\bv_2 + (\bv_2\cdot\nabla)\bv_2 -  \text{div}(\nu(\varphi_2)D\bv_2)  + \nabla \pi_2 &= \mu_2\nabla\varphi_2 &\quad&\text{ in }\Omega_T, \label{Eq:1'}\\
        \nabla\cdot\bv_2 &= 0 &\quad&\text{ in }\Omega_T, \label{Eq:2'}\\
        \partial_t\varphi_2 +  \bv_2\cdot\nabla\varphi_2 + \mu_2 - \overline{\mu}_2 &= 0 &\quad&\text{ in }\Omega_T, \label{Eq:3'}\\
         \mu_2 &= -\Delta \vphi_2+f'(\vphi_2)&\quad&\text{ in }\Omega_T, \label{Eq:3'01b}\\
        \bv_2= \partial_\mathbf{n}\varphi_2 &= 0 &\quad&\text{ on }\partial\Omega_T, \label{Eq:4'}\\
        \bv_2(0) = \bv_{0,2},\quad \varphi_2(0) &= \varphi_{0,2} &\quad&\text{ in }\Omega. \label{Eq:5'}
    \end{align}
    The ``natural'' approach is take the difference between the two systems and then test the resulting equations by $\bv$ and 
    $\mu - \overline{\mu}$. Nevertheless this is only formal, as $\bv_1\in C_w([0,T];\mathbf{L}^2_\sigma(\Omega))$, and thus the chain rule for $\bv=\bv_1-\bv_2$ in $\bLds$ does not apply. For a rigorous proof, we resort to a relative energy approach. Namely, we recall that the two solutions satisfy by assumption the energy inequality
 \begin{align}
                &\nonumber E(\bv_i(t),\varphi_i(t)) + \int_0^t\left[\|\sqrt{\nu(\varphi_i(\tau))}D\bv_i(\tau)\|_{\mathbf{L}^2(\Omega)}^2 + \|\mu_i(\tau)-\overline{\mu_i(\tau)}\|_{L^2(\Omega)}^2\right]\d\tau \\&\leq E(\bv_0,\varphi_0),\quad i=1,2.\label{ineq1A}
            \end{align}
Moreover, notice that it holds
\begin{align*}
   \onehalf \norm{\bv_1-\bv_2}^2_{\bLd}=\onehalf \norm{\bv_1}^2_{\bLd}+\onehalf \norm{\bv_2}^2_{\bLd}-\int_\Omega \bv_1\cdot \bv_2\dx.
\end{align*}
On the other hand, we have
\begin{align*}
    &-\int_\Omega \bv_1(t)\cdot \bv_2(t)\dx\\
    &=-\int_\Omega\bv_{0,1}\cdot \bv_{0,2}\dx-\int_0^t\int_\Omega\pt \bv_1(\tau)\cdot \bv_2(\tau)\dx\d\tau-\int_0^t\int_\Omega\bv_1(\tau)\cdot\pt\bv_2(\tau)\dx\d\tau.
\end{align*}
Similarly, we get
\begin{align*}
   \onehalf \norm{\nabla\vphi_1-\nabla\vphi_2}^2_{\bLd}=\onehalf \norm{\nabla\vphi_1}^2_{\bLd}+\onehalf \norm{\nabla\vphi_2}^2_{\bLd}-\int_\Omega\nabla\vphi_1\cdot \nabla\vphi_2\dx,
\end{align*}
and  
\begin{align*}
   & -\int_\Omega\nabla\vphi_1(t)\cdot \nabla\vphi_2(t)\dx\\&=-\int_\Omega{\nabla\vphi_{0,1}}\cdot\nabla \vphi_{0,2}\dx+\int_0^t\int_\Omega\pt \vphi_1(\tau)\Delta\vphi_2(\tau)\dx\d\tau+\int_0^t\int_\Omega\Delta\vphi_1(\tau)\pt\vphi_2(\tau)\dx\d\tau.
\end{align*}
The above identity can be obtained first for smooth solutions and then using a density argument.

Additionally, we observe that 
\begin{align*}
    \norm{\mu-\overline \mu}_{\Ld}^2= \norm{\mu_1-\overline \mu_1}_{\Ld}^2+ \norm{\mu_2-\overline \mu_2}_{\Ld}^2-2\int_\Omega (\mu_1-\overline \mu_1) (\mu_2-\overline \mu_2)\dx,
\end{align*}
as well as 
\begin{align*}
    \int_\Omega \nu(\vphi_2) \normmm{D\bv}^2\dx=\int_\Omega \nu(\vphi_2) \normmm{D\bv_1}^2\dx+\int_\Omega \nu(\vphi_2) \normmm{D\bv_2}^2\dx-2\int_\Omega \nu(\vphi_2) D\bv_1:D\bv_2\dx.
\end{align*}
Recalling the identities above, we now test, in the weak formulation integrated over $(0,t)$, \eqref{Eq:1'0} with $-\bv_2$, \eqref{Eq:3'0} with $-(\mu_2-\overline{\mu_2})$, \eqref{Eq:1'} with $-\bv_1$, and \eqref{Eq:3'} with $-(\mu_1-\overline{\mu_1})$. Then, we sum everything up with the energy inequalities \eqref{ineq1A}, for $i=1,2$. This gives
\begin{align}
    &\frac{1}{2}\Big(\|\bv(t)\|_{\mathbf{L}^2_\sigma(\Omega)}^2 + \|\nabla\varphi(t)\|_{\mathbf{L}^2(\Omega)}^2\Big) + \int_0^t\|\mu-\overline{\mu}\|_{L^2(\Omega)}^2\d s+\int_0^t\int_\Omega \nu(\vphi_2) \normmm{D\bv}^2\dx\d s\nonumber\\&
    \leq\nonumber \onehalf\norm{\bv_{0,1}-\bv_{0,2}}_{\bLds}^2+\onehalf\norm{\nabla\vphi_{0,1}-\nabla\vphi_{0,2}}_{\bLd}^2 \\&\quad +\underbrace{\int_0^t\int_\Omega (\nu(\vphi_2)-\nu(\vphi_1))\normmm{D\bv_1}^2\dx\d s+\int_0^t\int_\Omega (\nu(\vphi_1)-\nu(\vphi_2))D\bv_1:D\bv_2\dx\d s}_{A_1}\nonumber\\&
    \quad +\underbrace{\int_0^t\int_\Omega (\bv_1\cdot\nabla \bv_1)\bv_2\dx\d s+\int_0^t\int_\Omega (\bv_2\cdot\nabla \bv_2)\bv_1\dx\d s}_{A_2}\nonumber\\&
    \quad +\underbrace{\int_0^t\int_\Omega \big((\bv_1\cdot \nabla \vphi_1)\mu_2)+(\bv_2\cdot \nabla \vphi_2)\mu_1-(\mu_1\nabla \vphi_1)\cdot \bv_2-(\mu_2\nabla \vphi_2)\cdot\bv_1\big)\dx\d s}_{A_3}\nonumber\\&
    \quad +\underbrace{\int_0^t\int_\Omega (f'(\vphi_1)-f'(\vphi_2))\pt \vphi\dx\dt}_{A_4}.
    \label{Est:WSU:1}
\end{align}
Note that the difference of the initial data appears since, for instance,  after testing \eqref{Eq:1'0} with $-\bv_2$ and integrating in time over $(0,t)$, one has to integrate by parts to complete the squares, namely,
\begin{align*}
    -\int_0^t\int_\Omega \pt\bv_1\cdot\bv_2\dx= \int_0^t\int_\Omega \pt\bv_2\cdot\bv_1-\int_\Omega \bv_1(t)\cdot\bv_2(t)\dx +\int_\Omega \bv_{0,1}\cdot\bv_{0,2}\dx,
\end{align*}
and this completes the square on the right-hand side with $\frac12\norm{\bv_{0,1}}_{\bLds}^2+\frac12\norm{\bv_{0,2}}_{\bLds}^2$ coming from summing the energy inequality \eqref{ineq1A} for both $\bv_1$ and $\bv_2$. The same argument is used for $\nabla\vphi_{0,i}$, $i=1,2$. On the other hand, $A_4$ takes that form since in the left-hand side of \eqref{Est:WSU:1} we would have   
\begin{align*}
    \int_\Omega f(\vphi_1(t))\dx+ \int_\Omega f(\vphi_2(t))\dx-\int_0^t\int_\Omega \pt\vphi_1 f'(\vphi_2)\dx\d s-\int_0^t\int_\Omega \pt\vphi_2 f'(\vphi_1)\dx\d s,
\end{align*}
so that, summing and subtracting $\int_0^t\int_\Omega f'(\vphi_i)\pt\vphi_i\dx\d s$, $i=1,2$, we get
\begin{align*}
   & \int_\Omega f(\vphi_1(t))\dx+ \int_\Omega f(\vphi_2(t))\dx-\int_0^t\int_\Omega \pt\vphi_1 f'(\vphi_2)\dx\d s-\int_0^t\int_\Omega \pt\vphi_2 f'(\vphi_1)\dx\d s\\&= A_4+\int_\Omega f'(\vphi_{0,1})\dx+\int_\Omega f'(\vphi_{0,2})\dx,
\end{align*}
and the initial data simplify with the one arising from the sum of the energy inequalities \eqref{ineq1A}.

Observe now that 
\begin{align*}
    A_1=\int_0^t\int_\Omega (\nu(\vphi_2)-\nu(\vphi_1))D\bv_1:D\bv\dx\dt.
\end{align*}
Moreover, recalling that $\bv_i$ are divergence free and $\bv_i=\mathbf 0$ almost everywhere on $\partial\Omega\times (0,T_*)$ , we get
\begin{align*}
    A_2&=-\int_0^t\int_\Omega (\bv_1\cdot\nabla \bv_1)\bv\dx\d s+\int_0^t\int_\Omega (\bv_2\cdot\nabla \bv_2)\bv\dx\d s\\&
    =-\int_0^t\int_\Omega (\bv_1\cdot\nabla \bv)\bv\dx\d s-\int_0^t\int_\Omega (\bv\cdot\nabla \bv_2)\bv\dx\d s\\&
    =-\int_0^t\int_\Omega (\bv\cdot\nabla \bv_2)\bv\dx\d s.
\end{align*}
Concerning $A_3$, we have
\begin{align*}
A_3=\int_0^t\int_\Omega\mu_2\nabla\varphi\cdot\bv\dx -\int_0^t\int_\Omega \bv_2\cdot\nabla\varphi\mu\dx.
\end{align*}
% In conclusion, we finally rewrite $A_4$ as 
% \begin{align*}
%     A_4=-\int_0^t\int_\Omega(f'(\vphi_1)-f'(\vphi_2))\pt \vphi \dx\d s.
% \end{align*}
% \begin{align}\label{Est:WSU:1}
%     \begin{split}
%         &\qquad \frac{1}{2}\ddt\Big(\|\bv\|_{\mathbf{L}^2_\sigma(\Omega)}^2 + \|\nabla\varphi\|_{\mathbf{L}^2(\Omega)}^2\Big) + \|\mu-\overline{\mu}\|_{L^2(\Omega)}^2 \\
%         &= \int_\Omega \mu_1\nabla \varphi_1 \cdot\bv\dx - \int_\Omega\mu_2\nabla\varphi_2\cdot\bv\dx - \int_\Omega\big((\bv_1\cdot\nabla)\bv_1 - (\bv_2\cdot\nabla)\bv_2\big)\cdot\bv\dx \\
%         &- \int_\Omega\big(\bv_1\cdot\nabla\varphi_1 - \bv_2\cdot\nabla\varphi_2\big)\mu\dx - \int_\Omega\big(f^\prime(\varphi_1)-f^\prime(\varphi_2)\big)\partial_t\varphi\dx \\
%         &+ \int_\Omega\nu(\varphi_1)D\bv_1:D\bv\dx - \int_\Omega\nu(\varphi_2)D\bv_2:D\bv\dx.
%     \end{split}
%     \end{align}
%     For the term involving the viscosity, we first observe that
%     \begin{align*}
%         &\int_\Omega\nu(\varphi_1)D\bv_1:D\bv\dx - \int_\Omega\nu(\varphi_2)D\bv_2:D\bv\dx \\
%         &= \int_\Omega(\nu(\varphi_1) - \nu(\varphi_2))D\bv_2:D\bv\dx + \int_\Omega\nu(\varphi_2)|D\bv|^2\dx.
%     \end{align*}
Now, recalling \ref{ASS:Viscosity}, it follows that
    \begin{align}\label{Est:WSU:nu1}
        \int_0^t\int_\Omega\nu(\varphi_2)|D\bv|^2\dx\d s \geq \nu_*\int_0^t\|D\bv\|_{\mathbf{L}^2(\Omega)}^2\d s.
    \end{align}
    Consider now $A_1$. Exploiting the Lipschitz continuity of $\nu$, since $\bv_2\in L^\infty(0,T_*;\mathbf H^1_\sigma(\Omega))$, we infer that
\begin{align}\label{Est:WSU:nu2}
    \begin{split}
        \left|\int_\Omega(\nu(\varphi_1) - \nu(\varphi_2))D\bv_2:D\bv\dx\right|
        &\leq C\|\varphi\|_{L^6(\Omega)}\|D\bv_2\|_{\mathbf{L}^3(\Omega)}\|D\bv\|_{\mathbf{L}^2(\Omega)} \\
        &\leq C\|\nabla\varphi\|_{\mathbf{L}^2(\Omega)}\|\bv_2\|_{\mathbf{H}^2(\Omega)}^{1/2}\|D\bv\|_{\mathbf{L}^2(\Omega)} \\
        &\leq C\|\bv_2\|_{\mathbf{H}^2(\Omega)}\|\nabla\varphi\|_{\mathbf{L}^2(\Omega)}^2 + \frac{\nu_*}{8}\|D\bv\|_{\mathbf{L}^2(\Omega)}^2.
    \end{split}
    \end{align}
    This gives
    \begin{align*}
        \normmm{A_1}\leq C\int_0^t\|\bv_2\|_{\mathbf{H}^2(\Omega)}\|\nabla\varphi\|_{\mathbf{L}^2(\Omega)}^2\d s + \frac{\nu_*}{8}\int_0^t\|D\bv\|_{\mathbf{L}^2(\Omega)}^2\d s.
    \end{align*}
   Then, to estimate $A_2$, we use again that $\bv_2\in L^\infty(0,T_*;\mathbf{H}^1_\sigma(\Omega))$ as well as the Gagliardo--Nirenberg inequality and Young's inequality to deduce
    \begin{align*}
        \left|\int_\Omega(\bv\cdot\nabla)\bv_2\cdot\bv\dx\right| &\leq C\norm{\bv}_{\mathbf{L}^3(\Omega)}\norm{\bv}_{\mathbf{L}^6_\sigma(\Omega)}\norm{\nabla\bv_2}_{\mathbf{L}^2_\sigma(\Omega)}\\&\leq C\|\bv\|_{\mathbf{L}^2_\sigma(\Omega)}^{\frac12}\|\bv\|_{\mathbf{H}^1_\sigma(\Omega)}^{\frac32} \\&\leq C\|\bv\|_{\mathbf{L}^2_\sigma(\Omega)}^2 + \frac{\nu_*}{8}\|D\bv\|_{\mathbf{L}^2(\Omega)}^2,
    \end{align*}
    so that 
    \begin{align*}
        \normmm{A_2} \leq C\int_0^t\|\bv\|_{\mathbf{L}^2_\sigma(\Omega)}^2\d s + \frac{\nu_*}{8}\int_0^t\|D\bv\|_{\mathbf{L}^2(\Omega)}^2\d s.
    \end{align*}
   The first term in $A_3$ can be controlled as follows
\begin{align}\label{Est:WSU:2}
    \begin{split}
        \left|\int_\Omega\mu_2\nabla\varphi\cdot\bv\dx\right|
        &\leq \|\mu_2\|_{L^6(\Omega)}\|\nabla\varphi\|_{\mathbf{L}^2(\Omega)}\|\bv\|_{\mathbf{L}^3(\Omega)} \\
        &\leq C\|\mu_2\|_{H^1(\Omega)}\|\nabla\varphi\|_{\mathbf{L}^2(\Omega)}\|\bv\|_{\mathbf{L}^2(\Omega)}^{1/2}\|D\bv\|_{\mathbf{L}^2(\Omega)}^{1/2} \\
        &\leq C\|\mu_2\|_{H^1(\Omega)}^2\|\nabla\varphi\|_{\mathbf{L}^2(\Omega)}^2 + C\|\bv\|_{\mathbf{L}^2(\Omega)}^2 + \frac{\nu_*}{8}\|D\bv\|_{\mathbf{L}^2(\Omega)}^2.
    \end{split}
    \end{align}
    As far as the second term is concerned, we use that $\overline{\mu}$ is constant in space. Therefore, we have
    \begin{align}
    \begin{split}
        \left|\int_\Omega \bv_2\cdot\nabla\varphi\mu\dx\right|& = \left|\int_\Omega \bv_2\cdot\nabla\varphi(\mu-\overline{\mu})\dx\right|
       \\& \leq \|\bv_2\|_{\mathbf{L}^\infty(\Omega)}\|\nabla\varphi\|_{\mathbf{L}^2(\Omega)}\|\mu-\overline{\mu}\|_{L^2(\Omega)} \\
        &\leq C\|\bv_2\|_{\mathbf{H}^2(\Omega)}\|\bv_2\|_{\mathbf{H}^1(\Omega)}\|\nabla\varphi\|_{\mathbf{L}^2(\Omega)}^2 + \frac{1}{4}\|\mu-\overline{\mu}\|_{L^2(\Omega)}^2\\&\leq C\|\bv_2\|_{\mathbf{H}^2(\Omega)}\|\nabla\varphi\|_{\mathbf{L}^2(\Omega)}^2 + \frac{1}{4}\|\mu-\overline{\mu}\|_{L^2(\Omega)}^2,
    \end{split}
    \end{align}
    where we also used Agmon's inequality and $\bv_2\in L^\infty(0,T_*;\bH^1(\Omega))$. As a consequence, we deduce
    \begin{align*}
        \normmm{A_3}&\leq C\int_0^t\|\mu_2\|_{H^1(\Omega)}^2\|\nabla\varphi\|_{\mathbf{L}^2(\Omega)}^2\d s + C\int_0^t\big(\|\bv\|_{\mathbf{L}^2(\Omega)}^2 + \frac{\nu_*}{8}\|D\bv\|_{\mathbf{L}^2(\Omega)}^2\big)\d s\\&+C\int_0^t\|\bv_2\|_{\mathbf{H}^2(\Omega)}\|\nabla\varphi\|_{\mathbf{L}^2(\Omega)}^2\d s + \frac{1}{4}\int_0^t\|\mu-\overline{\mu}\|_{L^2(\Omega)}^2\d s.
    \end{align*}
    Finally, we need to control $A_4$.
    To this end, we need to estimate the $L^2$-norm of the difference of the potentials $f'(\vphi_i)$, $i=1,2$. In order to do that, we argue as in \cite{HKP}. This is possible since, thanks to the strict separation of $\varphi_{0,2}$, by Theorem \ref{THM:StrongWP} there exists $\delta\in(0,\tfrac12)$ such that
    \begin{align}
    \sup_{t\in[0,T_*)}\norm{\varphi_2(t)}_{C(\overline\Omega)}\leq 1-2\delta.
        \label{Est:WSU:SP}
    \end{align}
    For any $t\in[0,T_*)$, we introduce the following sets:
    \begin{align*}
        A_\delta(t) &:= \{x\in\Omega:\; |\varphi_2(x,t)|\geq 1-\delta\}, \\
        B_\delta(t) &:= \{x\in\Omega:\; |\varphi_2(x,t)-\varphi_1(x,t)|\geq \delta\}.
    \end{align*}
    Owing to \eqref{Est:WSU:SP}, we observe that
    \begin{align*}
        1- \delta \leq |\varphi_2(x,t)| \leq |\varphi_2(x,t)| + |\varphi_1(x,t)-\varphi_2(x,t)| \leq 1-2\delta + |\varphi_1(x,t)-\varphi_2(x,t)|
    \end{align*}
    for all $t\in[0,T_*)$ and $x\in A_\delta(t)$. In particular, this implies
    \begin{align*}
        |\varphi_2(x,t)-\varphi_1(x,t)|\geq \delta,
        \quad t\in[0,T_*),\; x\in A_\delta(t).
    \end{align*}
    Hence, we have
        $A_\delta(t) \subset B_\delta(t)$  for all $t\in[0,T_*)$. Thus, by Chebyshev's inequality, we deduce
    \begin{align}\label{Est:Chebyshev}
        |A_\delta(t)| \leq \int_{B_\delta(t)} 1\dx \leq \int_{B_\delta(t)} \frac{|\varphi_1(t)-\varphi_2(t)|^p}{\delta^p}\dx \leq \int_{\Omega} \frac{|\varphi_1(t)-\varphi_2(t)|^p}{\delta^p}\dx,
    \end{align}
    for all $t\in[0,T_*)$ and any $p\geq1$. Therefore, invoking Hölder's inequality as well as the Fundamental Theorem of Calculus, we deduce that
    \begin{align}
       \|f^\prime(\varphi_1)-f^\prime(\varphi_2)\|_{L^2(\Omega)} 
       &\leq \|f^\prime(\varphi_1)-f^\prime(\varphi_2)\|_{L^2(A_\delta)} + \|f^\prime(\varphi_1)-f^\prime(\varphi_2)\|_{L^2(\Omega\setminus A_\delta)} \nonumber\\
       &\leq \|f^\prime(\varphi_1)-f^\prime(\varphi_2)\|_{L^4(A_\delta)}|A_\delta|^{1/4} \nonumber \\
       &+ \left(\int_{\Omega\setminus A_\delta}\left|\int_0^1 f^{\prime\prime}(s\varphi_1 + (1-s)\varphi_2)(\varphi_1-\varphi_2)\d s\right|^2\dx\right)^{1/2}.  \label{Est:WSU:3}
    \end{align}
    Recalling \eqref{Est:WSU:SP} and the definition of $A_\delta$, we infer
    \begin{align*}
        |s\varphi_1(t)+(1-s)\varphi_2(t)| \leq s|\varphi_1(t)| + (1-s)|\varphi_2(t)| \leq 1-\delta \quad \text{ a.e. in }\Omega\setminus A_\delta(t)
    \end{align*}
    for all $t\in[0,T_*)$ and $s\in[0,1]$. Since $F^{\prime\prime}\in C(-1,1)$, we thus have
    \begin{align}\label{Est:WSU:4}
    \begin{split}
        &\left\vert\int_0^1 f^{\prime\prime}(s\varphi_1(t) + (1-s)\varphi_2(t))(\varphi_1(t)-\varphi_2(t))\d s \right\vert\\
        &\leq \left(\max_{|s|\leq 1-\delta}F^{\prime\prime}(s) +\theta_0\right)|\varphi_1(t)-\varphi_2(t)| \leq C_\delta|\varphi_1(t)-\varphi_2(t)| \quad \text{ a.e. in }\Omega\setminus A_\delta(t),
    \end{split}
    \end{align}
    where $C_\delta:=\max\{\frac1\delta, \max_{|s|\leq 1-\delta}F^{\prime\prime}(s) +\theta_0\}$.
    Moreover, crucially exploiting Chebyshev's inequality \eqref{Est:Chebyshev} with $p=4$ as well as Poincaré's inequality, we infer, by the embedding $H^1(\Omega)\hookrightarrow L^4(\Omega)$, that
    \begin{align}\label{Est:WSU:5}
        |A_\delta|^{1/4} \leq \Big(\int_{\Omega} \frac{|\varphi_1(t)-\varphi_2(t)|^4}{\delta^4}\dx\Big)^{\frac14} = \frac1\delta\|\varphi\|_{L^4(\Omega)} \leq C_\delta\|\nabla\varphi\|_{\mathbf{L}^2(\Omega)},
    \end{align}
   recalling the definition of $C_\delta>0$.
    Hence, inserting \eqref{Est:WSU:4} and \eqref{Est:WSU:5} into \eqref{Est:WSU:3}, we arrive at
    \begin{align}\label{Est:WSU:6}
        \|f^\prime(\varphi_1)-f^\prime(\varphi_2)\|_{L^2(\Omega)} \leq C\Big(1 + \|f^\prime(\varphi_1)\|_{L^4(\Omega)} + \|f^\prime(\varphi_2)\|_{L^4(\Omega)}\Big)\|\nabla\varphi\|_{\mathbf{L}^2(\Omega)}.
    \end{align}
    Observe now that, since it holds, in weak formulation,
    \begin{align*}
        \pt \vphi+\bv_1\cdot\nabla \vphi_1-\bv_2\cdot\nabla\vphi_2+\mu-\overline\mu=0,
    \end{align*}
    then we can write
    \begin{align*}
        A_4
        &= -\int_\Omega\big(f^\prime(\varphi_1)-f^\prime(\varphi_2)\big)(\mu-\overline{\mu})\dx - \int_\Omega\big(f^\prime(\varphi_1)-f^\prime(\varphi_2)\big)\bv_1\cdot\nabla\varphi_1\dx \\
        & \quad +\int_\Omega\big(f^\prime(\varphi_1)-f^\prime(\varphi_2)\big)\bv_2\cdot\nabla\varphi_2\dx.
    \end{align*}
    Then, by Young's inequality, we have
    \begin{align}
    \begin{split}
        &\left|\int_\Omega\big(f^\prime(\varphi_1)-f^\prime(\varphi_2)\big)(\mu-\overline{\mu})\dx\right| \\
        &\leq C\|f^\prime(\varphi_1)-f^\prime(\varphi_2)\|_{L^2(\Omega)}^2 + \frac{1}{4}\|\mu-\overline{\mu}\|_{L^2(\Omega)}^2 \\
        &\leq C\Big(1 + \|f^\prime(\varphi_1)\|_{L^4(\Omega)}^{2} + \|f^\prime(\varphi_2)\|_{L^4(\Omega)}^{2}\Big)\|\nabla\varphi\|_{\mathbf{L}^2(\Omega)}^2 + \frac{1}{4}\|\mu-\overline{\mu}\|_{L^2(\Omega)}^2.
    \end{split}
    \end{align}
   Concerning the remaining terms, we get
    \begin{align*}
        &-\int_\Omega\big(f^\prime(\varphi_1)-f^\prime(\varphi_2)\big)\bv_1\cdot\nabla\varphi_1\dx
        +\int_\Omega\big(f^\prime(\varphi_1)-f^\prime(\varphi_2)\big)\bv_2\cdot\nabla\varphi_2\dx \\
        &= -\int_\Omega\big(f^\prime(\varphi_1)-f^\prime(\varphi_2)\big)\bv_2\cdot\nabla\varphi\dx -\int_\Omega\big(f^\prime(\varphi_1)-f^\prime(\varphi_2)\big)\bv\cdot\nabla\varphi_1\dx.
    \end{align*}
    Then, recalling \eqref{Est:WSU:6}, using Agmon's inequality once more and recalling that $\bv_2\in L^\infty(0,T_*;\bH^1(\Omega))$, we find
    \begin{align}
&\normmm{\int_\Omega\big(f^\prime(\varphi_1)-f^\prime(\varphi_2)\big)\bv_2\cdot\nabla\varphi\dx}\\
        &\leq \|f^\prime(\varphi_1)-f^\prime(\varphi_2)\|_{L^2(\Omega)}\|\nabla\varphi\|_{\mathbf{L}^2(\Omega)}\|\bv_2\|_{\mathbf{L}^\infty(\Omega)} \\
        &\leq C\Big(1 + \|f^\prime(\varphi_1)\|_{L^4(\Omega)} + \|f^\prime(\varphi_2)\|_{L^4(\Omega)}\Big)\|\bv_2\|_{\mathbf{H}^2(\Omega)}^\frac12\|\nabla\varphi\|_{\mathbf{L}^2(\Omega)}^2.
    \end{align}
    Analogously, using also Young's inequality, we infer that
 \begin{align}
    \begin{split}
    &\normmm{\int_\Omega\big(f^\prime(\varphi_1)-f^\prime(\varphi_2)\big)\bv\cdot\nabla\varphi_1\dx}
        \\&\leq \|f^\prime(\varphi_1)-f^\prime(\varphi_2)\|_{L^2(\Omega)}\|\nabla\varphi_1\|_{\mathbf{L}^6(\Omega)}\|\bv\|_{\mathbf{L}^3(\Omega)} \\
        &\leq C\norm{ \varphi_1}_{ H^2(\Omega)}\Big(1 + \|f^\prime(\varphi_1)\|_{L^4(\Omega)} + \|f^\prime(\varphi_2)\|_{L^4(\Omega)}\Big)\|\bv\|_{\mathbf{L}^2_\sigma(\Omega)}^\frac12\|D\bv\|_{\mathbf{L}^2(\Omega)}^\frac12\|\nabla\varphi\|_{\mathbf{L}^2(\Omega)}\\&
        \leq
        C\norm{\varphi_1}_{H^2(\Omega)}^\frac43\Big(1 + \|f^\prime(\varphi_1)\|_{L^4(\Omega)}^{4/3} + \|f^\prime(\varphi_2)\|_{L^4(\Omega)}^{4/3}\Big)(\|\bv\|_{\mathbf{L}^2_\sigma(\Omega)}^2+\|\nabla\varphi\|_{\mathbf{L}^2(\Omega)}^2)+\frac{\nu_*}{8}\|D\bv\|_{\mathbf{L}^2(\Omega)}^2
        \\&
        \leq C\Big(1+\norm{\varphi_1}_{H^2(\Omega)}^2+ \|f^\prime(\varphi_1)\|_{L^4(\Omega)}^4 + \|f^\prime(\varphi_2)\|_{L^4(\Omega)}^4\Big)(\|\bv\|_{\mathbf{L}^2_\sigma(\Omega)}^2+\|\nabla\varphi\|_{\mathbf{L}^2(\Omega)}^2)+\frac{\nu_*}{8}\|D\bv\|_{\mathbf{L}^2(\Omega)}^2
        .
    \end{split}
    \end{align}
    % For the last term, we use the inequalities of Hölder and Gagliardo--Nirenberg to deduce
    % \begin{align}
    %     \int_\Omega f^\prime(\varphi_2)\bv\cdot\nabla\varphi\dx
    %     &\leq C\|f^\prime(\varphi_2)\|_{L^6(\Omega)}\|\nabla\varphi\|_{\mathbf{L}^2(\Omega)}\|\bv\|_{\mathbf{L}^2_\sigma(\Omega)}^{1/2}\|D\bv\|_{\mathbf{L}^2(\Omega)}^{1/2} \\
    %     &\leq \frac{\nu_*}{8}\|D\bv\|_{\mathbf{L}^2(\Omega)}^2 + C\|\bv\|_{\mathbf{L}^2_\sigma(\Omega)}^2 + C\|f^\prime(\varphi_2)\|_{L^6(\Omega)}^2\|\nabla\varphi\|_{\mathbf{L}^2(\Omega)}^2.
    % \end{align}
    Eventually, combining all the estimates above, we arrive at
    \begin{align*}
        &\frac{1}{2}\Big(\|\bv(t)\|_{\mathbf{L}^2_\sigma(\Omega)}^2 + \|\nabla\varphi(t)\|_{\mathbf{L}^2(\Omega)}^2\Big) + \frac{1}{2}\int_0^t\|\mu-\overline{\mu}\|_{L^2(\Omega)}^2\d s + \frac{\nu_*}{4} \int_0^t\|D\bv\|_{\mathbf{L}^2(\Omega)}^2\d s
        \\&\leq \onehalf\norm{\bv_{0,1}-\bv_{0,2}}_{\bLds}^2+\onehalf\norm{\nabla\vphi_{0,1}-\nabla\vphi_{0,2}}_{\bLd}^2+\int_0^t\Lambda(s)\Big(\|\bv\|_{\mathbf{L}^2_\sigma(\Omega)}^2 + \|\nabla\varphi\|_{\mathbf{L}^2(\Omega)}^2\Big)\d s,
    \end{align*}
    where we set
    \begin{align*}
        \Lambda(t) := C\Big(1 + \|\bv_2(t)\|_{\mathbf{H}^2(\Omega)} +\|\mu_2(t)\|_{H^1(\Omega)}^2+\norm{\varphi_1(t)}_{H^2(\Omega)}^2+ \|f^\prime(\varphi_2(t))\|_{L^4(\Omega)}^4 + \|f^\prime(\varphi_1(t))\|_{L^4(\Omega)}^4\Big)
    \end{align*}
    for all $t\in[0,T_*)$.
    Thanks to the properties of $\bv_i$, together with  \eqref{L4}, it follows that $\Lambda\in~ L^1(0,T_*)$. Thus, applying the Gronwall inequality we find \eqref{contdep2}. This concludes the proof.

    \section{Proof of Theorem \ref{THM:Unique:2D}: Weak uniqueness in two dimensions}
    \label{weakuniq}
    We consider two sets of initial data $(\bv_{0,1},\varphi_{0,1})$ and $(\bv_{0,2},\varphi_{0,2})$ which comply with Theorem \ref{THM:WeakExistence} and such that $\overline{\vphi_{0,1}}=\overline{\vphi_{0,2}}$. For $i=1,2$, let $(\bv_i,\varphi_i,\mu_i)$ denote a weak solution corresponding to the initial data $(\bv_{0,i},\varphi_{0,i})$, whose existence is ensured by Remark \ref{2d}. We set
    $$
        (\bv_{0},\varphi_{0}) := (\bv_{0,1},\varphi_{0,1}) - (\bv_{0,2},\varphi_{0,2}), \quad (\bv,\varphi,\mu):= (\bv_1,\varphi_1,\mu_1) - (\bv_2,\varphi_2,\mu_2).
    $$
    Thus, $(\bv,\varphi,\mu)$ fulfills the following system of equations in a weak sense:
    \begin{align}
        \partial_t\bv + (\bv_1\cdot\nabla)\bv_1 - (\bv_2\cdot\nabla)\bv_2 - \text{div}(\nu(\varphi_1)D\bv_1) \nonumber\\
        +  \text{div}(\nu(\varphi_2)D\bv_2) + \nabla \pi &= \mu_1\nabla\varphi_1 -\mu_2\nabla\varphi_2 &\quad&\text{ in }\Omega_T, \label{WEq:1}\\
        \text{div}\ \bv &= 0 &\quad&\text{ in }\Omega_T, \label{WEq:2}\\
        \partial_t\varphi + \bv_1\cdot\nabla\varphi_1 - \bv_2\cdot\nabla\varphi_2 + \mu - \overline{\mu} &= 0 &\quad&\text{ in }\Omega_T, \label{WEq:3}\\
        \bv = \partial_\mathbf{n}\varphi &= 0 &\quad&\text{ on }\partial\Omega_T, \label{WEq:4}\\
        \bv(0) = \bv_0,\quad \varphi(0) &= \varphi_0 &\quad&\text{ in }\Omega. \label{WEq:5}
    \end{align}
     Note also that, by mass conservation, $\overline{\vphi}(t)=0$ for any $t\geq0$.
    We now test \eqref{WEq:1} with $A_S^{-1}\bv$ and \eqref{WEq:3} with $\varphi$ and add the resulting equations. Integrating by parts and using the identity
    \begin{align*}
        \frac{1}{2}\ddt\|\bv\|_{\mathbf{H}^1_\sigma(\Omega)^\prime}^2 = \frac{1}{2}\ddt\|\nabla A_S^{-1}\bv\|_{\mathbf{L}^2(\Omega)}^2 = \langle\partial_t\bv, A_S^{-1}\bv\rangle_{\mathbf{H}^1_\sigma(\Omega)',\bHus},
    \end{align*}
    we then infer that
    \begin{align}\label{Est:WU:1}
    \begin{split}
        &\frac{1}{2}\ddt\Big(\|\bv\|_{\mathbf{H}^1_\sigma(\Omega)^\prime}^2 + \|\varphi\|_{L^2(\Omega)}^2\Big) = \int_\Omega\big(\mu_1\nabla\varphi_1 - \mu_2\nabla\varphi_2\big)\cdot A_S^{-1}\bv \dx\\\
        &- \int_\Omega\big((\bv_1\cdot\nabla)\bv_1 - (\bv_2\cdot\nabla)\bv_2\big)\cdot A_S^{-1}\bv\dx -\int_\Omega\big(\nu(\varphi_1)D\bv_1 - \nu(\varphi_2)D\bv_2\big):DA_S^{-1}\bv\dx\\
        &- \int_\Omega(\bv_1\cdot\nabla\varphi_1 - \bv_2\cdot\nabla\varphi_2)\varphi\dx
        - \int_\Omega(\mu-\overline{\mu})\varphi\dx .
    \end{split}
    \end{align}
    Recalling that $\bv_i$ are divergence-free, integrating by parts, using Lady\v{z}enskaja's inequality and Young's inequality, we deduce, since $\bv_i\in L^\infty(0,T;\bLds)\cap  L^2(0,T;\bHus) \hookrightarrow L^4(0,T;\mathbf{L}^4(\Omega))$,
    \begin{align}\label{Est:WU:2}
    \begin{split}
        &\quad\left|\int_\Omega\big((\bv_1\cdot\nabla)\bv_1 - (\bv_2\cdot\nabla)\bv_2\big)\cdot A_S^{-1}\bv\dx\right| \\
        &\leq  \int_\Omega|(\bv_1\otimes\bv):\nabla A_S^{-1}\bv|\dx +\int_\Omega|(\bv\otimes\bv_2):\nabla A_S^{-1}\bv|\dx \\
        &\leq \|\bv_1\|_{\mathbf{L}^4(\Omega)}\|\nabla A_S^{-1}\bv\|_{\mathbf{L}^4(\Omega)}\|\bv\|_{\mathbf{L}_\sigma^2(\Omega)} + \|\bv\|_{\mathbf{L}_\sigma^2(\Omega)}\|\nabla A_S^{-1}\bv\|_{\mathbf{L}^4(\Omega)}\|\bv_2\|_{\mathbf{L}^4(\Omega)} \\
        &
        \leq C(\|\bv_1\|_{\mathbf{L}^4(\Omega)}+\|\bv_2\|_{\mathbf{L}^4(\Omega)})\norm{\nabla A_S^{-1}\bv}_{\mathbf{L}^2(\Omega)}^\frac12\norm{\bv}^\frac32_{\mathbf{L}^2(\Omega)}
        \\&
        \leq \frac{\nu_*}{28} \|\bv\|_{\mathbf{L}_\sigma^2(\Omega)}^2 + C(\|\bv_1\|_{\mathbf{L}^4(\Omega)}^4+\|\bv_2\|_{\mathbf{L}^4(\Omega)}^4)\|\bv\|_{\mathbf{H}^1_\sigma(\Omega)^\prime}^2.
        \end{split}
        \end{align}
    Next, we note that
    \begin{align*}
            \int_\Omega\big(\nu(\varphi_1)D\bv_1 - \nu(\varphi_2)D\bv_2\big):DA_S^{-1}\bv\dx
            &= \int_\Omega\big(\nu(\varphi_1)-\nu(\varphi_2)\big)D\bv_1:DA_S^{-1}\bv\dx \\
            &+ \int_\Omega\nu(\varphi_2)D\bv:DA_S^{-1}\bv\dx.
    \end{align*}
    Thanks to the Lipschitz continuity of $\nu$ and using again Lady\v{z}enskaja's inequality, we get
    \begin{align}\label{Est:WU:3}
        \begin{split}
            \left|\int_\Omega\big(\nu(\varphi_1)-\nu(\varphi_2)\big)D\bv_1:DA_S^{-1}\bv\dx\right|
            &\leq C\|\varphi\|_{L^4(\Omega)}\|D\bv_1\|_{\mathbf{L}^2(\Omega)}\|DA_S^{-1}\bv\|_{\mathbf{L}^4(\Omega)} \\
            &\leq C\|\varphi\|_{L^2(\Omega)}^{1/2}\|\varphi\|_{H^1(\Omega)}^{1/2}\|D\bv_1\|_{\mathbf{L}^2(\Omega)}\|\bv\|_{\mathbf{L}^2_\sigma(\Omega)}^{1/2}\|DA_S^{-1}\bv\|_{\mathbf{L}^2(\Omega)}^{1/2} \\
            &\leq \frac{1}{8}\|\nabla\varphi\|_{\mathbf{L}^2(\Omega)}^2 + \frac{\nu_*}{28}\|\bv\|_{\mathbf{L}_\sigma^2(\Omega)}^2 \\
            &\quad +C\|D\bv_1\|_{\mathbf{L}^2(\Omega)}^2\big(\|\varphi\|_{L^2(\Omega)}^2 + \|\bv\|_{\mathbf{H}^1_\sigma(\Omega)^\prime}^2\big).
        \end{split}
    \end{align}
    In the second term, we use integration by parts and we find
    \begin{align*}
            \int_\Omega\nu(\varphi_2)D\bv:DA_S^{-1}\bv\dx = -\int_\Omega\nu^\prime(\varphi_2)D A_S^{-1}\bv\nabla\varphi_2\cdot \,\bv\dx - \frac12\int_\Omega\nu(\varphi_2)\bv\cdot\Delta A_S^{-1}\bv\dx.
        \end{align*}
    By Lady\v{z}enskaja's inequality, we then infer that
        \begin{align}\label{Est:WU:4}
        \begin{split}
            \normmm{\int_\Omega\nu^\prime(\varphi_2)D A_S^{-1}\bv\nabla\varphi_2\cdot \,\bv\dx }
            &\leq C\|\nabla\varphi_2\|_{\mathbf{L}^4(\Omega)}\|\bv\|_{\mathbf{L}^2_\sigma(\Omega)}\|DA_S^{-1}\bv\|_{\mathbf{L}^4(\Omega)} \\&
            \leq C\|\nabla\varphi_2\|_{\mathbf{L}^4(\Omega)}\|\bv\|_{\mathbf{L}^2_\sigma(\Omega)}^\frac32\|DA_S^{-1}\bv\|_{\mathbf{L}^2(\Omega)}^\frac12 \\
            &\leq \frac{\nu_*}{28}\|\bv\|_{\mathbf{L}_\sigma^2(\Omega)}^2 + C\|\nabla\varphi_2\|_{\mathbf{L}^4(\Omega)}^4\|\bv\|_{\mathbf{H}^1_\sigma(\Omega)^\prime}^2.
        \end{split}
        \end{align}
    In the next step, we use regularity theory for the Stokes operator developed in \cite{GGGP}. In particular, there exists some $\pi\in C^0([0,T];H^1(\Omega))$ such that $-\Delta A_S^{-1}\bv + \nabla \pi = \bv$ holds almost everywhere in $\Omega$. This yields
        \begin{align*}
            -\frac12\int_\Omega\nu(\varphi_2)\bv\cdot\Delta A_S^{-1}\bv\dx
            &= \frac12\int_\Omega\nu(\varphi_2)|\bv|^2\dx - \frac12\int_\Omega\nu(\varphi_2)\bv\cdot\nabla \pi\dx \\
            &\geq \frac{\nu_*}{2}\|\bv\|_{\mathbf{L}_\sigma^2(\Omega)}^2 - \frac12\int_\Omega\nu(\varphi_2)\bv\cdot\nabla \pi\dx
        \end{align*}
    due to the properties of $\nu$. From the fundamental \cite[Lemma 3.1]{GGGP}, we have
    \begin{align*}
    \norm{\pi}_{L^4(\Omega)}\leq C\norm{\nabla A_S\bv}_{\mathbf L^2(\Omega)}^\frac12\norm{\bv}_{\mathbf L^2(\Omega)}^\frac12,
    \end{align*}
    and thus we obtain
    \begin{align}\label{Est:WU:6}
    \begin{split}
        &\left|-\frac12\int_\Omega\nu(\varphi_2)\bv\cdot\nabla \pi\dx\right| = \left|\frac12\int_\Omega\nu^\prime(\varphi_2)\nabla\varphi_2\cdot\bv \pi\dx\right|\\
        &\leq C\|\nabla\varphi_2\|_{\mathbf{L}^4(\Omega)}\|\bv\|_{\mathbf{L}^2_\sigma(\Omega)}\|\pi\|_{L^4(\Omega)} \\
        &\leq C\|\nabla\varphi_2\|_{\mathbf{L}^4(\Omega)}\|\bv\|_{\mathbf{L}^2_\sigma(\Omega)}^{3/2}\|\nabla A_S^{-1}\bv\|_{\mathbf{L}^2(\Omega)}^{1/2} \\
        &\leq \frac{\nu_*}{28}\|\bv\|_{\mathbf{L}_\sigma^2(\Omega)}^2 + C\|\nabla\varphi_2\|_{\mathbf{L}^4(\Omega)}^4\|\bv\|_{\mathbf{H}^1_\sigma(\Omega)^\prime}^2.
    \end{split}
    \end{align}
    In the Korteweg stress, thanks to the fact that $A_S^{-1}\bv$ is divergence-free, we recall the definition of $\mu_i$ to get
    \begin{align*}
        \int_\Omega\big(\mu_1\nabla\varphi_1 - \mu_2\nabla\varphi_2\big)\cdot A_S^{-1}\bv \dx = -\int_\Omega\text{div}(\nabla\varphi_1\otimes\nabla\varphi)\cdot A_S^{-1}\bv \dx - \int_\Omega\text{div}(\nabla\varphi\otimes\nabla\varphi_2)\cdot A_S^{-1}\bv \dx.
    \end{align*}
    Integration by parts and Lady\v{z}enskaja's inequality then yield
    \begin{align*}
        \left|\int_\Omega\text{div}(\nabla\varphi_1\otimes\nabla\varphi)\cdot A_S^{-1}\bv \dx\right|
        &= \left|-\int_\Omega\nabla\varphi_1\otimes\nabla\varphi:\nabla A_S^{-1}\bv \dx\right| \\
        &\leq C\|\nabla\varphi_1\|_{\mathbf{L}^4(\Omega)}\|\nabla\varphi\|_{\mathbf{L}^2(\Omega)}\|\nabla A_S^{-1}\bv\|_{\mathbf{L}^2(\Omega)}^\frac{1}{2}\|\bv\|_{\mathbf{L}^2(\Omega)}^\frac{1}{2} \\
        &\leq C\|\nabla\varphi_1\|_{\mathbf{L}^4(\Omega)}^4\|\nabla A_S^{-1}\bv\|_{\mathbf{L}^2(\Omega)}^2 + \frac{1}{4}\|\nabla\varphi\|_{\mathbf{L}^2(\Omega)}^2 + \frac{\nu_*}{28}\|\bv\|_{\mathbf{L}^2(\Omega)}^2.
    \end{align*}
    In a similar way, we can treat the second term. Indeed, we have
    \begin{align*}
        \left|\int_\Omega\text{div}(\nabla\varphi\otimes\nabla\varphi_2)\cdot A_S^{-1}\bv \dx\right|  \leq C\|\nabla\varphi_2\|_{\mathbf{L}^4(\Omega)}^4\|\nabla A_S^{-1}\bv\|_{\mathbf{L}^2(\Omega)}^2 + \frac{1}{4}\|\nabla\varphi\|_{\mathbf{L}^2(\Omega)}^2 + \frac{\nu_*}{28}\|\bv\|_{\mathbf{L}^2(\Omega)}^2.
    \end{align*}
    For the convective term in the Allen--Cahn equation, it holds
    \begin{align*}
        \int_\Omega(\bv_1\cdot\nabla\varphi_1 - \bv_2\cdot\nabla\varphi_2)\varphi\dx = \int_\Omega(\bv\cdot\nabla\varphi_1)\varphi\dx + \int_\Omega(\bv_2\cdot\nabla\varphi)\varphi\dx.
    \end{align*}
    Since $\bv_2$ is divergence-free, integration by parts yields
    \begin{align*}
        \int_\Omega(\bv_2\cdot\nabla\varphi)\varphi\dx = 0.
    \end{align*}
    For the remaining term, we obtain
    \begin{align}\label{Est:WU:9}
    \begin{split}
        \left|\int_\Omega(\bv\cdot\nabla\varphi_1)\varphi\dx\right|
        &\leq \|\bv\|_{\mathbf{L}^2_\sigma(\Omega)}\|\nabla\varphi_1\|_{\mathbf{L}^4(\Omega)}\|\varphi\|_{L^4(\Omega)} \\
        &\leq C\|\bv\|_{\mathbf{L}^2_\sigma(\Omega)}\|\nabla\varphi_1\|_{\mathbf{L}^4(\Omega)}\|\varphi\|_{L^2(\Omega)}^{1/2}\|\varphi\|_{H^1(\Omega)}^{1/2} \\
        &\leq \frac{\nu_*}{28}\|\bv\|_{\mathbf{L}^2_\sigma(\Omega)}^2 + \frac{1}{8}\|\nabla\varphi\|_{\mathbf{L}^2(\Omega)}^2 + C\|\nabla\varphi_1\|_{\mathbf{L}^4(\Omega)}^4\|\varphi\|_{L^2(\Omega)}^2.
    \end{split}
    \end{align}
    Since $\overline{\varphi}=0$, Fubini's theorem and the definition of $\mu$ imply
    \begin{align*}
        \int_\Omega(\mu-\overline{\mu})\varphi\dx = \int_\Omega\mu\,\varphi\dx = \|\nabla\varphi\|_{\mathbf{L}^2(\Omega)}^2 + \int_\Omega\big(f^\prime(\varphi_1)-f^\prime(\varphi_2)\big)\varphi\dx.
    \end{align*}
    Finally, recalling \ref{ASS:S1}, we observe that
    \begin{align}\int_\Omega\big(F^\prime(\varphi_1)-F^\prime(\varphi_2)\big)\varphi\dx \geq \theta\|\varphi\|_{L^2(\Omega)}^2,
    \end{align}
    and
    \begin{align}
        \label{Est:WU:10}
        \normmm{\int_\Omega \frac{\theta_0}2\varphi^2\dx}\leq C\norm{\varphi}_{L^2(\Omega)}^2.
    \end{align}
    Combining estimates \eqref{Est:WU:1}--\eqref{Est:WU:10}, we get
    \begin{align*}
        &\frac{1}{2}\ddt\Big(\|\bv\|_{\mathbf{H}^1_\sigma(\Omega)^\prime}^2 + \|\varphi\|_{L^2(\Omega)}^2\Big) + \frac{1}{4}\|\nabla\varphi\|_{\mathbf{L}^2(\Omega)}^2 + \frac{\nu_*}{4}\|\bv\|_{\mathbf{L}^2_\sigma(\Omega)}^2
        \leq \Lambda_1\Big(\|\bv\|_{\mathbf{H}^1_\sigma(\Omega)^\prime}^2 + \|\varphi\|_{L^2(\Omega)}^2\Big),
    \end{align*}
    where
    \begin{align*}
        \Lambda_1 := C\Big(\|\bv_1\|_{\mathbf{L}^4(\Omega)}^4 + \|\bv_2\|_{\mathbf{L}^4(\Omega)}^4 + \|D\bv_1\|_{\mathbf{L}^2(\Omega)}^2 + \|\nabla\varphi_1\|_{\mathbf{L}^4(\Omega)}^4 + \|\nabla\varphi_2\|_{\mathbf{L}^4(\Omega)}^4\Big).
    \end{align*}
    Since $H^1(\Omega)\hookrightarrow L^4(\Omega)$ and $\vphi_i\in L^\infty(0,T;H^1(\Omega))$, $i=1,2$, as well as $\bv_i\in L^4(0,T;\mathbf L^4(\Omega))$, $i=1,2$, we have $\Lambda_1\in L^1(0,T)$, Gronwall's lemma gives the proof.

\section{Proof of Lemma \ref{convaaa}}
\label{proofequil}
The proof is divided into three steps.

{\bf (i)}.
First we notice that in this case, since we do not have a separate energy inequality only concerning the kinetic energy (see Remark \ref{nosplit}), it is nontrivial to show in the first step that $\bv(t)\to\mathbf 0$ in $\bLds$ as $t\to\infty$, like in \cite{GP}. Therefore, we will need a different approach to obtain this convergence at the very end of the proof of Theorem \ref{finn}.

Let us now consider a sequence $t_n\to \infty$ such that $\bv(t_n)\rightharpoonup  \bv_\infty$ weakly in $\bLd$ and $\varphi(t_n)\rightharpoonup {\varphi}_\infty$ weakly in $H^1(\Omega)$,  with $(\bv_\infty,\vphi_\infty)\in \omega(\bv,\varphi)$. We then define the sequence of trajectories $\bv_n(t):=\bv(t+t_n)$, $\varphi_n(t):=\varphi(t+t_n)$ and $\mu_n(t):=\mu(t+t_n)$. Thanks to Corollary \ref{thmglobal}, recalling the energy inequality \eqref{energyineq}, we get that ${{E}}(\bv(t_n),\vphi(t_n))\leq {{E}}({\bv_0,\vphi_0})$ for any $n$, and thus, for any $T>0$, uniformly in $n$,
				\begin{align*}
					&\| \partial_t\bv_n\|_{L^\frac43(0,T;\mathbf H^1_\sigma(\Omega)')}+\| \bv_n\|_{L^\infty(0, T;\bLd)}+\| \bv_n\|_{L^2(0,T;\bH_\sigma^1(\Omega))}\\&+\| \partial_t\vphi_n\|_{L^\frac43(0,T;\Ld)}+\| \vphi_n\|_{L^2(0,T;H^2(\Omega))}+\| \vphi_n\|_{L^\infty(0,T;H^1(\Omega))}+\| \mu_n\|_{L^2(0,T;L^2(\Omega))}\leq C(T).
				\end{align*}
				From the estimates above, we deduce that there exists $(\bv^*,\vphi^*)$ (and thus $\mu^*$) such that, for any fixed $T>0$,
				\begin{align}
		\nonumber	&\bv_n\rightharpoonup\bv^*\quad\text{weakly in } L^2(0,T;\bHus)\cap W^{1,\frac43}(0,T;\bHus'),\\	&
        \bv_n\to\bv^*\quad\text{strongly in } L^2(0,T;\bLds)\label{essenziale0}\\&
        \label{est2b}\vphi_n\rightharpoonup \vphi^*\quad\text{weakly in } L^2(0,T;H^2(\Omega))\cap H^1(0,T;\Hu'),\\&
					\vphi_n\overset{\ast}{\rightharpoonup} \vphi^*\quad\text{weakly-$*$ in } L^\infty(\Omega\times(0,T)),
					\\&
					\vphi_n\to \vphi^*\quad \text{in } L^2(0,T;H^s(\Omega)),\, \forall\: s\in[0,2),\label{essenziale}\\&
					\mu_n\rightharpoonup  \mu^*\quad\text{weakly in } L^2(0,T;\Ld).
				\end{align}
				These convergences are of course enough to pass to the limit in the equations for $(\bv_n,\vphi_n)$. In particular, the limit pair $(\vphi^*,\mu^*)$ satisfies, for any $T>0$, the equations
				\begin{align*}
					\langle \partial_t\vphi^*,v\rangle_{\Hu',\Hu}-(\bv^*\cdot \vphi^*,\nabla \bv)+(\mu^*-\overline {\mu^*}, v)=0,\quad \forall \: v\in \Hu,\quad \text{a.e. in }(0,T),\\
					\mu^*=-\Delta\vphi^*+f'(\vphi^*),\quad\text{a.e. in }\Omega\times(0,T),
				\end{align*}
				with initial datum $\vphi^*(0)={\vphi_\infty}$. This follows immediately from the fact that $\vphi_n(0)=\vphi(t_n){\rightharpoonup} {\vphi_\infty}$ weakly in $\Hu$. Furthermore, we have \[
				\lim_{n\to \infty}{{E}}(\bv_n(t),\vphi_n(t))={{E}}(\bv^*(t),\vphi^*(t))
				\] for {almost any} $t\geq 0$. By the energy inequality \eqref{energyineq}, we infer that the energy ${{E}}(\bv(\cdot),\vphi(\cdot))$ is nonincreasing, thus there exists ${{E}}_\infty$ such that
                \begin{align}
				\lim_{t\to+ \infty}{{E}}(\bv(t),\vphi(t))={{E}}_\infty.
				\label{ene}
                \end{align}
				 Hence, {for almost any $t\geq0$}, we have
				$${{E}}(\bv^*(t),\vphi^*(t))=\lim_{n\to \infty}{E}(\bv_n(t),\vphi_n(t))=\lim_{n\to \infty}{{E}}(\bv(t+t_n),\vphi(t+t_n))={{E}}_\infty,$$entailing that $E(\bv^*(\cdot),\vphi^*(\cdot))$ is equal to $ E_\infty$ almost everywhere in time. Passing then to the limit in the energy inequality, which is valid for each $(\bv_n,\vphi_n)$ thanks again to \eqref{energyineq}, we obtain
			\begin{align}\label{zeros}{{E}}_\infty+\nu_*\int_s^t\norm{D\bv^*}_{\bLd}^2\d \tau+\int_s^{t}\|\mu^*(\tau)-\overline{\mu^*(t)}\|_{\bLd}^2 \: \d \tau\leq {{E}}_\infty\quad \text{ for a.a. } 0\leq s\leq t<\infty,\end{align}
            with $s=0$ included.
            This first entails that $D\bv^*=\mathbf 0$  almost everywhere in $\Omega$, and thus, by Korn's inequality and the boundary condition $\bv^*=\mathbf 0$ on $\partial\Omega\times(0,\infty)$, gives $\bv^*\equiv \mathbf 0$ for any $t\geq0$. Since it holds $\bv_n(0)=\bv(t_n)\rightharpoonup \bv_\infty$ in $\bLd$ as $n\to\infty$, we have $\bv^*(0)=\bv_\infty$. Thus $\bv_\infty=\mathbf 0$ almost everywhere in $\Omega$.  Additionally, \eqref{zeros} entails $\mu^*=\overline{\mu^*}$ almost everywhere in $\Omega$. By comparison, it also holds $\partial_t\vphi^*=0$ in $\Hu'$, for almost every $t\geq 0$. Therefore, we infer that $$\vphi^*(t)=\int_0^t\partial_t\vphi^*(s) \: \d s+{\vphi_\infty}={\vphi_\infty}$$ in $\Hu'$, for all $t\geq 0$. Therefore, ${\vphi_\infty}$ satisfies \eqref{conv1t} for some constant $\mu_\infty\in \R$, and then $(\mathbf{0},{\vphi_\infty})\in \mathcal{S}$. This gives $\omega(\bv,\vphi)\subset \mathcal{S}$.
    %     Now we can prove

				% Furthermore, convergence \eqref{essenziale} implies that, up to subsequences,
				% \begin{align}
				% 	\vphi(t+t_n)\to {\vphi_\infty} \text{ strongly in }\Hu
				% 	\label{ttt}
				% \end{align}
				% for almost any $t\in(0,\infty)$, from which we deduce \eqref{ttt1}.
We can now show the uniform strict separation properties of the $\vphi$ projection of the $\omega$-limit. This is indeed known in the literature (see, e.g., \cite{GGPJ, AW}). It is enough to show that, given $(\mathbf 0,\varphi_\infty)\in \mathcal S$, there exists $\delta_{\vphi_\infty}>0$, possibly depending on $\varphi_\infty$ such that
\begin{align*}
    \norm{\vphi_\infty}_{L^\infty(\Omega)}\leq 1-\delta_{\vphi_\infty},
\end{align*}
which is trivially seen from \eqref{conv1t}, since we easily get $\norm{F'(\vphi_\infty)}_{L^\infty(\Omega)}\leq C(1+\normmm{\mu_\infty})$. Then, since $\omega(\bv,\vphi)\subset \mathcal S$, the same property holds for any element of the $\omega$-limit. To prove that the separation property is actually uniform over  ${\bf P}_2\omega(\bv,\vphi)$, where ${\bf P}_2$ is the projection over the second component, first need to show that  ${\bf P}_2\omega(\bv,\vphi)$ is bounded in $H^2(\Omega)$. Now, by the energy inequality \eqref{energyineq} it is immediate to deduce that there exists $C>0$, only depending on the initial data, such that
\begin{align}
\norm{\nabla \varphi_\infty}_{\bLd}\leq \sup_{t\geq0}\norm{\nabla \vphi(t)}_{\bLd}\leq C(E(\bv_0,\vphi_0)),\quad \forall \varphi_\infty\in {\bf P}_2\omega(\bv,\vphi).\label{f}
\end{align}
Then, given $\varphi_\infty\in {\bf P}_2\omega(\bv,\vphi)$, by multiplying \eqref{conv1t} by $-\Delta\varphi$ and integrating over $\Omega$ we get, after an integration by parts,
\begin{align*}
    \norm{\Delta \varphi_\infty}_{\Ld}^2+\int_\Omega F''(\varphi_\infty)\normmm{\nabla \varphi_\infty}^2\dx=\theta_0\int_\Omega \normmm{\nabla \varphi_\infty}^2\dx +\int_\Omega \nabla \mu_\infty\cdot \nabla \varphi_\infty\dx ,
\end{align*}
so that, since $\mu_\infty$ is a constant and $F''\geq 0$, we have
$$
\norm{\Delta \varphi_\infty}_{\Ld}^2\leq \theta_0 \normm{\nabla \varphi_\infty}^2_{\bLd}\leq \theta_0C^2,\quad \forall \varphi_\infty\in {\bf P}_2\omega(\bv,\vphi),
$$
where we used \eqref{f}. Since, by mass conservation, $\overline {\vphi}_\infty=\overline{\vphi}_0$ for any $\varphi_\infty\in {\bf P}_2\omega(\bv,\vphi)$, this estimate ensures that ${\bf P}_2\omega(\bv,\vphi)$ is bounded in $H^2(\Omega)$. Since $H^2(\Omega)\hookrightarrow \hookrightarrow C^\alpha(\overline\Omega)$ for some $\alpha\in(0,1)$, by a simple contradiction argument we get \eqref{sepaglobal} (see, for instance, \cite{GGPJ}).

{\bf (ii)}. We do not discuss the validity of the asymptotic strict separation property \eqref{asympt}, as this can be obtained by means of De Giorgi's iterations arguing word by word as in \cite[Lemma 3.11]{GP2}. Indeed, the additional divergence-free advective term does not play any role in De Giorgi's iterations, thanks to \eqref{vanish}.
% As this proof is quite technical, for the sake of readability we postpone it to Section \ref{asymptsep}.

{\bf (iii)}. We need to finally prove the further characterization of the $\omega$-limit, together with \eqref{convergence1b}. But the former is a trivial consequence of convergences \eqref{essenziale0}-\eqref{essenziale} (see also \cite{GP2}), whereas the latter can be shown by contradiction exploiting precompactness of trajectories (see, e.g., \cite{GP,GP2}), after recalling that $\vphi\in BC_w([0,\infty); H^1(\Omega))$ and $\bv\in BC_w([0,\infty);\mathbf L^2_\sigma(\Omega))$.

\section{Proof of Lemma \ref{twoparts1}: Properties of good equilibrium points}
\label{goodeq}
The proof is very similar to the one of \cite[Lemma 3.13]{GP2}, after an adaptation to include the velocity term. Recall that, for fixed $T>0$, we have
\begin{align}
    \sup_{t\in A_M(T)}\left(\nu_*\norm{D\bv(t)}^2_{\bLd}+\norm{\mu(t)-\overline{\mu}(t)}_{L^2(\Omega)}^2\right)\leq M^2.\label{M1b}
\end{align}
Then, following \cite{GP2}, since the velocity $\bv$ does not appear in the equation for $\mu$, we find
\begin{align*}
          \normmm{\int_\Omega{\mu(t)}\dx}\leq C\int_\Omega\normmm{f'(\vphi(t))}\dx\leq C(1+\norm{\mu(t)-\overline\mu(t)}_{L^2(\Omega)}),
                \end{align*}
which entails that (see \eqref{M1b})
\begin{align}
    \sup_{t\in A_M(T)}\norm{\mu(t)}_{L^2(\Omega)}\leq C(M).\label{M2b}
\end{align}
We can also immediately show that
\begin{align}
    \sup_{t\in A_M(T)}\norm{F'(\vphi(t))}_{L^2(\Omega)}\leq C(M),\label{M2bA}
\end{align}
and, by elliptic regularity, 
\begin{align}
   \sup_{t\in A_M(T)} \norm{\vphi(t)}_{H^2(\Omega)}\leq C(M),\quad \forall t\in A_M(T),\label{prb}
\end{align}
and for any $T>0$. Concerning the velocity, we have by Korn's inequality
\begin{align}
   \sup_{t\in A_M(T)} \norm{\bv(t)}_{\bH^1_\sigma(\Omega)}\leq C(M),\quad \forall t\in A_M(T),\label{velox}
\end{align}
Therefore, a standard contradiction argument and the embeddings $\bH^1_\sigma(\Omega)\hookrightarrow\hookrightarrow \bH^s(\Omega)$, $s\in(0,1)$, and $H^2(\Omega)\hookrightarrow\hookrightarrow H^r(\Omega)$, $r\in(0,2)$, yield \eqref{precomp1a} as well as the compactness properties of $\omega_g(\bv,\vphi)$. This ends the proof.

\section{Proof of Theorem \ref{uniqueeq}}
\label{uniqeq}
We first recall the \L ojasiewicz-Simon inequality, which holds thanks to \cite[Proposition 6.1]{AW}:
	\begin{proposition}
    \label{LSineq}
				 Assume that $F$ is additionally real analytic in $(-1,1)$ {and let } $\vphi\in H^2(\Omega)\cap \mathcal H_m$ {be } such that $-1+\gamma\leq \vphi(x)\leq 1-\gamma$, for {almost }any $x\in \Omega$ {and } for some $\gamma\in(0,1)$. {Furthermore, let } $\vphi_\infty\in {\bf P}_2\mathcal{S}$ {be fixed}. Then there exist $\vartheta\in \left(0,\frac{1}{2}\right]$, $\eta>0$ and a positive constant $C$ such that
				\begin{align}
					\vert{E_{free}}(\vphi)-{E_{free}}(\vphi_\infty)\vert^{1-\vartheta}\leq C\| \delta E_{free}(\varphi)\|_{\Hu'},
					\label{ener}
				\end{align}
				{provided that } $\|\vphi-\vphi_\infty\|_{\Hu}\leq \eta$,
				\label{Lojaw}
where $\delta E:H^1_{(m)}(\Omega)\to H^1_{(0)}(\Omega)'$ is the Frechét derivative of $E_{free}:H^1_{(m)}(\Omega)\to \R$.
\end{proposition}

The proof is now very similar to \cite[Theorem 3.15]{GP2}, so that we only give the main details. Let us choose $\tilde \gamma$ coinciding with the value of $\delta$ in \eqref{asympt}, so that, as $\tilde \gamma\leq\delta_1$, it holds (see \eqref{sepaglobal}) $-1+\tilde\gamma\leq \vphi_\infty\leq 1-\tilde\gamma$ in $\Omega$, for any $\vphi_\infty\in\omega_g(\bv,\vphi)$. Furthermore, for any $\vphi_{\infty,m}\in \omega_g(\bv,\vphi)$ we can find $\vartheta_m\in \left(0,\frac{1}{2}\right]$ and $\eta_m>0$, given {by }Proposition \ref{Lojaw}, for which \eqref{ener} is valid with constant $C_m$. Now, from Lemma \ref{twoparts1} we get that $\omega_g(\bv,\vphi)\subset \omega(\bv,\vphi)$ is compact in $\{\mathbf 0\}\times H^1(\Omega)$, so that we can find a finite family of $H^1(\Omega)$-open balls, say $\{B_{\eta_m}\}_{m=1}^{M_1}$, centered at $\{\vphi_{\infty,m}\}_{m=1}^{M_1}\subset {\bf P}_2\omega_g(\bv,\vphi)$ (recall that ${\bf P}_2$ is the projection over the second component) and with radii $\eta_m$ (depending on the center $\vphi_{m,\infty}\in {\bf P}_2\omega_g(\bv,\vphi)$), such that
				\begin{equation*}
					\bigcup_{ \varphi_\infty \in {\bf P}_2\omega_{g}(\bv,\varphi)} \{\varphi_\infty\} \subset W:= \bigcup_{m=1}^{M_{1}}B_{\eta_m}.
				\end{equation*}
	Recalling \eqref{E2}, the energy functional $E(\cdot,\cdot)$ is constant over $\omega(\bv,\vphi)$. Also, since {the centers  $\{\vphi_m\}_{m = 1}^{M_1}$ are in finite number, }we can infer that \eqref{ener} holds \textit{uniformly}, for any $\vphi\in W$ such that $\|\vphi\|_{L^\infty(\Omega)}\leq 1-\tilde\gamma$, and we can replace $E_{free}(\vphi_\infty)$ with $E_{\infty}$.

    Then, by  \eqref{asympt} and \eqref{precomp1a}, there exist ${t_*}>0$ such that $\vphi(t)\in W$ for any good time $t\in A_M(t_*)$ and the uniform strict separation property holds for any $t\geq t_*$, i.e.,
\begin{align*}
                    \sup_{t\geq t*}\norm{\vphi(t)}_{L^\infty(\Omega)}\leq 1-\tilde \gamma.
                \end{align*}
               Therefore, we get from \eqref{ener} that, for any $t\in A_M(t_*)$,
					\begin{align}
						\label{pp}\left(E_{free}(\vphi(t))-{{E}}_\infty\right)^{1-\vartheta}\leq C\| \mu(t)-\overline {\mu(t)}\|_{\Ld},
					\end{align}
                    where we can choose $\vartheta\in(0,\frac12)$, as we have $\sup_{t\geq 0}E(\bv(t),\varphi(t))<\infty$.
Fix now $s,t\geq t^*$, $0\leq s\leq t$. Then we have, by the energy inequality \eqref{energyineq},
\begin{align*}
   &\nonumber \int_s^t \norm{\mu(\tau)-\overline{\mu}(\tau)}_{\Ld}^2\dtau+ \int_s^t\norm{\sqrt{\nu(\vphi(\tau))}D\bv(\tau)}_{\bLd)}^2\dtau\\&
   \leq E(\bv(s),\vphi(s))-E(\bv(t),\vphi(t)).
\end{align*}
This gives
\begin{align}
   &\nonumber \left(\int_s^t \norm{\mu(\tau)-\overline{\mu}(\tau)}_{\Ld}^2\dtau+ \int_s^t\norm{\sqrt{\nu(\vphi(\tau))}D\bv(\tau)}_{\bLd)}^2\dtau\right)^{2(1-\vartheta)}\\&
   \leq (E(\bv(s),\vphi(s))-E(\bv(t),\vphi(t)))^{2(1-\vartheta)}.\label{cv}
\end{align}
% Now we observe that
% \begin{align*}
%    \lim_{t\to\infty} \normmm{\frac12 \norm{\bv(s)}^2_{\bLd}- \frac12 \norm{\bv(t)}^2_{\bLd}}\leq \frac12\normmm{\bv(s)}^2_{\bLd},
% \end{align*}
% by \eqref{veloxconv}.
We now let $t\to\infty$ in \eqref{cv}, and obtain, recalling that $E(\bv(t),\vphi(t))\to E_\infty$ as $t\to \infty$,
\begin{align}
   &\nonumber \left(\int_s^\infty \norm{\mu(\tau)-\overline{\mu}(\tau)}_{\Ld}^2\dtau+ \int_s^\infty\norm{\sqrt{\nu(\vphi(\tau))}D\bv(\tau)}_{\bLd)}^2\dtau\right)^{2(1-\vartheta)}\\&
   \leq (E(\bv(s),\vphi(s))-E_\infty)^{2(1-\vartheta)}\nonumber\\&\leq C\normmm{E_{free}(\bv(s),\vphi(s))-E_\infty}^{2(1-\vartheta)}+{C}\norm{\bv(s)}^{4(1-\vartheta)}_{\bLd}\nonumber\\&
 \leq  C\normmm{E_{free}(\bv(s),\vphi(s))-E_\infty}^{2(1-\vartheta)} (\chi_{A_M(t_*)}(s)\nonumber+\chi_{(t_*,\infty)\setminus A_M(t_*)}(s))\\&\quad+C\norm{\sqrt{\nu(\vphi(s))}D\bv(s)}_{\bLd)}^2
   ,\label{cv2}
\end{align}
for some $C>0$, where in the last inequality we used Korn's inequality, and the bound $\sup_{t\geq 0}\norm{\bv(t)}_{\bLd}\leq C_1$, for some $C_1>0$, and $4(1-\vartheta)>2$.

Now, recalling that $E_{free}(\vphi(t))\leq E(\bv_0,\vphi_0)$ for any $t\geq0$, we infer, since $0<\nu_*\leq \nu(\cdot)$ by assumption,
					\begin{align*}&\normmm{E(\vphi(s))-E_\infty}^{2(1-\vartheta)}\chi_{(t_*,\infty)\setminus A_M(t_*)}(s)\\&\leq (2E(\bv_0,\vphi_0))^{2(1-\vartheta)}\frac{\nu_*\norm{D\bv(s)}_{\bLd}^2+\norm{\mu(s)-\overline{\mu(s)}}_{L^2(\Omega)}^2}{M^2}\chi_{(t_*,\infty)\setminus A_M(t_*)}(s)\\&
                    \leq (2E(\bv_0,\vphi_0))^{2(1-\vartheta)}\frac{\norm{\sqrt{\nu(\vphi(s))}D\bv(s)}_{\bLd}^2+\norm{\mu(s)-\overline{\mu(s)}}_{L^2(\Omega)}^2}{M^2}\chi_{(t_*,\infty)\setminus A_M(t_*)}(s),
				\end{align*}
				for almost any $s\in (t_*,\infty)\setminus A_M(t_*)$.
				On the other hand, in the set of good times, thanks to \eqref{pp}, we have
					$$
				\normmm{E_{free}(\vphi(s))-E_\infty}^{2(1-\vartheta)}\chi_{ A_M(t_*)}(s)\leq C^2{\norm{\mu(s)-\overline{\mu(s)}}_{\bLd}^2}\chi_{A_M(t_*)}(s),
				$$
				for almost any $s\in  A_M(t_*)$.
			Putting everything together, we obtain from \eqref{cv2} that 
\begin{align}
   &\nonumber \left(\int_s^\infty \norm{\mu(\tau)-\overline{\mu}(\tau)}_{\Ld}^2\dtau+ \int_s^\infty\norm{\sqrt{\nu(\vphi(\tau))}D\bv(\tau)}_{\bLd)}^2\dtau\right)^{2(1-\vartheta)}\\&
   \leq C_2\left(\norm{\mu(s)-\overline\mu(s)}^2+\norm{\sqrt{\nu(\vphi(s))}D\bv(s)}_{\bLd)}^2\right),\label{cv2b}
\end{align}
    where $C_2>0$ is a suitable constant.  Applying \cite[Lemma 7.1]{FS} to \eqref{cv2b}, where, using the notation of the lemma, we set
$$
Z(\cdot)=\sqrt{\norm{\sqrt{\nu(\vphi(\tau))}D\bv(\tau)}_{\bLd)}^2+\norm{\mu(\cdot)-\overline{\mu(\cdot)}}_{L^2(\Omega)}^2}, \quad \alpha=2(1-\vartheta)\in(1,2), \quad \zeta=C_2>0,
$$ 
and $\mathcal M=(t_*,\infty)$, we find, using also Korn's inequality,
\begin{align}
  \mu-\overline\mu\in L^1(t^*,\infty;L^2(\Omega)) , \quad \bv\in L^1(t^*,\infty;\mathbf H^1_\sigma(\Omega)). \label{bA}
\end{align}
Moreover, recalling that $\vert\vphi\vert\leq 1$ and $\bv\in L^\infty(0,\infty;\bHus)$, we get
			$$
		\left\vert\int_\Omega \vphi\nabla\zeta\cdot \bv \: \d x\right\vert\leq \| \bv\|_{\bLd} \| \zeta\|_{\Hu},\quad \forall \zeta\in \Hu,
			$$
			so that, by \eqref{bA},
			$$
			\| \bv\cdot \nabla \vphi\|_{L^1(t^*;\infty;\Hu')} \leq C\| D\bv\|_{L^1(t^*;\infty;\bLd)}\leq C.
			$$
			By comparison in the conserved advective Allen-Cahn equation, we then deduce that $\partial_t\vphi\in L^1(t^*,\infty;\Hu')$. As a consequence, we have
            $$
			\vphi(t)=\vphi(t^*)+ \int_{t^*}^t\partial_t\vphi(\tau) \: \d\tau \to {\vphi_\infty}\quad \text{ in }\Hu', \text{ as } t \to \infty,$$for some ${\vphi_\infty}\in \Hu'$. Therefore, $\vphi(t)$ converges in $\Hu'$ as $t\to\infty$ and, by uniqueness of the limit in $\Hu'$, we conclude that $\omega(\bv,\vphi)$ is a singleton, say $\omega(\bv,\vphi)=\{(\textbf{0},{\vphi_\infty})\}$. The convergences \eqref{equil} are then a consequence of \eqref{convergence1b}.
\bigskip

\section{Proof of Theorem \ref{finn}: Asymptotic regularization of weak solutions}
\label{asymreg}
 Here we show that any weak solution asymptotically regularizes, if $\nu\in W^{1,\infty}(\R)$.  
 % the trajectories are actually precompact in $\mathbf H^s(\Omega)\times H^q(\Omega)$ for $s\in(0,1)$, $q\in[1,2)$. 
 %First, from
%  From we have that $\varphi(t_n)\rightharpoonup \vphi_\infty$ weakly in $H^1(\Omega)$, and $\bv(t_n)\rightharpoonup \mathbf 0$ weakly in $\bLds$, where $(\mathbf 0,\vphi_\infty)$ is an equilibrium point according to Definition \ref{defequil}. 
Again we will make use of the notion of good and bad times. 
Let us now consider an initial datum $(\bv_0,\vphi_0)\in \bLds\times \mathcal H_m$ and set $C_0:=E(\bv_0,\varphi_0)$.
First, we observe that by Bihari's Lemma used in the proof of Theorem \ref{THM:StrongWP}, given an initial datum $(\tilde\bv_0,\tilde\vphi_0)$, the maximal time of existence of the strong solution depends only on the $L^2$-norms of the initial data $D\tilde\bv_0$, $\tilde\mu_0-\overline{\tilde\mu_0}$ and the initial energy $E(\tilde\bv_0,\tilde\varphi_0)$ (see estimate \eqref{maxtime}). Let us then fix $\epsilon>0$ and consider any initial datum $(\tilde\bv_0,\tilde\varphi_0)$, satisfying $E(\tilde\bv_0,\tilde\varphi_0)\leq C_0$ and
                \begin{align}
                 \norm{D\tilde\bv_0}_{\bLd}^2+\norm{\tilde\mu_0-\overline{\tilde\mu}_0}^2_{\Ld}\leq \epsilon.
                    \label{initial}
                \end{align}
                As a consequence, for any such initial data there exists a unique strong solution, and the maximal time of existence $T^*$ only depends on $\epsilon>0$ (and of course on $C_0$). Not only this, but the (higher-order) norms of the solutions only depend on $\epsilon$ and $C_0$, as it is seen from \eqref{max}.  We fix the time horizon as $T_M=T_M(\epsilon,C_0)\in (0,T^*)$, which is assumed as strictly less than the maximal time of existence, whereas we denote the maximal (positive) constant controlling all higher-order norms on $(0,T_M)$ by $C_M=C_M(\epsilon,C_0)$. Namely, \textit{for any} strong solution $(\tilde\bv,\tilde\varphi)$ departing from an initial datum in the class above it holds
                \begin{align}
&\nonumber\norm{\tilde\bv}_{H^1(0,T_M;\bLds)\cap L^2(0,T_M;\mathbf H^2(\Omega)
)\cap L^\infty(0,T_M;\bH^1_\sigma(\Omega)))}\\&
+\norm{\tilde\varphi}_{H^1(0,T_M;H^1(\Omega))\cap L^\infty(0,T_M;H^2(\Omega))}+\norm{\tilde\mu}_{L^\infty(0,T_M;L^2(\Omega))\cap L^2(0,T_M;H^1(\Omega))}\leq C_M.
                    \label{holdercontrol}
                \end{align}
Let us consider a general time $s>0$. For $r>0$ and $t_0>0$, given the ball $\widehat B_{r}(t_0):=\{t\in [0,\infty):\ \normmm{t-t_0}\leq r\}$, setting $r=\frac {T_M}4$, we introduce the set of bad times
$$
C_s:=\{t\in [0,\infty):\ \norm{D\bv(t)}_{\bLd}^2+\norm{\mu(t)-\overline\mu(t)}_{\Ld}^2\geq \epsilon\}\cap \widehat B_r(s).
$$
Now, the energy inequality \eqref{ineq1} gives
\begin{align}
\int_0^\infty \norm{D\bv}_{\bLd}^2\d t+\int_0^\infty \norm{\mu-\overline\mu}_{\Ld}^2\d t\leq C
    \label{ctrl1}.
\end{align}
Thus, we can find $s_0=s_0(\epsilon,C_0)>0$ such that 
\begin{align}
    \normmm{C_s}\leq \frac1\epsilon\int_{s-r}^\infty \norm{D\bv}_{\bLd}^2\d t+\frac1\epsilon\int_{s-r}^\infty \norm{\mu-\overline\mu}_{\Ld}^2\d t\leq \frac r 4,\quad \forall s>s_0,\label{basic1}
\end{align}
where we recall that we fixed $r=\frac{T_M}4$. This entails of course a lower bound on the measure of good times:
\begin{align*}
 \normmm{\widehat B_r(s)\setminus C_s}\geq \frac 74 r,\quad \forall s>s_0,
\end{align*}
so that, for any $s>s_0$, we can find a Lebesgue point $t_s$, with  $s-r<t_{s}<s-\frac r2$, such that
\begin{align}
\norm{D\bv(t_{s})}_{\bLd}^2+\norm{\mu(t_{s})-\overline{\mu(t_{s})}}_{\Ld}^2< \epsilon.\label{trajectory}
\end{align}
Note also that, by \eqref{energyineq}, we have $E(\bv(t_{s}),\varphi(t_{s}))\leq E(\bv_0,\varphi_0)=C_0$ for any $s>s_0$.
% we surely have
% \begin{align}
%     \liminf_{t\to\infty} \norm{D\bv(t)}_{\bLd}=0,\quad \liminf_{t\to\infty}\norm{\mu(t)-\overline\mu(t)}_{\bLd}=0.\label{liminf}
% \end{align}

 We can thus consider the initial datum $(\tilde\bv_{0,s},\tilde\varphi_{0,s}):=(\bv(t_{s}),\varphi(t_{s}))$. Thanks to \eqref{trajectory}, \eqref{initial} is satisfied. Hence, we infer that, for any $s>s_0$, there exists a unique strong solution $(\tilde\bv_s,\tilde\varphi_s)$, with $(\tilde\bv_{0,s},\tilde\varphi_{0,s})$ as initial datum, defined in $(0,T_M)$ and such that (see \eqref{holdercontrol})
  \begin{align}
&\label{hol1}\sup_{s>s_0}\norm{\tilde\varphi_s}_{H^1(0,T_M;H^1(\Omega))\cap L^\infty(0,T_M;H^2(\Omega))}\leq C_M,\\&
\sup_{s>s_0}\norm{\tilde\mu_s}_{L^\infty(0,T_M;L^2(\Omega))\cap L^2(0,T_M;H^1(\Omega))}\leq C_M,\\&
\sup_{s>s_0} \norm{\tilde\bv_s}_{H^1(0,T_M;\bLds)\cap L^2(0,T_M;\mathbf H^2(\Omega)
)}\leq C_M,
\
                    \label{holdercontrol2}
                \end{align}
                and the right-hand side only depends on $\epsilon$ and $C_0$. Notice now that, by the fundamental asymptotic strict separation property of weak solutions \eqref{asympt}, it holds, up to enlarging $s_0$,
\begin{align*}
    F'(\varphi)\in L^4(s_0,\infty;L^4(\Omega)),
\end{align*}
as well as there exists $\delta\in(0,1)$ such that
\begin{align}
\norm{\tilde \varphi_{0,s}}_{L^\infty(\Omega)}=\norm{\varphi(t_{s})}_{L^\infty(\Omega)}\leq 1-\delta, \quad \forall s>s_0,
\end{align}
 so that, by property \eqref{sepa2} of strong solutions, recalling the uniform estimates \eqref{holdercontrol}-\eqref{hol1}, there exists $\delta_1\in(0,\delta)$ such that it holds
 \begin{align}
\sup_{t\in(0,T_M)}\norm{\tilde\varphi_s(t)}_{L^\infty(\Omega)}\leq 1-\delta_1, \quad \forall s>s_0.
\end{align}
As a consequence, by the conditional weak-strong uniqueness principle of Theorem \ref{THM:WSU}, the weak solution coincides with each of the strong solutions, namely it holds
\begin{align}
\tilde\bv_s(t)\equiv\bv(t+t_{s}),\quad \tilde\varphi_s(t)\equiv\varphi(t+t_{s}),\quad\text{for a.a. } t\in(0,T_M),\quad \forall s>s_0. \label{coincidence}
\end{align}
Therefore, thanks to the uniform estimates \eqref{hol1}-\eqref{holdercontrol2}, recalling that we chose $t_s$ in such a way that $ \normmm{(t_s,s)}\in[\frac r2,r)=[\frac{T_M}{8},\frac{T_M}4)$, we can choose the sequence $s_n=s_0+n\frac r 2$, $n\in \N$, and thus see that in $[s_0+\frac r2,\infty)$ the solution $(\bv,\vphi)$ becomes strong. This means that there exists  $T_R>s_0$ such that 
  \begin{align*}
&\norm{\varphi}_{H^1_{uloc}([T_R,\infty),H^1(\Omega))\cap L^\infty(T_R,\infty;H^2(\Omega))}\leq \tilde{C}_M,\\&
\norm{\mu}_{L^\infty(T_R,\infty;L^2(\Omega))\cap L^2_{uloc}([T_R,\infty);H^1_{uloc}(\Omega))}\leq \tilde{C}_M,\\&
 \norm{\bv}_{H^1_{uloc}([T_R,\infty);\bLds)\cap L^\infty(T_R,\infty;\bH^1(\Omega))\cap L^2_{uloc}([T_R,\infty);\mathbf H^2(\Omega)
)}\leq \tilde{C}_M,
                \end{align*}
for some $\tilde C_M>0$, entailing that the solution $(\bv,\vphi)$ is a strong solution for $t\geq T_R$. This concludes the proof, since \eqref{equil1} directly comes from \eqref{equil} and by precompactenss, after noticing that $\bv\in BC([T_R,\infty);\bH^1(\Omega))$ and $\vphi\in BC_w([T_R,\infty);H^2(\Omega))$. The proof is finished.
\bigskip
% Since by construction $t_n-t_{1,n}<r=\frac{T_M}4$, we have $t_{n}\in(t_{1,n},t_{1,n}+T_M)$, and thus by \eqref{holdercontrol2} and \eqref{coincidence} we can conclude that
% \begin{align}
% \norm{\bv(t_n)}_{\mathbf H^1_\sigma(\Omega)}+\norm{\vphi(t_n)}_{H^\frac32(\Omega)}\leq C_EC_M,\quad \forall n>N,
%     \label{ff}
% \end{align}
% from which, by compactness, we can extract a further subsequence such that, recalling $\varphi(t_n)\rightharpoonup \vphi_\infty$ weakly in $H^1(\Omega)$, and $\bv(t_n)\rightharpoonup \mathbf 0$ weakly in $\bLds$,
% \begin{align}
% &\vphi(t_n)\to \vphi_\infty,\quad  \text{ strongly in }H^1(\Omega),\nonumber\\&
% \bv(t_n)\to \mathbf 0\quad  \text{ strongly in }\bLds,
%     \label{cmpct}
% \end{align}
% as $n\to\infty$. In other words we have shown that the trajectories are precompact in $\bLds\times H^1(\Omega)$. As a consequence, we infer that the $\omega$-limit set can be expressed as in point (iii) of Lemma \ref{convaaa}, and also that \eqref{convergence1b} holds, which can be shown easily by a contradiction argument exploiting the precompactness of trajectories (see, e.g., \cite{GP}).

\noindent
\textbf{Acknowledgments.}
This research was funded in part by the Austrian Science Fund (FWF) \href{https://doi.org/10.55776/ESP552}{10.55776/ESP552}.
MG and AP are members of Gruppo Nazionale per l'Ana\-li\-si Matematica, la Probabilit\`{a} e le loro Applicazioni (GNAMPA), Istituto Nazionale di Alta Matematica (INdAM). MG's research is part of the activities of ``Dipartimento di Eccellenza 2023--2027'' of Politecnico di Milano.
CH was partially supported by the RTG 2339 ``Interfaces, Complex Structures, and Singular Limits'' of the Deutsche Forschungsgemeinschaft (DFG, German Research Foundation).
%\bibliographystyle{siam}
%\bibliography{NSCAC}

\end{document}